\input ./style/arxiv-general.cfg
\documentclass[aop,MSNbibl,seceqn,dvips]{arximspdf}
\makeatletter
   \@ifpackageloaded{graphicx}{}{\usepackage{graphicx}}
\makeatother
\usepackage{stfloats}
\doi{10.1214/14-AOP948} 
\volume{43}
\issue{6}
\pubyear{2015}
\firstpage{2859}
\lastpage{2948}
\docsubty{FLA}

\makeatletter
\newcommand{\ex}{\operatorname{ex}}
\newcommand{\rrvert}{\vert}
\newcommand{\rrVert}{\Vert}
\newcommand{\llvert}{\vert}
\newcommand{\llVert}{\Vert}
\renewcommand{\mid}{|}
\newtheorem{theorem}{Theorem}[section]
\newtheorem{lemma}{Lemma}[section]
\newtheorem{proposition}{Proposition}[section]
\newtheorem{corollary}{Corollary}[section]
\newproclaim{remark}{Remark}[section]
\def\erm{\mathrm{e}}
\newcommand{\dz}{d}
\newcommand{\dt}{\mathrm{d}}
\newcommand{\dist}{\mathrm{d}}
\newcommand{\diam}{\operatorname{diam}}
\newcommand{\dL}{\mathrm{d}_L}
\newcommand{\dLk}{\mathrm{d}_{L+r}}
\newcommand{\dtL}{\tilde{\mathrm{d}}}
\newcommand{\sh}{\operatorname{Sh}}
\newcommand{\ctwo}{\mathbf{Cond}}
\newcommand{\good}{\operatorname{Good}_L}
\newcommand{\onebad}{\operatorname{OneBad}_L}
\newcommand{\manybad}{\operatorname{ManyBad}_L}
\newcommand{\bbad}{\operatorname{BdBad}_{L,r}}
\newcommand{\badP}{\mathcal{B}}
\newcommand{\taubm}{\tilde{\tau}}
\newcommand{\Ph}{\hat{\Pi}}
\newcommand{\Pho}{\breve{\Pi}}
\newcommand{\Pt}{\tilde{\Pi}}
\newcommand{\Gh}{\hat{G}}
\newcommand{\Gho}{\breve{G}}
\newcommand{\Gt}{\tilde{G}}
\newcommand{\Prw}{\mathrm{P}}
\newcommand{\Erw}{\mathrm{E}}
\newcommand{\pP}{\mathbb{P}}
\newcommand{\pE}{\mathbb{E}}
\newcommand{\pF}{\mathcal{F}}
\newcommand{\mP}{\mathcal{P}}
\newcommand{\bj}{\mathbf{j}}
\newcommand{\mJ}{\mathcal{J}_m}
\makeatother

\begin{document}
\begin{frontmatter}

\title{Exit laws from large balls of (an)isotropic random~walks in
random environment\thanksref{T1}}
\runtitle{Exit laws of (an)isotropic RWRE}

\begin{aug}
\author[A]{\fnms{Erich}~\snm{Baur}\corref{}\ead[label=e1]{erich.baur@math.uzh.ch}}
\and
\author[A]{\fnms{Erwin}~\snm{Bolthausen}\ead[label=e2]{eb@math.uzh.ch}}
\runauthor{E. Baur and E. Bolthausen}
\affiliation{Universit\"at Z\"urich}
\address[A]{Institut f\"ur Mathematik\\
Universit\"at Z\"urich\\
Winterthurerstrasse 190\\
CH-8057 Z\"urich\\
Switzerland\\
\printead{e1}\\
\phantom{E-mail: }\printead*{e2}}
\end{aug}
\thankstext{T1}{Supported in part by the Swiss National Science
Foundation contract nos. PDFMP2\_132099 and 200020\_138141.}

\received{\smonth{9} \syear{2013}}
\revised{\smonth{5} \syear{2014}}

%
\begin{abstract}
We study exit laws from large balls in $\mathbb{Z}^d$, $d\geq3$, of
random walks in an i.i.d. random environment that is a small perturbation
of the environment corresponding to simple random walk. Under a
centering condition on the measure governing the environment, we prove
that the exit laws are close to those of a symmetric random walk, which
we identify as a perturbed simple random walk.
We obtain bounds on total variation distances as well as local results comparing
exit probabilities on boundary segments.
As an application, we prove transience of the
random walks in random environment.

Our work includes the results on isotropic random walks in random
environment of Bolthausen and Zeitouni
[\textit{Probab. Theory Related Fields} \textbf{138} (2007) 581--645]. Since several proofs
in Bolthausen and Zeitouni (2007) were incomplete, a somewhat different
approach was given in
the first author's thesis
[Long-time behavior of random walks in random environment (2013) Z\"urich Univ.]. Here, we extend this approach to
certain anisotropic walks and provide a further step towards a fully
perturbative theory of random walks in random environment.
\end{abstract}

%
\begin{keyword}[class=AMS]
\kwd[Primary ]{60K37}
\kwd[; secondary ]{82C41}
\end{keyword}
\begin{keyword}
\kwd{Random walk}
\kwd{random environment}
\kwd{exit measure}
\kwd{perturbative regime}
\kwd{nonballistic behavior}
\end{keyword}
\end{frontmatter}

\tableofcontents[level=2]

\setcounter{footnote}{1}

\section{Introduction and main results}\label{intro}
\subsection{The model and main results}
\subsubsection{Our model of random walks in random environment}
Consider the integer lattice $\mathbb{Z}^d$ with unit vectors $e_i$, whose
$i$th component equals $1$. We let $\mathcal{P}$ be the set of probability
distributions on $\{\pm e_i\dvtx  i = 1,\ldots,d\}$. Given a probability
measure $\mu$ on $\mathcal{P}$, we\vspace*{1pt} equip
$\Omega=\mathcal{P}^{\mathbb{Z}^d}$ with its natural product $\sigma$-field
$\pF$ and the product measure $\pP= \mu^{\otimes\mathbb{Z}^d}$. Each
element $\omega\in\Omega$ yields transition probabilities of a nearest
neighbor Markov chain on $\mathbb{Z}^d$, the \textit{random walk in random
environment} (RWRE for short), via
\[
p_\omega(x,x+e) = \omega_x(e),\qquad e\in\{\pm
e_i\dvtx  i = 1,\ldots,d\}. %
\]
We write $\Prw_{x,\omega}$ for the ``quenched'' law of the
canonical Markov chain $(X_n)_{n\geq0}$ with these transition
probabilities, starting at $x\in\mathbb{Z}^d$.

We study asymptotic properties of the RWRE in dimension $d\geq3$ when
the underlying environments are small perturbations of the fixed
environment $\omega_x(\pm e_i) =1/(2d)$ corresponding to simple random walk.
\begin{itemize}
\item Let $0<\varepsilon<1/(2d)$. We say that $\mathbf{A0}(\varepsilon)$ holds if
$\mu(\mathcal{P}_\varepsilon) = 1$, where
\[
\mathcal{P}_\varepsilon= \bigl\{q\in\mathcal{P}\dvtx  \bigl\llvert q(\pm
e_i) - 1/(2d)\bigr\rrvert\leq\varepsilon\mbox{ for all }i=1,\ldots,d
\bigr\}. %
\]
\end{itemize}
The perturbative behavior concerns the behavior of the RWRE
when $\mathbf{A0}(\varepsilon)$ holds for small $\varepsilon$.
However, even for arbitrarily small
$\varepsilon$, such walks can behave in very different manners. This motivates
a further ``centering''
restriction on $\mu$.
\begin{itemize}
\item We say that $\mathbf{A1}$ holds if $\mu$ is invariant under reflection
in the coordinate hyperplanes,
that is, under all $d$ reflections
$O_i\dvtx \mathbb{R}^d\rightarrow\mathbb{R}^d$
with
$O_ie_i=-e_i$ and $O_ie_j=e_j$ for $j\neq i$.
\end{itemize}
Condition $\mathbf{A1}$ is weaker than the isotropy condition
introduced by
Bricmont and Kupiainen~\cite{BK}, which requires that $\mu$ is invariant
under all orthogonal transformations
$O\dvtx \mathbb{R}^d\rightarrow\mathbb{R}^d$ fixing the lattice
$\mathbb{Z}^d$. This stronger condition was also assumed in Bolthausen and
Zeitouni~\cite{BZ}, in the first author's thesis~\cite{BA} and in a similar
form in Sznitman and Zeitouni~\cite{SZZ}, who consider isotropic
diffusions. Weaker than $\mathbf{A1}$ is the requirement that $\mu$
is invariant
under $(\omega_0(e))_{\llvert e \rrvert =1}\rightarrow(\omega
_0(-e))_{\llvert e \rrvert =1}$, which
is used in Bolthausen, Sznitman and Zeitouni~\cite{BSZ}, (2.1).

\subsubsection{Our main results}
Write $V_L=\{y\in\mathbb{Z}^d\dvtx  \llvert y \rrvert \leq L\}$
for the discrete ball of radius $L$. Given $\omega\in\Omega$, denote by
$\Pi_L=\Pi_L(\omega)$ the exit distribution from $V_L$ of the random
walk with law $\Prw_{x,\omega}$, that is,
\[
\Pi_L(x,z) = \Prw_{x,\omega} (X_{\tau_L}=z ),
\]
where $\tau_L = \inf\{n\geq0\dvtx  X_n\notin V_L\}$. For probability measures
$\nu_1$ and $\nu_2$, we let $\llVert \nu_1-\nu_2\rrVert _1$ be the
total variation
distance between $\nu_1$ and $\nu_2$. Denote by $\pE$ the
expectation with
respect to $\pP$, and let $p_o(\pm e_i)=p_o(x,x\pm
e_i) = 1/(2d)$ be the transition kernel of simple random walk.
%
\begin{proposition}
\label{main-propkernel}
Assume $\mathbf{A1}$. There is $\varepsilon_0>0$ such that for
$\varepsilon\leq\varepsilon_0$, under
$\mathbf{A0}(\varepsilon)$ the limit
\[
2p_\infty(\pm e_i) = \lim_{L\rightarrow\infty}\sum
_{y\in\mathbb{Z}^d}\pE\bigl[\Pi_L(0,y) \bigr]
\frac{y_i^2}{\llvert y \rrvert ^2} %
\]
exists for $i=1,\ldots,d$. Moreover, $\llVert p_\infty-p_o\rrVert
_1\rightarrow
0$ as $\varepsilon\downarrow0$.
\end{proposition}
From now on, $p_\infty$ is always given by the limit above. The proposition
suggests that for large radii $L$, the RWRE exit measure should be
close to
that of a~symmetric random walk with transition kernel $p_\infty$.
Write $\pi_L^{(p)}(x,\cdot)$ for the exit distribution from $V_L$ of
a random walk
with homogeneous nearest neighbor kernel $p$, started at
$x\in\mathbb{Z}^d$. Recall that we
assume $d\geq3$.
%
\begin{theorem}
\label{main-theorem1}
Assume $\mathbf{A1}$. For $\delta>0$ small enough, there exists
$\varepsilon_0 =
\varepsilon_0(\delta) > 0$ such that if $\mathbf{A0}(\varepsilon)$
is satisfied for some
$\varepsilon\leq
\varepsilon_0$, then
\[
\pP\Bigl(\sup_{x\in V_{L/5}}\bigl\llVert\bigl(\Pi_L-
\pi_L^{(p_\infty
)} \bigr) (x,\cdot)\bigr\rrVert_1>
\delta\Bigr) \leq\exp\bigl(-(\log L)^2 \bigr). %
\]
\end{theorem}
The difference in total variation of the exit laws of the RWRE and the
random walk with kernel $p_\infty$ does not tend to zero as $L\to\infty$, due to localized
perturbations near the boundary. However, with an additional smoothing,
convergence occurs. Let $\rho$ be a random variable that is
independent of
the environment and has a smooth density compactly supported in~$(1,2)$. For $m>0$ and $y\in\mathbb{Z}^d$, write
\[
\Sigma_m^{(p)}(y,\cdot) = E \bigl[\pi_{\rho\cdot m}^{(p)}(0,
\cdot-y) \bigr] %
\]
for the averaged exit distribution from balls $y+V_t$, $t\in(m,2m)$,
of a
random walk with kernel $p$, where $E$ is the
expectation with respect to $\rho$.
%
\begin{theorem}
\label{main-theorem2}
Assume $\mathbf{A1}$. There exists $\varepsilon_0 > 0$ such that if
$\mathbf{A0}(\varepsilon)$ is
satisfied for some $\varepsilon\leq\varepsilon_0$, then for any
$\eta> 0$, we can find
$L_\eta$ and a smoothing radius $m_\eta$ such that for $m\geq m_\eta$,
$L\geq L_\eta$,
\[
\pP\Bigl(\sup_{x\in V_{L/5}}\bigl\llVert\bigl(\Pi_L-
\pi_L^{(p_\infty
)} \bigr)\Sigma_m^{(p_\infty)}(x,
\cdot)\bigr\rrVert_1 > \eta\Bigr) \leq\exp\bigl(-(\log
L)^2 \bigr). %
\]
\end{theorem}
%
%
\begin{remark}\label{remark-after-thm}
(i) As an easy consequence of the last theorem, if one increases the
smoothing scale with $L$, that is, if $m = m_L \uparrow\infty$ (arbitrary
slowly) as $L\rightarrow\infty$, then
\[
\sup_{x\in
V_{L/5}}\bigl\llVert\bigl(\Pi_L-
\pi_L^{(p_\infty)} \bigr)\Sigma_m^{(p_\infty)}(x,
\cdot)\bigr\rrVert_1\rightarrow0\qquad\pP\mbox{-almost surely}.
\]

(ii) The averaging over the radius ensures that the smoothing kernel is
smooth enough: we have, uniformly in $y,y',z$ and for some constant $C$
depending only on the dimension,
\begin{eqnarray*}
\Sigma_m^{(p_\infty)}(y,z) &\leq& Cm^{-d},
\\
\bigl
\llvert\Sigma_m^{(p_\infty)}(y,z)-\Sigma_m^{(p_\infty)}
\bigl(y',z\bigr)\bigr\rrvert
&\leq& C\bigl\llvert y-y'
\bigr\rrvert m^{-(d+1)}\log m; %
\end{eqnarray*}
see Lemma~\ref{app-kernelest} (there, $\hat{\pi}^{(p)}_\psi$ with
$\psi\equiv
m$ takes
the role of $\Sigma_m^{(p)}$). Theorem~\ref{main-theorem2} does still
hold if
$\Sigma_m^{(p_\infty)}$ is replaced by another probability kernel sharing
these properties. However, our particular choice of the smoothing kernel
simplifies the presentation of the proof.
\end{remark}

Our methods enable us to compare the exit measures in a more local
way. Denote by $\partial V_L =\{y\in\mathbb{Z}^d\dvtx  \dist(y, V_L) =1\}
$ the
outer boundary of $V_L$. For positive $t$ and $z\in\partial V_L$ let
$W_t(z)=V_t(z)\cap\partial V_L$, where $V_t(z)=z+V_t$. Then $\llvert
W_t(z)\rrvert $ is
of order $t^{d-1}$. We obtain the following:
%
\begin{theorem}
\label{local-thm-exitmeas}
Assume $\mathbf{A1}$. There exist $\varepsilon_0 > 0$ and $L_0>0$
such that if
$\mathbf{A0}(\varepsilon)$ is satisfied for some $\varepsilon\leq
\varepsilon_0$, then for $L\geq L_0$,
there exists an event $A_L\in\pF$ with $\pP(A^c_L) \leq\exp
(-(1/2)(\log
L)^2)$ such that on $A_L$, the following holds true. If $0<\eta<1$ and
$x\in V_{\eta L}$, then for all $z\in\partial V_L$:
\begin{longlist}[(ii)]
\item For $t\geq L/(\log L)^{15}$, there
exists $C=C(\eta)$ with
\[
\Pi_L\bigl(x,W_t(z)\bigr) \leq C
\pi_L^{(p_o)}\bigl(x,W_t(z)\bigr). %
\]
\item There exists a homogeneous symmetric nearest neighbor kernel
$p_L$ such that for $t\geq L/(\log L)^{6}$,
\[
\Pi_L\bigl(x,W_t(z)\bigr) = \pi_L^{(p_L)}
\bigl(x,W_t(z)\bigr) \bigl(1+O \bigl((\log L)^{-5/2} \bigr)
\bigr). %
\]
Here, the constant in the $O$-notation depends only on $d$ and $\eta$.
\end{longlist}
\end{theorem}
We give one possible choice for the kernel $p_L$ in~(\ref{kernelpL}). Our
results can also be used to deduce transience of the RWRE.

%
\begin{corollary}
\label{main-transience}
Assume $\mathbf{A1}$. There exist $\varepsilon_0$ such that if
$\mathbf{A0}(\varepsilon)$ is
satisfied for some $\varepsilon\leq\varepsilon_0$, then on $\pP
$-almost all
$\omega\in\Omega$ the RWRE $(X_n)_{n\geq0}$ is transient.
\end{corollary}

%
\begin{remark}
Let us mention the simplest nontrivial example of a RWRE under conditions
$\mathbf{A0}(\varepsilon)$ and $\mathbf{A1}$, with no isotropy. Fix
a coordinate
direction, say~$e_1$, and define two symmetric kernels $q,q'\in\mathcal
{P}_\varepsilon$ by setting
\begin{eqnarray*}
q(e) &=& \frac{1}{2d} + \cases{ \varepsilon, &\quad for $e=e_1$,
\cr
-\varepsilon,&\quad for $e=-e_1$,
\cr
0, &\quad for $e\neq\pm
e_1$,}
\\
q'(e)&=& \frac{1}{2d} + \cases{ -
\varepsilon, &\quad for $e=e_1$,
\cr
\varepsilon, &\quad for
$e=-e_1$,
\cr
0, &\quad for $e\neq\pm e_1$.} %
\end{eqnarray*}
Then the law $\mu$ on $\mathcal{P}_\varepsilon$ with $\mu(q) = \mu
(q') = 1/2$
satisfies $\mathbf{A0}(\varepsilon)$ and $\mathbf{A1}$. With this
choice of $\mu$,
Corollary~\ref{main-transience} settles the generalization of Problem~4 in Kalikow~\cite{KAL} to dimensions $d\geq3$ (for small disorder).
\end{remark}

\subsection{Discussion of this work}
This paper is inspired by the work of Bolt\-hausen and
Zeitouni~\cite{BZ}. There, Theorems~\ref{main-theorem1}
and~\ref{main-theorem2} appeared in a similar form for the case of
isotropic RWRE in dimension $d\geq3$. A corrected and extended version
of~\cite{BZ} forms part of the first author's thesis~\cite{BA}. Our work
should in turn be understood as an extension of~\cite{BA} to the case of
certain anisotropic random walks in random environment.

Here, the main difficulty stems from the fact that the kernel
$p_{\infty}$
is not explicitly computable and depends in a complicated way on $\mu
$. We
will \textit{not} first prove the existence of $p_\infty$ and then deduce our
results about the exit measures---in fact, it will be a side effect of our
multiscale analysis of exit laws that $p_\infty$ exists and is the right
object of comparison. The idea of its construction starts with the
observation that if the statements of Theorems~\ref{main-theorem1}
and~\ref{main-theorem2} are true for \textit{some} kernel $p_\infty$, then
the averaged exit distribution on a global scale will be the same as the
exit distribution of the random walk with kernel $p_\infty$, when
$L\rightarrow\infty$. Therefore, it is natural to choose for any
scale $L$
a symmetric transition kernel $p_L$ which has the property that the
covariance matrix of the averaged exit distribution from $V_{L}$, scaled
down by $L^{2}$, is the covariance matrix of $p_L$ [in fact, we will choose
$p_L$ in a slightly different way; see~(\ref{kernelpL}) for the precise
definition]. The difficult task is then to show that
$p_{\infty}=\lim_{L\rightarrow\infty}p_L $ exists. In the isotropic case,
this problem is absent since one can choose $p_L =p_o$ for every $L$.

The thesis~\cite{BA} develops a somewhat new approach to the isotropic case
covered in~\cite{BZ}, which is, as we hope, easier to understand. Since
we follow here the same strategy, let us explain the main
changes compared to~\cite{BZ} and point at some of the problems which
appeared there.

Our focus lies on (coarse grained) Green's function estimates on a large
class of environments, so-called \textit{goodified} environments. These concepts
are developed in Section~\ref{super}. In contrast with~\cite{BZ}, we state
our core Green's function estimates (Lemma~\ref{superlemma}, often used
in the version of Lemma~\ref{superlemma2}) in terms of an appropriate
notion of domination of kernels, and also employ basic operations on
kernels; see, for example, Propositions~\ref{super-concatenating}
and~\ref{super-keyest}. While Lemma~\ref{superlemma} requires some
effort to be set up, it then yields in a relatively straight forward and
systematic way controls on both smoothed and nonsmoothed estimates, see,
for example, Lemma~\ref{lemma-badpart}. In contrast, the goodified
Green's function estimates in~\cite{BZ}, namely, (4.24) and (4.25),
are weaker
than our Lemma~\ref{superlemma}. (We note in passing that fleshing out
the missing details in the proof of~\cite{BZ}, (4.24), without a
version of
Lemma~\ref{superlemma} seems challenging; see, e.g., the end of
Section~4.3, page 606 there.) Likewise, the lack in~\cite{BZ} of a
statement like Lemma~\ref{superlemma} makes the derivation of (4.46)
there incomplete. The same issue arises, for dimensions $d=3,4$, in the
derivation of (4.48) and (4.49) in~\cite{BZ}.

With our Green's function estimates and the concept of goodified
environments, we give proofs of the main results in our
Sections~\ref{smv-exits} and~\ref{nonsmv-exits}, which differ even in the
mere isotropic case in many details from the derivation in~\cite{BZ}. We
believe that our proofs are more transparent. The reader who is primarily
interested in the isotropic case is, however, advised to consult the
thesis~\cite{BA} first.

Finally, our \hyperref[appe]{Appendix} includes the results in~\cite{BZ} on simple random
walk and standard Brownian motion as special cases. We include the proofs
both because our statements are more general, and also because the proofs
of different cases are only sketched or altogether omitted in~\cite{BZ},
for example, in the proof of Lemma 3.4 there; we also provide a lower
bound on exit
probabilities [Lemma~\ref{hittingprob}(iii)] which is implicitly used
in~\cite{BZ}, but not proved there.

For a better reading, a rough overview over this paper is given in
Section~\ref{readingguide}.

\subsection{Some relevant literature}
Let us comment on some further literature which is relevant for our
study. For a detailed survey on RWRE, we refer to the
lecture notes of Sznitman~\cite{SZ-LN,SZ-LN2} and
Zeitouni~\cite{ZT-StF,ZT-TRev}, and also to the overview
article of
Bogachev~\cite{BG}.

Assuming $\mathbf{A0}(\varepsilon)$ for small $\varepsilon$ and the
stronger isotropy condition
that was mentioned at the beginning, Bricmont and Kupiainen~\cite{BK}
prove a (quenched) invariance principle, showing that in dimensions
$d\geq
3$, the RWRE is asymptotically Gaussian, on $\pP$-almost all
environments. A continuous counterpart, isotropic diffusions in a random
environment which are small perturbations of Brownian motion, has been
investigated by Sznitman and Zeitouni in~\cite{SZZ}. They prove transience
and a full quenched invariance principle in dimensions $d\geq3$.

Our centering condition $\mathbf{A1}$ excludes so-called ballistic behavior,
that is, the regime where the limit velocity $v=\lim_{n\rightarrow
\infty}
X_n/n$ is an almost sure constant vector different from zero. Ballistic
behavior has been studied extensively, for example, by Kalikow~\cite{KAL},
Sznitman~\cite{SZ1,SZ2,SZ4}, Bolthausen and
Sznitman~\cite{BOLTSZ}, or more recently by Berger~\cite{BER} and Berger,
Drewitz and Ram\'irez~\cite{BERDRERAM}.

In the perturbative regime when $d\geq3$, Sznitman~\cite{SZ4} shows that
some strength of the mean local drift $m=\pE[\sum_{\llvert e \rrvert
=1}e\omega
_0(e)]$ is
enough to deduce ballisticity. However, as examples in Bolthausen,
Sznitman and Zeitouni~\cite{BSZ} for dimensions $d\geq7$ demonstrate,
ballisticity can also occur with $m=0$, and one can even construct
examples exhibiting ballistic behavior with $v=-cm$ and $c>0$.
Note that our condition $\mathbf{A1}$ implies $m=v=0$.

The work of Bolthausen, Sznitman and Zeitouni~\cite{BSZ} provides also
examples and results for nonballistic behavior. They consider the special
class of multidimensional RWRE for which the projection onto at least $d_1
\geq5$ components behaves as a standard random walk. In particular, if
$d_1 \geq7$ and the law of the environment is invariant under the
antipodal transformation (\cite{BSZ}, (2.1)), mentionned at the
beginning, a
quenched invariance principle is proved.

Much is also known for the class of \textit{balanced} RWRE when
$\pP(\omega_0(e_i) = \omega_0(-e_i))=1$ for all $i=1,\ldots,d$. Employing
the method of environment viewed from the particle, Lawler proves
in~\cite{LAWbalanced} that for $\pP$-almost all $\omega$,
$X_{\lfloor
n\cdot\rfloor}/\sqrt{n}$ converges in $\Prw_{0,\omega
}$-distribution to
a nondegenerate Brownian motion with diagonal covariance matrix, even in
the nonperturbative regime. Moreover, the RWRE is recurrent in dimension
$d=2$ and transient when $d\geq3$; see~\cite{ZT-StF}. Recently, within
the i.i.d. setting, diffusive behavior has been shown in the mere elliptic
case by Guo and Zeitouni~\cite{GUOZT} and in the nonelliptic case by
Berger and Deuschel~\cite{BERDEU}.

\section{Basic notation and main techniques}\label{prel}
\subsection{Basic notation}
Our purpose here is to cover the most relevant notation which will be used
throughout this text. Further notation will be introduced later on when needed.
\subsubsection{Sets and distances}
We let $\mathbb{N} = \{0,1,2,3,\ldots\}$ and $\mathbb{R}_+ =
\{x\in\mathbb{R}\dvtx  x \geq0\}$. For a set $A$, its complement is
denoted by
$A^c$. If $A\subset\mathbb{R}^d$ is measurable and nondiscrete, we write
$\llvert A \rrvert $ for its $d$-dimensional Lebesgue measure.
Sometimes, $\llvert A \rrvert $ denotes
the surface measure instead, but this will be clear from the context. If
$A\subset\mathbb{Z}^d$, then $\llvert A \rrvert $ denotes its cardinality.

For $x\in\mathbb{R}^d$, $\llvert x \rrvert $ is the Euclidean norm.
If $A, B \subset
\mathbb{R}^d$, we set $\dist(A,B) = \inf\{\llvert x-y\rrvert\dvtx  x\in A,
y\in B\}
$ and
$\diam(A) = \sup\{\llvert x-y\rrvert\dvtx  x,y \in A\}$. Given $L > 0$, let
$V_L=\{x\in\mathbb{Z}^d\dvtx  \llvert x \rrvert \leq L\}$ and for $x\in
\mathbb{Z}^d$,
$V_L(x) = x + V_L$. For Euclidean balls in $\mathbb{R}^d$ we write
$C_L=\{x\in\mathbb{R}^d\dvtx  \llvert x \rrvert < L\}$ and for $x\in\mathbb{R}^d$,
$C_L(x) = x + C_L$.

If $V\subset\mathbb{Z}^d$, then $\partial V =\{x\in V^c\cap\mathbb
{Z}^d\dvtx  \dist(\{x\},V) = 1\}$ is the outer boundary, while in the case
of a
nondiscrete set $V\subset\mathbb{R}^d$, $\partial V$ stands for the usual
topological boundary of $V$ and $\overline{V}$ for its closure. For
$x\in\overline{C}_L$, we set $\dL(x) = L-\llvert x \rrvert $.
Finally, for $0\leq
a<b\leq L$, the ``shell'' is defined by
\[
\sh_L(a,b)=\bigl\{x\in V_L\dvtx  a\leq\dL(x) < b\bigr\},\qquad
\sh_L(b) = \sh_L(0,b). %
\]

\subsubsection{Functions}
If $a,b$ are two real numbers, we set $a\wedge b = \min\{a,b\}$,
$a\vee
b=\max\{a,b\}$. The largest integer not greater than $a$ is denoted by
$\lfloor a\rfloor$. As usual, set $1/0 =\infty$. For us, $\log$ is the
logarithm to the base $\erm$, and $\log_a$ is then the logarithm to
the base
$a$. For $x,z\in\mathbb{R}^d$, the Delta function $\delta_{x}(z)$ is
defined to be equal to one for $z=x$ and zero otherwise.

Given two functions $F,G\dvtx
\mathbb{Z}^d\times\mathbb{Z}^d\rightarrow\mathbb{R}$, we write
$FG$ for
the (matrix) product $FG(x,y) = \sum_{u\in\mathbb{Z}^d}F(x,u)G(u,y)$,
provided the right-hand side is absolutely summable. $F^k$ is the $k$th
power defined in this way, and $F^0(x,y) = \delta_x(y)$. $F$ can also
operate on functions $f\dvtx \mathbb{Z}^d\rightarrow\mathbb{R}$ from the
left via $Ff(x) = \sum_{y\in\mathbb{Z}^d}F(x,y)f(y)$.\vspace*{1pt}

We use the symbol $1_W$ for the indicator
function of the set $W$. By an abuse of notation, $1_W$ will also
denote the kernel
$(x,y)\mapsto1_W(x)\delta_x(y)$. If $f\dvtx \mathbb{Z}^d\rightarrow
\mathbb{R}$, $\llVert f\rrVert _1 = \sum_{x\in\mathbb{Z}^d}\llvert
f(x) \rrvert \in[0,\infty
]$ is
its $L^1$-norm. When $\nu\dvtx  \mathbb{Z}^d\rightarrow\mathbb{R}$ is a
(signed) measure, $\llVert \nu\rrVert _1$ is its total variation norm.

Let $U\subset\mathbb{R}^d$ be a bounded open set, and let
$k\in\mathbb{N}$. For a real-valued function $f$ with $f\mid_U\in
C^k(U)$, that is, $f$ is $k$-times continuously differentiable in $U$, we
define for $i=0,1,\ldots,k$,
\[
{\bigl\llVert D^if\bigr\rrVert}_{U} = \sup
_{\llvert \beta\rrvert =
i}\sup_{U}\biggl\llvert
\frac{\partial^{i}}{\partial x_1^{\beta_1}\cdots
\partial x_d^{\beta_d}}f\biggr\rrvert, %
\]
where the first supremum is over all multi-indices $\beta=
(\beta_1,\ldots,\beta_d)$, $\beta_j\in\mathbb{N}$, with
$\llvert \beta\rrvert =\sum_{j=1}^d\beta_j$. We also write $\nabla$
for the
gradient of a
function, and in some proofs, $\Delta$ denotes the Laplace operator.

Let $L>0$, and put $\mathcal{U}_L = \{x\in\mathbb{R}^d\dvtx  L/2<
\llvert x \rrvert <
2L \}$. We denote by $\mathcal{M}_L$ the set of functions $\psi$,
whose restrictions to $\mathcal{U}_L$ satisfy the following properties:
\begin{itemize}
\item$\psi\mid_{\mathcal{U}_L}\dvtx  \mathcal{U}_L\rightarrow(L/10,5L)$,
\item$\psi\mid_{\mathcal{U}_L}\in C^4(\mathcal{U}_L)$ with
${\llVert D^i\psi\mid_{\mathcal{U}_L}\rrVert }_{\mathcal{U}_L}\leq10$ for
$i=1,2,3,4$.
\end{itemize}
Functions in $\mathcal{M}_L$ will be used to define smoothing kernels with
good smoothing properties.

\subsubsection{Transition probabilities and exit distributions}
Given (not necessarily nearest neighbor) transition probabilities $p =
{(p(x,y))}_{x,y\in\mathbb{Z}^d}$, we write $\Prw_{x,p}$ for the law
of the
canonical Markov chain ${(X_n)}_{n\geq0}$ on
$({(\mathbb{Z}^d)}^{\mathbb{N}},{\mathcal{G}})$,
${\mathcal{G}}$ the $\sigma$-algebra generated
by cylinder functions, with transition probabilities $p$ and starting point
$X_0=x$ $\Prw_{x,p}$-a.s. The expectation with respect to $\Prw_{x,p}$
is denoted by $\Erw_{x,p}$. The simple random walk kernel $p_o(x,x\pm e_i)
= 1/(2d)$ will play a prominent role. Clearly, every $p\in\mathcal{P}$
gives rise to a homogeneous nearest neighbor kernel, which by a small
abuse of notation we again denote by $p$.

If $V\subset\mathbb{Z}^d$, we denote by $\tau_V = \inf\{n\geq0\dvtx
X_n\notin V\}$ the first exit time from $V$, with $\inf\varnothing=
\infty$, whereas $T_V = \tau_{V^c}$ is the first hitting time of $V$. Given
$x,z\in\mathbb{Z}^d$ and $p, V$ as above, we define
\[
\ex _{V}(x,z; p) = \Prw_{x,p} (X_{\tau_V}=z ).
\]
Notice that for $x\in V^c$, $\ex _{V}(x,z;p) = \delta_{x}(z)$.

For $p\in\mathcal{P}$, we write
\[
\pi^{(p)}_V(x,z) = \ex _{V} (x,z;p ),
\]
and for $\omega\in\Omega$, we set
\[
\Pi_{V,\omega}(x,z) = \ex _{V}(x,z;p_{\omega}).
\]
We usually suppress $\omega$ from the notation and simply write $\Pi_V$. Mostly, we shall
interpret $\Pi_V$ as a \textit{random} exit distribution. However, sometimes we work with
a fixed environment $\omega\in \Omega$, and then we still write
$\Pi_V$ instead of $\Pi_{V,\omega}$.

Recall the definitions of the sets $\mP$ and $\mP_\varepsilon$ from the
\hyperref[intro]{Introduction}. For $0<\kappa<1/(2d)$, let
\[
\mathcal{P}^{\mathrm{s}}_\kappa=\bigl\{p\in\mP_\kappa\dvtx
p(e_i)=p(-e_i), i=1,\ldots,d\bigr\}, %
\]
that is, $\mathcal{P}^{\mathrm{s}}_\kappa$ is the subset of $\mP
_{\kappa}$ which contains all
symmetric
probability distributions on $\{\pm e_i\dvtx  i = 1,\ldots,d\}$. At various
places, the parameter
$\kappa$ bounds the range of the symmetric transition kernels we work
with.

\subsubsection{Coarse grained transition kernels}
\label{prel-cgrw}
Fix once for all a probability density $\varphi\in
C^\infty(\mathbb{R}_+,\mathbb{R}_+)$ with compact support in
$(1,2)$. Given\vspace*{1pt}
a transition kernel $p\in\mathcal{P}$ and a strictly positive
function $\psi= (m_x)_{x\in W}$, where
$W\subset\mathbb{R}^d$ with $W\cap\mathbb{Z}^d\neq\varnothing$, we define
the coarse grained transition kernels on $W\cap\mathbb{Z}^d$
associated to
$(\psi, p)$,
%
%
\begin{equation}
\label{prel-cgpsi} \hat{\pi}^{(p)}_{\psi}(x,\cdot) =
\frac{1}{m_x}\int_{\mathbb{R}_+} \varphi\biggl(\frac{t}{m_x}
\biggr)\ex _{V_t(x)}(x,\cdot;p)\,\dt t,\qquad x\in W\cap
\mathbb{Z}^d.
\end{equation}
Mostly, we will take $\psi\in\mathcal{M}_L$, and then~(\ref{prel-cgpsi})
yields a collection of transition kernels on at least
$\mathcal{U}_L\cap\mathbb{Z}^d$. Often, we consider for $m>0$ the constant
function $\psi\equiv m$ (sometimes denoted $\psi_m$), and
then~(\ref{prel-cgpsi}) gives coarse grained transition kernels on the
whole grid $\mathbb{Z}^d$.

\subsubsection{Coarse graining schemes in the ball}
Similarly to~(\ref{prel-cgpsi}), we will now introduce coarse grained
transition kernels for the motion inside the ball $V_L$, for both symmetric
random walk and RWRE.

We use a particular function $\psi$. Once for all, let
\[
s_L = \frac{L}{(\log L)^3} \quad\mbox{and}\quad r_L=
\frac{L}{(\log
L)^{15}}. %
\]
Our coarse graining schemes in the ball are indexed by a parameter $r$,
which can either be a constant $\geq$100, but much smaller than $r_L$, or,
in most of the cases, $r=r_L$. We fix a smooth function $h\dvtx
\mathbb{R}_+\rightarrow\mathbb{R}_+$ satisfying
\[
h(x)= \cases{ x, &\quad for $x\leq1/2$,
\cr
1, &\quad for $x\geq2$,} %
\]
such that $h$ is concave and increasing on $(1/2,2)$.
Define $h_{L,r}\dvtx  \overline{C}_L \rightarrow\mathbb{R}_+$ by
%
%
\begin{equation}
\label{smoothbad-hlr} h_{L,r}(x) = \frac{1}{20}\max\biggl
\{s_L h \biggl(\frac{\dL(x)}{s_L} \biggr), r \biggr\}.
\end{equation}
Since we mostly work with $r=r_L$, we use the abbreviation $h_L = h_{L,r_L}$.
We write $\Ph_{L,r}(=\Ph_{L,r,\omega})$ for the coarse grained RWRE transition kernel inside
$V_L$ associated to $(\psi= (h_{L,r}(x) )_{x\in
V_L},p_\omega)$,
\[
\Ph_{L,r}(x,\cdot) = \frac{1}{h_{L,r}(x)}\int_{\mathbb{R}_+}
\varphi\biggl(\frac{t}{h_{L,r}(x)} \biggr)\Pi_{V_t(x)\cap V_L}(x,\cdot)\,\dt t,
\]
and $\hat{\pi}^{(p)}_{L,r}$ for the coarse grained kernel coming
from symmetric
random walk with transition
kernel $p\in\mathcal{P}$, where in the definition $\Pi$ is replaced
by $\pi^{(p)}$.
Most of the time we view $\Ph_{L,r}$ as a random transition kernel, but we shall also
write $\Ph_{L,r}$ if the underlying environment $\omega$ is fixed.
For convenience, we set $\Ph_{L,r}(x,\cdot)
=\hat{\pi}^{(p)}_{L,r}(x,\cdot) = \delta_x(\cdot)$ for $x\in
\mathbb
{Z}^d\setminus
V_L$. By the strong Markov property, the exit measures from the
ball $V_L$ remain unchanged under these transition kernels, that is,
%
%
\begin{equation}
\label{prel-exitunchanged} \ex _{V_L} (x,\cdot;\Ph_{L,r} ) =
\Pi_L(x,\cdot) \quad\mbox{and}\quad\ex _{V_L}
\bigl(x,\cdot;\hat{\pi}^{(p)}_{L,r} \bigr) = \pi^{(p)}_L(x,\cdot).
\end{equation}
See Figure \ref{fig1} for a visualization of the coarse graining
scheme.

%
%
\begin{remark}
(i) Later on, we will also work with slightly modified transition kernels
$\Pho$ and ${\breve{\pi}^{(p)}}$, which depend on the environment.
We elaborate on
this
in Section~\ref{super-cgmod}.

%
%
\begin{figure}

\includegraphics{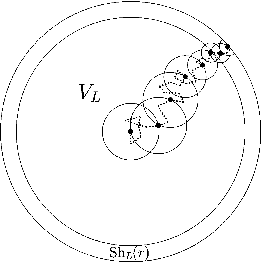}

\caption{The coarse graining scheme in $V_L$. In the bulk $\{x\in
V_L\dvtx  \dL (x) \geq2s_L\}$, the exit distributions are taken from balls
of radii between $(1/20)s_L$ and $(1/10)s_L$. When entering
$\sh_L(2s_L)$, the coarse graining radii start to shrink, up to the
boundary layer $\sh_L(r)$, where the exit distributions are taken from
intersected balls $V_t(x)\cap V_L$, $t\in[(1/20)r, (1/10)r]$. The
dotted lines indicate a corresponding random walk sample path.}\label{fig1}
\end{figure}

(ii) Due to the lack of the last smoothing step outside $V_L$, we need to
zoom in near the boundary in order to handle nonsmoothed exit
distributions in Section~\ref{nonsmv-exits}. The parameter $r$ allows
us to adjust the step size in the boundary region.

(iii) For every choice of $r$,
\[
h_{L,r}(x) = \cases{ \dL(x)/20, &\quad for $x \in V_L$
with $r_L\leq\dist_L(x) \leq s_L/2$,
\vspace*{2pt}\cr
s_L/20, &\quad for $x \in V_L$ with $
\dist_L(x) \geq2s_L$.} %
\]
\end{remark}

\subsubsection{Abbreviations}
If it is clear from the context which transition
kernel~$p$ we are working with, we often drop\vspace*{1pt} the sub- or superscript
$p$ from
notation. Then, for example, we write $\pi_V$ for $\pi^{(p)}_V$,
$\Prw_x$
instead of $\Prw_{x,p}$ or $\Erw_x$ for $\Erw_{x,p}$. Given transition
probabilities $p_\omega$ coming from an environment $\omega$, we use
the notation $\Prw_{x,\omega}$,~$\Erw_{x,\omega}$.

If $V=V_L$ is the ball around zero of radius $L$, we usually write $\pi_L$
instead of $\pi_V$, $\Pi_L$ for $\Pi_V$ and $\tau_L$ for $\tau_{V}$.

Many of our quantities, for example, the transition kernels $\Ph_{L,r}$,
$\hat{\pi}_{L,r}$ or the kernel $\Gamma_{L,r}$ which
is introduced in
Section~\ref{super}, are indexed by both $L$ and $r$. While we always keep
the indices in the statements, we normally drop both of them in the proofs.

Finally, we will often use the
abbreviations $\dist(y,B)$ for $\dist(\{y\},B)$, $T_x$ for $T_{\{x\}
}$ and
$\pP(A; B )$ for $\pP(A\cap B )$.

\subsubsection{Some words about constants, $O$-notation and large $L$ behavior}
All our constants are positive. They only depend on the dimension
$d\geq3$
unless stated otherwise.
In particular, constants do \textit{not} depend on $L$, on $\delta$, on
$\omega$ or on any point $x\in\mathbb{Z}^d$, and they are also independent
of the parameter $r$.

At some places, one might have the impression that constants depend
on the transition kernel $p$. However, we only work with
$p\in\mathcal{P}^{\mathrm{s}}_\kappa$, and $\kappa$ can be
chosen (arbitrarily) small. Such kernels
$p$ are therefore small perturbations of the simple
random walk kernel $p_o$, and since all dependencies emerge in a
continuous way (in $p$), we may always assume that constants are
uniform in $p$.

We use $C$ and $c$ for generic positive constants
whose values can change in different expressions, even in the same
line. In
the proofs, we often use other constants like $K, C_1, c_1$; their values
are fixed throughout the proofs. Lower-case constants usually indicate
small (positive) values.

Given two functions $f,g$ defined on some subset of $\mathbb{R}$, we write
$f(t) = O (g(t) )$ if there exists a positive $C > 0$ and a real
number $t_0$ such that $\llvert f(t)\rrvert \leq C \llvert
g(t)\rrvert$ for $t\geq t_0$.

If a statement holds for ``$L$ large (enough),'' this means that there exists
$L_0>0$ depending only on the dimension
such that the statement is true for all $L\geq L_0$. This applies
analogously to expressions like ``$\delta$
(or $\varepsilon$, or $\kappa$) small (enough).''

One should always keep in mind that we are interested in asymptotics
when $L\rightarrow\infty$ and the perturbation parameter $\varepsilon
$ is
arbitrarily small, but fixed. Even though some of our statements are valid
only for large $L$ and $\varepsilon$ (or $\delta$, or $\kappa$)
sufficiently small,
we do not mention this every time.

\subsection{Perturbation expansion for Green's functions}
Our approach of comparing the RWRE exit distribution with that of an
appropriate symmetric random walk is based on a perturbation argument.
Namely, the resolvent equation allows us to express Green's functions
of the RWRE in
terms of Green's functions of homogeneous random walks. More generally,
let $p =
(p(x,y) )_{x,y\in\mathbb{Z}^d}$ be a family of finite range
transition probabilities on $\mathbb{Z}^d$, and let $V\subset\mathbb{Z}^d$
be a finite set. The corresponding Green's kernel or Green's function for
$V$ is defined by
\[
g_V(p) (x,y) = \sum_{k=0}^\infty
(1_Vp )^k(x,y). %
\]
The connection with the exit measure is given by the fact that for
$z\notin
V$, we have
%
%
\begin{equation}
\label{prel-pbe0} g_V(p) (\cdot,z) = \ex _{V}(
\cdot,z;p).
\end{equation}
Now write $g$ for $g_V(p)$, and let $P$ be another transition
kernel with corresponding Green's function $G$ for $V$. With $\Delta=
1_V (P-p )$, we have by the resolvent equation
%
%
\begin{equation}
\label{prel-pbe1} G -g = g\Delta G = G\Delta g.
\end{equation}
In order to get rid of $G$ on the right-hand side, we
iterate~(\ref{prel-pbe1}) and obtain
%
%
\begin{equation}
\label{prel-pbe2} G -g = \sum_{k=1}^\infty(g
\Delta)^kg,
\end{equation}
provided the infinite series converges, which is always the case in our
setting. Writing~(\ref{prel-pbe2}) as
\[
G = g\sum_{k=0}^\infty(\Delta g
)^k, %
\]
replacing the rightmost $g$ by $g(x,\cdot)= \delta_x(\cdot) +
1_Vpg(x,\cdot)$
and reordering terms, we get
%
%
\begin{equation}
\label{prel-pbe3} G = g\sum_{m=0}^\infty{ (Rg
)}^m\sum_{k=0}^\infty
\Delta^k,
\end{equation}
where $R= \sum_{k=1}^\infty\Delta^kp$.

Two Green's functions for the ball $V_L$ will play a particular role: the
(coarse grained) RWRE Green's function $\Gh_{L,r}$ corresponding to
$\Ph_{L,r}$, and the Green's function $\hat{g}_{L,r}$
corresponding to
$\hat{\pi}_{L,r}$,
\[
\Gh_{L,r}(x,y)=\sum_{k=0}^\infty
(1_{V_L}\Ph_{L,r} )^k(x,y),\qquad
\hat{g}_{L,r}(x,y)=\sum_{k=0}^\infty
(1_{V_L}\hat{\pi}_{L,r} )^k(x,y). %
\]
A ``goodified'' version of $\Gh_{L,r}$ will be introduced in
Section~\ref{smoothbad}.

\subsection{Main technical statement}
We will deduce our main results from Proposition~\ref{main-prop}
below. The latter involves a technical condition, which we will propagate
from one level to the next. This condition depends on the deviation
$\delta$ (cf. Theorem~\ref{main-theorem1}) and on a parameter
$L_0\geq
3$ which will finally be chosen sufficiently large.

Recall the coarse graining schemes on $V_L$. Even though we use the
``final'' kernel $p_\infty$ in the formulation of our main theorems, we
will work in the proofs with intermediate kernels $p_L$ depending on the
radius of the ball.
More precisely, we assign to each $L>0$ the symmetric transition kernel
$(i=1,\ldots,d)$
%
%
\begin{equation}
\label{kernelpL} p_L(\pm e_i) =\cases{ 1/(2d),&\quad for
$0<L\leq L_0$,
\vspace*{2pt}\cr
\displaystyle\frac{1}{2}\sum
_{y\in\mathbb{Z}^d}\pE\bigl[\Ph_{L,r}(0,y) \bigr]
\frac{y_i^2}{\llvert y \rrvert ^2},&\quad for $L>L_0$.}
\end{equation}
Since $h_{L,r}(0)=s_L/20$, the definition of $p_L$ does not
depend on the parameter $r$. In
words, for radii $0<L\leq L_0$, $p_L$ agrees with the simple random walk
kernel $p_o$, while for $L>L_0$ the kernel $p_L$ is defined as an average
of variances of normalized mean exit distributions from balls of radii
$t\in[(1/20)s_L,(1/10)s_L]$.

For $\psi\in\mathcal{M}_t$ and $p,q\in\mathcal{P}$, define
\begin{eqnarray*}
D_{t,p,\psi,q}^{\ast} &=&\sup_{x\in
V_{t/5}}\bigl\llVert
\bigl(\Pi_{V_t}-\pi^{(p)}_{V_t} \bigr){\hat{
\pi}}_{\psi}^{(q)}(x,\cdot)\bigr\rrVert_1,
\\
D_{t,p}^\ast&=& \sup_{x\in V_{t/5}}\bigl\llVert
\bigl(\Pi_{V_t}-\pi^{(p)} _{V_t} \bigr) (x,\cdot)
\bigr\rrVert_1. %
\end{eqnarray*}
With $\delta> 0$, we set for $i=1,2,3$
\begin{eqnarray*}
&& b_i(L,p,\psi,q,\delta)
\\
&&\qquad = \pP\bigl( \bigl\{(\log L)^{-9+9(i-1)/4} < D_{L,p,\psi
,q}^\ast
\leq(\log L)^{-9 +9i/4} \bigr\}\cap\bigl\{D_{L,p}^\ast
\leq\delta\bigr\} \bigr)
\end{eqnarray*}
and
\[
b_4(L,p,\psi,q,\delta)= \pP\bigl( \bigl\{D_{L,p,\psi,q}^\ast>
(\log L)^{-3+3/4} \bigr\}\cup\bigl\{D_{L,p}^\ast>
\delta\bigr\} \bigr). %
\]
Put $\iota= (\log L_0)^{-7}$, and let us now formulate the following:
\subsubsection{Condition $\mathbf{Cond}$}
Let $\delta> 0$ and $L_1\geq L_0\geq3$. We say that\break  $\ctwo(\delta,L_0,L_1)$
holds if:
\begin{itemize}
\item For all $3\leq L\leq2L_0$, all $\psi\in\mathcal{M}_L$ and
all $q\in\mathcal{P}^{\mathrm{s}}_\iota$,
\[
\pP\bigl( \bigl\{D_{L,p_o,\psi,q}^\ast> (\log L)^{-9}
\bigr\}\cup\bigl\{D_{L,p_o}^\ast> \delta\bigr\} \bigr) \leq
\exp\bigl(-\bigl(\log(2L_0)\bigr)^2 \bigr). %
\]
\item For all $L_0< L\leq L_1$, $L'\in[L,2L]$, $\psi\in\mathcal
{M}_{L'}$ and
$q\in\mathcal{P}^{\mathrm{s}}_\iota$,
\[
b_i\bigl(L',p_L,\psi,q,\delta\bigr)\leq
\tfrac{1}{4}\exp\bigl(- \bigl((3+i)/4 \bigr) \bigl(\log L'
\bigr)^2 \bigr)\qquad\mbox{for }i=1,2,3,4. %
\]
\end{itemize}
Let us summarize this condition in words.

The first point controls the total variation distance of the RWRE exit
measure to the exit measure of simple random
walk on balls of radii $3\leq L\leq2L_0$. Note that the bound on the
probability is given in terms of $L_0$, for all such $L$.

The second point concerns radii $L_0< L\leq L_1$ and gives control over
the deviation of the RWRE exit measure from
that of a symmetric random walk with kernel $p_L$. It
also includes a continuity property of RWRE exit measures when $L'$ varies
(note that we use $p_L$ on the left-hand side, not $p_{L'}$), which
will be
crucial to compare the distance between two kernels for different
radii; see
Lemma~\ref{diff-admiss}.

The main technical statement of this paper is the following:
%
\begin{proposition}
\label{main-prop}
Assume $\mathbf{A1}$. For $\delta>0$ small enough, there exists
$\varepsilon_0 =
\varepsilon_0(\delta) > 0$ with the following property: if
$\varepsilon\leq\varepsilon_0$
and $\mathbf{A0}(\varepsilon)$ holds, then:
\begin{longlist}[(ii)]
\item There exists $L_0 = L_0(\delta)$ such that for $L_1\geq L_0$,
\[
\ctwo(\delta,L_0,L_1 ) \Rightarrow\ctwo\bigl(
\delta,L_0,L_1(\log L_1)^2 \bigr).
\]
\item There exist $L_0=L_0(\delta)$ and sequences $\ell_n$,
$m_n\rightarrow\infty$ with the following property: if $L_1\geq\ell
_n$ and
$L_1\leq L
\leq L_1(\log L_1)^2$, then for $q\in{\mathcal{P}^{\mathrm{s}}_\iota
}$ and $m\geq
m_n$, with $\psi\equiv m$,
\[
\ctwo(\delta,L_0,L_1) \Rightarrow\bigl(\pP
\bigl(D_{L,p_L,\psi,q}^\ast> 1/n \bigr) \leq\exp\bigl(-(\log
L)^2 \bigr) \bigr). %
\]
\end{longlist}
\end{proposition}

%
\begin{remark}
\label{remark-ctwo}
(i) It is important to realize that for every choice of $\delta$ and $L_0$,
we can make sure that $\ctwo(\delta,L_0,L_0)$ is fulfilled, simply by choosing
the perturbation $\varepsilon$ small enough. This observation provides
us with the
base step of the induction in Proposition~\ref{main-prop}(i): once we
know that $\ctwo$ propagates for properly chosen $\delta$ and $L_0$,
we can
choose $\varepsilon$ so
small such that $\ctwo(\delta,L_0,L)$ holds for \textit{all} $L\geq
L_0$.

(ii) One should note that under $\ctwo(\delta,L_0,L_1)$, if $L \leq
L_1(\log L_1)^2$, then $h_{L,r}(x)\leq s_L \leq L_1/2$, so that
$\ctwo(\delta,L_0,L_1)$ can be used to control the exit distributions of
the coarse grained walks inside $V_L$.

(iii) The number $\iota$ defined above the condition bounds the range
of symmetric transition kernels $q$
from which smoothing kernels $\hat{\pi}_\psi^{(q)}$
are built. In
Lemma~\ref{diff-admiss} we will see that under
$\ctwo(\delta,L_0,L_1)$, for $L\leq L_1(\log L_1)^2$, the kernels
$p_L$ are
elements of $\mathcal{P}^{\mathrm{s}}_\iota$.

(iv) If $\ctwo(\delta,L_0,L_1)$ is satisfied, then for any
$3\leq L\leq L_1$ and for all $L'\in[L,2L]$, all
$\psi\in\mathcal{M}_{L'}$ and all $q\in{\mathcal{P}^{\mathrm
{s}}_\iota}$,
\[
\pP\bigl( \bigl\{D_{L',p_L,\psi,q}^\ast> \bigl(\log L'
\bigr)^{-9} \bigr\} \cup\bigl\{D_{L',p_L}^\ast>
\delta\bigr\} \bigr) \leq\exp\bigl(-\bigl(\log L'
\bigr)^2 \bigr). %
\]
\end{remark}

For the rest of this paper, if we write ``assume $\ctwo(\delta,L_0,L_1)$,''
this means that we assume $\ctwo(\delta,L_0,L_1)$ for some $\delta
>0$ and
some $L_1\geq L_0$, where $\delta$ can be chosen arbitrarily small and
$L_0$ arbitrarily large.

\subsection{A short reading guide}\label{readingguide}
The key idea behind our proofs is to compare exit
measures by means of the expansion
%
%
\begin{equation}
\label{readingguide-pert} \Pi_L-\pi_L = \Gh_{L,r}1_{V_L}(
\Ph_{L,r}-\hat{\pi}_{L,r})\pi_L,
\end{equation}
which results
from~(\ref{prel-exitunchanged}),~(\ref{prel-pbe0}) and~(\ref{prel-pbe1}).
Our coarse grained transition kernels are given by exit distributions from
smaller balls inside $V_L$, and we obtain our results by transferring
inductively information on smaller scales to the scale $L$.
The notion of good and bad points, introduced in Section~\ref{smoothbad},
allows us to classify the exit behavior on smaller scales.
If inside $V_L$ all points are good, then the estimates on smaller balls
can be transferred to a (globally smoothed) estimate on the larger ball~$V_L$ (Lemma~\ref{lemma-goodpart}). But \textit{bad} points can appear,
and in fact we have to distinguish four different levels of badness
(Section~\ref{s1}). When bad points are present, it is convenient to
``goodify'' the environment, that is, to replace
bad points by good ones. This important
concept is first explained in Section~\ref{smoothbad} and then further
developed in Section~\ref{super}.

However, for the globally smoothed estimate, we only
have to deal with the case where all bad points are enclosed in a
comparably small region; two or more such regions are too unlikely
(Lemma~\ref{smoothbad-lemmamanybad}). Some special care is required for
the worst class of bad points in the interior of the ball. For environments
containing such points, we slightly modify the coarse graining scheme
inside $V_L$, as described in Section~\ref{super-cgmod}.

In Lemma~\ref{lemma-badpart}, we prove the smoothed estimates on
environments with bad points and show that the degree of badness decreases
by one from one scale to the next.

For exit measures where no or only a local last smoothing step is added
(Section~\ref{nonsmv-exits}, Lemmata~\ref{nonsmv-lemma1}
and~\ref{nonsmv-lemma2}, resp.), bad points near the boundary of
$V_L$ are much more delicate to handle, since we have to take into account
several possibly bad regions. However, they do not occur too frequently
(Lemma~\ref{smoothbad-lemmabbad}) and can be controlled by capacity arguments.

All these estimates require precise bounds on coarse grained Green's
functions, which are developed in Section~\ref{super}. Roughly
speaking, we show
that on environments with no bad points, the coarse grained RWRE Green's
function for the ball is dominated from above by the analogous quantity
coming from simple random walk (or some symmetric perturbation).

In Section~\ref{estimates}, we present various bounds on hitting
probabilities for both symmetric random walk and Brownian motion, and
difference estimates of smoothed exit measures. One main difficulty is that
we have to work with a whole family $(p_L)$ of nearest neighbor transition
kernels. For example, we have to control the total variation distance of
exit measures corresponding to two different kernels. Here, the crucial
statement is Lemma~\ref{exitbmpbmq}, which is formulated in terms of
Brownian motion and then transferred to random walks via coupling
arguments.

The statements from Section~\ref{exit-meas} are
finally used in Section~\ref{proofmain} to prove the main results. In the
\hyperref[appe]{Appendix} we prove the main statements from Section~\ref{estimates}, as well
as a local central limit theorem for the coarse grained symmetric random
walk.

\section{Transition kernels and notion of badness}\label{smoothbad}
Here, we look closer at the family of kernels defined
in~(\ref{kernelpL}) and introduce the concept of ``good'' and
``bad'' points. Furthermore, we define ``goodified'' transition kernels and
prove two estimates ensuring that we do not have to consider environments
with bad points that are widely spread out in the ball or densely
packed in
the boundary region.
\subsection{Some properties of the kernels $p_L$}
The first general statement exemplifies how to extract information
about a symmetric kernel
$p\in\mathcal{P}^{\mathrm{s}}_\kappa$, $0<\kappa<1/(2d)$, out of
the corresponding exit
measure on $\partial V_L$.

%
\begin{lemma}
\label{onestep-exitlaw}
For $i=1,\ldots,d$,
\[
p(e_i) = \frac{1}{2}\sum_{y\in\partial V_L}\pi^{(p)}_L(0,y) \biggl(\frac
{y_i}{L}
\biggr)^2 + O\bigl(L^{-1}\bigr). %
\]
\end{lemma}
\begin{pf}
Recall that under $\Prw_{0,p}$, $(X_n)_{n\geq0}$ denotes the canonical
random walk on $\mathbb{Z}^d$ with transition kernel $p$ starting at the
origin. Write $\mathcal{G}_m = \sigma(X_1,\ldots,X_m)$ for the
filtration up
to time $m$, and denote by $X_{n,i}$ the $i$th component of $X_n$. Due to
the symmetry of $p$, the process $X_{n,i}^2-2p(e_i)n$, $n\geq0$, is a
martingale with respect to $\mathcal{G}_n$. By the optional stopping
theorem,
\[
\Erw_{0,p} \bigl[X_{\tau_L,i}^2 \bigr] =
2p(e_i)\Erw_{0,p}[\tau_L]. %
\]
Since\vspace*{1.5pt} $X_{\tau_L,1}^2+\cdots+X_{\tau_L,d}^2 =(L+O(1))^2$, it follows that
$\Erw_{0,p}[\tau_L] = (L+O(1))^2$, and the claim is proved.
\end{pf}

Now let us turn to the kernel $p_L$.
%
\begin{lemma}
\label{diff-admiss}
Assume $\ctwo(\delta,L_0,L_1)$. There exists a constant $C>0$ such that:
\begin{longlist}[(ii)]
\item
For $3\leq L\leq L_1(\log L_1)^2$,
\[
\llVert p_{s_L/20}-p_{L}\rrVert_1 \leq C(\log
L)^{-9}. %
\]
In particular, if $L_0$ is sufficiently large, we have $p_L\in
\mathcal{P}^{\mathrm{s}}_\iota$.
\item
Let $3\leq L\leq L_1$ and
$L'\in[L/2,L]$. Then
\[
\llVert p_{L'}-p_L\rrVert_1 \leq C(\log
L)^{-9}. %
\]
\end{longlist}
\end{lemma}
\begin{pf}
(i) For $L\leq L_0$, there is nothing to show since $p_L=p_o$. Now
assume $L_0<L\leq L_1(\log L_1)^2$. We apply
Lemma~\ref{onestep-exitlaw} to $t\in[s_L/20,s_L/10]$ in place of
$L$. Writing $p$ for $p_{s_L/20}$ and ${\hat{\pi}}$ for
$\hat{\pi}_{L,r}$, we then
obtain for each $i=1,\ldots,d$ after an integration
%
%
\begin{equation}
\label{onestep-exitlaw-1} 2p(e_i) = \sum_{y\in V_L}{
\hat{\pi}^{(p)}}(0,y)\frac{y_i^2}{\llvert y \rrvert ^2}+ O\bigl(s_L^{-1}
\bigr).
\end{equation}
Therefore, by the definition of $p_L$,
\[
2\bigl\llvert p(e_i)-p_L(e_i)\bigr\rrvert
= \sum_{y} \bigl(\hat{\pi}^{(p)}-\pE[
\Ph] \bigr) (0,y)\frac
{y_i^2}{\llvert y \rrvert ^2} + O\bigl(s_L^{-1}
\bigr). %
\]
We now use the elementary fact that for centered random variables $Y, Y',Z$
with $Y,Y'$ independent of $Z$, we have $E[Y^2]-E[Y'^2] =
E[(Y+Z)^2]-E[(Y'+Z)^2]$. Moreover, $\hat{\pi}^{(p)}(0,\cdot)$ and
$\pE[\Ph
](0,\cdot)$
have support in $V_{s_L/10}$, where the kernel ${\hat{\pi
}}^{(p_o)}$ is
homogeneous. We can therefore write
\[
\sum_{y} \bigl(\hat{\pi}^{(p)}-\pE[
\Ph] \bigr) (0,y)\frac
{y_i^2}{\llvert y \rrvert ^2}=\sum_{y}
\bigl({\hat{\pi}}^{(p)}-\pE[\Ph] \bigr){\hat{\pi}}^{(p_o)}(0,y)
\frac{y_i^2}{\llvert y \rrvert ^2}. %
\]
We next note that by definition of the coarse-grained transition kernels,
we have
\[
\bigl\llVert\bigl(\pE[\Ph]-\hat{\pi}^{(p)} \bigr){\hat{\pi
}}^{(p_o)}(0,\cdot)\bigr\rrVert_1\leq\sup
_{t\in[s_L/20, s_L/10]}\pE\bigl[\bigl\llVert\bigl(\Pi_{V_t}-
\pi^{(p)} _{V_t} \bigr){\hat{\pi}}^{(p_o)}(0,\cdot)
\bigr\rrVert_1 \bigr]. %
\]
We apply condition $\ctwo(\delta,L_0,L_1)$ in order to bound the right-hand
side. Clearly, $s_L/10 \leq L_1$. Moreover, the function $h_{L,r}$ defined
in~(\ref{smoothbad-hlr}) lies in $\mathcal{M}_t$ for each
$t\in[s_L/20,s_L/10]$. Recalling the last point of
Remark~\ref{remark-ctwo}, we obtain under $\ctwo(\delta,L_0,L_1)$
(with $t$
in place of $L'$ and $p=p_{s_L/20}$ in place of $p_L$ in this remark)
\[
\pE\bigl[\bigl\llVert\bigl(\Pi_{V_t}-\pi^{(p)}_{V_t}
\bigr) {\hat{\pi}}^{(p_o)}(0,\cdot)\bigr\rrVert_1 \bigr] \leq
C(\log L)^{-9} %
\]
for some constant $C$ which is uniform in $t\in[s_L/20,s_L/10]$. Putting
the pieces together, we have shown that for each $i$,
$\llvert p(e_i)-p_L(e_i)\rrvert \leq C(\log L)^{-9}$ and the first
part of (i) follows.

In order to see that $p_L\in\mathcal{P}^{\mathrm{s}}_\iota$, we
put $\ell_1=L$,
$\ell_{k+1}=s_{\ell_k}/20$ and then apply the above bound repeatedly
to the
differences $\llVert p_{\ell_k}-p_{\ell_{k+1}}\rrVert _1$, until $\ell
_{k+1}\leq L_0$
and hence $p_{\ell_{k+1}} = p_o$. With $K=\lfloor\log
_2(L/L_0)\rfloor$, we
obtain the bound
\[
\llVert p_L-p_o\rrVert_1\leq C\sum
_{i=1}^K\bigl(\log\bigl(2^{-i}L
\bigr)\bigr)^{-9}\leq C(\log L_0)^{-8},
\]
which implies the second part of (i).

(ii) Let $\psi\equiv L \in\mathcal{M}_L$. By Lemma~\ref{onestep-exitlaw}
and the same variance additivity property as in the proof of (i),
\begin{eqnarray*}
\llVert p_{L'}-p_L\rrVert_1&=&
\frac{1}{L^2}\biggl\llvert\sum_{y\in\partial
V_L} \bigl(
\pi_L^{(p_{L'})}-\pi_L^{(p_L)} \bigr)
(0,y)y_i^2\biggr\rrvert+ O\bigl(L^{-1}\bigr)
\\
&=&\frac{1}{L^2}\biggl\llvert\sum_{y\in\mathbb{Z}^d} \bigl(
\pi_L^{(p_{L'})}-\pi_L^{(p_L)} \bigr){
\hat{\pi}}^{(p_o)}_\psi(0,y)y_i^2\biggr
\rrvert+ O\bigl(L^{-1}\bigr)
\\
&\leq& C\bigl\llVert\bigl(\pi_{L}^{(p_{L'})}-
\pi_L^{(p_L)} \bigr) {\hat{\pi}}^{(p_o)}_\psi(0,
\cdot)\bigr\rrVert_1 +O\bigl(L^{-1}\bigr).
\end{eqnarray*}
Since under $\ctwo(\delta, L_0,L_1)$,
\begin{eqnarray*}
\bigl\llVert\bigl(\pi_{L}^{(p_{L'})}-\pi_L^{(p_L)}
\bigr) {\hat{\pi}}^{(p_o)}_\psi(0,\cdot)\bigr\rrVert
_1&\leq& C \bigl(\pE\bigl[D^\ast_{L,p_L,\psi,p_o} \bigr] +
\pE\bigl[D^\ast_{L,p_{L'},\psi,p_o} \bigr] \bigr)
\\
&\leq& C(\log L)^{-9},
\end{eqnarray*}
the second claim is proved.
\end{pf}

\subsection{Good and bad points}\label{goodandbad}
We shall partition the grid points inside $V_L$ according to
their influence on the exit behavior. Recall assignment~(\ref{kernelpL}),
and fix an environment $\omega\in\Omega$.
We say that a point $x\in V_L$ is \textit{good} (with respect to $\omega$, $L$,
$\delta>0$ and $r$, $100\leq r\leq r_L$) if:
\begin{itemize}
\item For all $t\in[h_{L,r}(x),2h_{L,r}(x)]$, with $q=p_{h_{L,r}(x)}$,
\[
\bigl\llVert\bigl(\Pi_{V_t(x)}-\pi^{(q)}_{V_t(x)}
\bigr) (x,\cdot)\bigr\rrVert_1 \leq\delta.
\]
\item If $\dL(x) > 2r$, then additionally
\[
\bigl\llVert\bigl(\Ph_{L,r} - {\hat{\pi}^{(q)}}_{L,r}
\bigr){\hat{\pi}^{(q)}}_{L,r}(x,\cdot)\bigr\rrVert
_1 \leq\bigl(\log h_{L,r}(x) \bigr)^{-9}.
\]
\end{itemize}
A point $x\in V_L$ which is not good is called \textit{bad}.
We denote by $\badP_{L,r} = \badP_{L,r}(\omega)$ the set of all bad points
inside $V_L$ and write $\badP_L = \badP_{L,r_L}$ for short.\vspace*{1pt} Furthermore,
set $\badP_{L,r}^\partial= \badP_{L,r}\cap\sh_L(r_L)$
and $\mathcal{B}^\star_{L,r} =
\badP_{L,r}\cup\badP_{L} = \badP_{L,r}^\partial\cup\badP_L$.
Of course, the set of bad points depends also on $\delta$, but we do
not indicate this.

%
\begin{remark}
\label{remark-smoothbad}
(i) For the coarse graining scheme associated to $r=r_L$, we have by
definition $\mathcal{B}^\star_{L,r_L} = \badP_L$. When performing
the nonsmoothed
estimates in Section~\ref{nonsmv-exits}, we work with constant $r$. In
this
case, $\mathcal{B}^\star_{L,r}$ can contain more points than $\badP
_L$.

(ii) Assume $L$ large. If $x\in V_L$ with $\dist_L(x) > 2r$, then the
function \mbox{$h_{L,r}(x+\cdot)$} lies in
$\mathcal{M}_t$ for each $t\in[h_{L,r}(x),2h_{L,r}(x)]$. Thus, for all
$x \in V_L$, we can use $\ctwo(\delta,L_0,L_1)$ to control the event
$\{x\in\badP_{L,r}\}$, provided $2h_{L,r}(x) \leq L_1$.
\end{remark}

We shall replace the RWRE transition kernels at bad points by those of a
symmetric random walk. Write $p$ for $p_{s_L/20}$. For all environments, we introduce
the ``goodified'' transition kernels as follows:
%
%
\begin{equation}
\label{smoothbad-goodifiedkernel} {\hat{\Pi}^g}_{L,r}(x,\cdot) = \cases{
\Ph_{L,r}(x,\cdot), &\quad for $x \in V_L\setminus\mathcal{B}^\star_{L,r}$,
\vspace*{2pt}\cr
\hat{\pi}^{(p)}_{L,r}(x,
\cdot), &\quad for $x \in\mathcal{B}^\star_{L,r}$.}
\end{equation}
We write ${\hat{G}^g}_{L,r}$ for the
corresponding (random) Green's function. Note that the transition kernel
$q$ used in the definition of a good point $x\in V_L$ does
depend on the location of $x$ inside the ball, whereas the
goodifying-procedure uses the same transition kernel $p$ for all points
[which agrees with $q$ for $x\in V_L$ with $\dL(x)\geq2s_L$, since
in this region $h_{L,r}\equiv(1/20)s_L$].
Goodified transition kernels and Green's functions will play a major role
from Section~\ref{super} onwards.

\subsection{Bad regions in the case $r=r_L$}
\label{s1}
The next lemma shows that with high probability, all bad points
with respect to $r=r_L$ are contained in a ball of radius $4h_L(x)$. Let
\[
\mathcal{D}_L = \bigl\{V_{4h_L(x)}(x)\dvtx  x\in
V_L \bigr\}. %
\]

%
%
\begin{figure}

\includegraphics{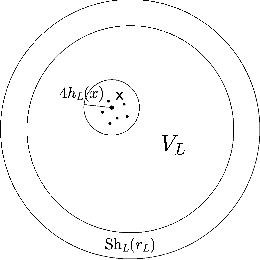}

\caption{On environments $\omega\in\onebad$, all bad points are
enclosed in a ball $V_{4h_L(x)}(x)$.}\label{fig2}
\end{figure}

\noindent We will look at the events $\onebad= \{\badP_L\subset D$ for
some $ D\in\mathcal{D}_L \}$, see Figure~\ref{fig2}, and
$\manybad={ (\onebad)}^c$. It is also useful to define
the set
of \textit{good} environments, $\good= \{\badP_L =
\varnothing\}\subset\onebad$.
%
\begin{lemma}
\label{smoothbad-lemmamanybad}
Assume $\ctwo(\delta,L_0,L_1)$. Then for
$L_1\leq L \leq L_1(\log L_1)^2$,
\[
\pP(\manybad) \leq\exp\bigl(-\tfrac{19}{10}(\log L)^2
\bigr). %
\]
\end{lemma}
\begin{pf}
Let $x \in V_L$ with $\dL(x) > 2r_L$. Set $q=p_{h_L(x)}$
and $\Delta= 1_{V_L}(\Ph_{L,r_L} - {\hat{\pi}^{(q)}}_{L,r_L})$. Put
$D_{t,q,h_L,q}(x)=\llVert (\Pi_{V_t(x)}-\pi^{(q)}_{V_t(x)}){\hat{\pi
}^{(q)}}_{h_L}(x,\cdot
)\rrVert _1$,
$D_{t,q}(x)=\llVert (\Pi_{V_t(x)}-\pi^{(q)}_{V_t(x)})(x,\cdot
)\rrVert _1$.
Using\vspace*{1.5pt} $r_L/20 \leq
h_L(x) \leq s_L\leq L_1/2$ and the second point of Remark~\ref{remark-smoothbad},
\begin{eqnarray*}
\pP(x\in\badP_{L} )
&\leq&\pP\biggl(\bigcup_{t\in[h_L(x),2h_L(x) ]} \bigl\{
D_{t,q,h_L,q}(x)>\bigl(\log h_L(x)\bigr)^{-9} \bigr\}
\cup\bigl\{D_{t,q}(x) > \delta\bigr\} \biggr)
\\
&\leq& Cs_L^{d}\exp\bigl(- \bigl(\log(r_L/20)
\bigr)^{2} \bigr),
\end{eqnarray*}
and a similar estimate holds when $\dL(x) \leq2r_L$. On the event
$\manybad$, there exist $x,y \in\badP_L$ with $\llvert x-y\rrvert > 2h_L(x)
+2h_L(y)$. But for such $x,y$, the events $\{x\in\badP_L\}$ and $\{
y\in
\badP_L\}$ are independent, whence for $L$ large
\begin{eqnarray*}
\pP(\manybad) &\leq& C L^{2d}s_L^{2d} \bigl[
\exp\bigl(-\bigl(\log(r_L/20)\bigr)^2 \bigr)
\bigr]^2
\\
&\leq&\exp\bigl(-(19/10) (\log L)^2 \bigr).
\end{eqnarray*}\upqed
\end{pf}

The estimate is good enough for our inductive procedure, so we only
have to
deal with the case where all possibly bad points are enclosed in a ball
$D\in\mathcal{D}_L$. However, inside $D$ we need to look closer at the
degree of badness.

We say that $\omega\in\onebad$ is \textit{bad on level}
$i$, $i=1,2,3$, if the following holds:
\begin{itemize}
\item For all $x\in V_L$, for all $t\in[h_L(x),2h_L(x)]$, with
$q=p_{h_{L}(x)}$,
\[
\bigl\llVert\bigl(\Pi_{V_t(x)} - \pi^{(q)}_{V_t(x)}
\bigr) (x,\cdot)\bigr\rrVert_1 \leq\delta. %
\]
\item For all $x\in V_L$ with $\dL(x) > 2r_L$, additionally
\[
\bigl\llVert\bigl(\Ph_{L,r_L} - {\hat{\pi}^{(q)}}_{L,r_L}
\bigr){\hat{\pi}^{(q)}} _{L,r_L}(x,\cdot)\bigr\rrVert
_1 \leq\bigl(\log h_L(x) \bigr)^{-9+9i/4}.
\]
\item There exists $x\in\badP_L(\omega)$ with $\dL(x) > 2r_L$ such that
\[
\bigl\llVert\bigl(\Ph_{L,r_L} - {\hat{\pi}^{(q)}}_{L,r_L}
\bigr){\hat{\pi}^{(q)}} _{L,r_L}(x,\cdot)\bigr\rrVert
_1 > \bigl(\log h_L(x) \bigr)^{-9+9(i-1)/4}.
\]
\end{itemize}
If $\omega\in\onebad$ is neither bad on level $i=1,2,3$ nor good, we call
$\omega$ \textit{bad on level} $4$. In this case, $\badP_L(\omega)$ contains
``really bad'' points. We write $\operatorname{OneBad}_L^{(i)}\subset
\onebad$ for the
subset of
all those $\omega$ which are bad on level $i=1,2,3,4$. Observe that we have
the partition
\[
\onebad= \good\cup\,\bigl(\onebad^{(1)}\cup\cdots\cup\onebad
^{(4)} \bigr). %
\]
On $\good$, ${\hat{\Pi}^g}_{L,r_L}=\Ph_{L,r_L}$ and therefore
${\hat{G}^g}
_{L,r_L}=\Gh_{L,r_L}$.

\subsection{Bad regions when $r$ is a constant}
\label{s2}
When estimating nonsmoothed exit measures, we cannot stop the
refinement of the coarse graining in the boundary region $\sh_L(r_L)$.
Instead, we
will choose $r$ as a (large) constant. However, now it\vadjust{\goodbreak} is no longer true
that essentially all bad points are contained in one single region
$D\in\mathcal{D}_L$. For example, if $x\in V_L$ is such that
$\dist_L(x)$ is of order $\log L$, we only have a bound of the form
\[
\pP(x\in\badP_{L,r} ) \leq\exp\bigl(-c{ (\log\log L
)}^2 \bigr), %
\]
which is clearly not enough to get an estimate as in
Lemma~\ref{smoothbad-lemmamanybad}. We therefore choose a different
strategy to handle bad points within $\sh_L(r_L)$. We split the boundary
region into layers of an appropriate size and use independence to show that
with high probability, bad regions are rather sparse within those
layers. Then the Green's function estimates of
Corollary~\ref{super-cor} will ensure that on such environments, there
is a high chance to never hit points in $\badP_{L,r}^\partial$ before
leaving the ball.\vspace*{1pt}

%
%
\begin{figure}

\includegraphics{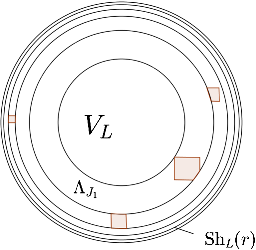}

\caption{The layers $\Lambda_j$, $0\leq j\leq J_1$, with $\Lambda
_0=\sh_L(2r)$. Subsets $D_{\mathbf k}^{(j)}\subset\Lambda_j$ containing bad points are shaded.}\label{fig3}
\end{figure}

To begin with the first part, fix $r$ with $r\geq r_0\geq
100$, where $r_0 = r_0(d)$ is a constant that will be chosen
below. Let $L$ be large enough such that $r < r_L$, and set $J_1 = J_1(L)
= \lfloor\log_2(r_L/r) \rfloor+1$. We define
layers $\Lambda_0 = \sh_L(2r)$ and $\Lambda_j = \sh
_L(r2^j,r2^{j+1})$ for integers
$1\leq j\leq J_1$. Then
\[
\sh_L(2r_L) \subset\bigcup
_{0\leq j\leq J_1}\Lambda_j\subset\sh_L(4r_L).
\]
Let $1\leq j\leq J_1$. For $k\in\mathbb{Z}$, consider the interval
$I_k^{(j)} = (kr2^j,(k+1)r2^j]\cap\mathbb{Z}$. We divide $\Lambda_j$
into subsets by setting $D_{\mathbf k}^{(j)} = \Lambda_j\cap
(I_{k_1}\times\cdots\times I_{k_d} )$, where ${\mathbf k} =
(k_1,\ldots,k_d) \in\mathbb{Z}^d$, cf. Figure~\ref{fig3}. Denote by $\mathcal{Q}_{j,r}$
the set
of those subsets which are not empty. Setting $N_{j,r} =
\llvert \mathcal{Q}_{j,r}\rrvert $, it follows that
\[
\frac{1}{C}{ \biggl(\frac{L}{r2^j} \biggr)}^{d-1} \leq
N_{j,r} \leq C{ \biggl(\frac{L}{r2^j} \biggr)}^{d-1}.
\]
We\vspace*{-1pt} say that a set $D\in\mathcal{Q}_{j,r}$ is \textit{bad} if
$\badP^\partial_{L,r}\cap D \neq\varnothing$. As we want to make use of
independence, we partition $\mathcal{Q}_{j,r}$ into disjoint sets
$\mathcal{Q}_{j,r}^{(1)},\ldots, \mathcal{Q}_{j,r}^{(R)}$, such that
for each
$1\leq m \leq R$, we have:
\begin{itemize}
\item$\dist(D,D') > 4\max_{x\in\Lambda_j}h_{L,r}(x)$ for all
$D\neq
D'\in\mathcal{Q}_{j,r}^{(m)}$,

\item$N_{j,r}^{(m)} = \llvert \mathcal{Q}_{j,r}^{(m)}\rrvert \geq
\frac{N_{j,r}}{2R}$.
\end{itemize}
Note that since $h_{L,r}$ is proportional to $\dL$ on $\Lambda_j$ for
$1\leq j\leq J_1$, the number \mbox{$R\in\mathbb{N}$} can be chosen to
depend on
the dimension only. Then the events $\{D\mbox{ is bad}\}$,
$D\in\mathcal{Q}_{j,r}^{(m)}$, are\vspace*{1.5pt} independent. Furthermore, if
$L_1\leq L \leq
L_1(\log L_1)^2$, it follows that under $\ctwo(\delta,L_0,L_1)$,
\begin{eqnarray*}
\pP(D\mbox{ is bad} ) &\leq& C\bigl(r2^j\bigr)^{2d}\exp
\bigl(- \bigl(\log\bigl(r2^j/20\bigr) \bigr)^2 \bigr)
\\
&\leq&\exp\bigl(- (\log r + j )^{5/3} \bigr)= p_{j,r},
\end{eqnarray*}
for all $r \geq r_0$ and $j\in\mathbb{N}$, if $r_0$ is big enough. Let
$Y_{j,r}$ and $Y_{j,r}^{(m)}$ be the number of bad sets in
$\mathcal{Q}_{j,r}$ and $\mathcal{Q}_{j,r}^{(m)}$, respectively. For
$r\geq
5$, we have $p_{j,r} \leq(\log r +j)^{-3/2}\leq1/2$. A standard large
deviation estimate for Bernoulli random variables yields
\[
\pP\bigl(Y_{j,r}^{(m)} \geq(\log r + j)^{-3/2}N_{j,r}^{(m)}
\bigr) \leq\exp\bigl(-N_{j,r}^{(m)} I \bigl((\log r +
j)^{-3/2} \mid p_{j,r} \bigr) \bigr), %
\]
with $I(x \mid p) = x\log(x/p) + (1-x)\log((1-x)/(1-p))$. By enlarging $r_0$
if necessary, we get $I ((\log r + j)^{-3/2} \mid p_{j,r} )
\geq2
R(\log r+j)^{1/7}$ for $r\geq r_0$, whence
\begin{eqnarray*}
&& \pP\bigl(Y_{j,r} \geq(\log r + j)^{-3/2}N_{j,r}
\bigr)
\\
&&\qquad\leq R\max_{m=1,\ldots,R}\pP\bigl(Y_{j,r}^{(m)}
\geq(\log r + j)^{-3/2}N_{j,r}^{(m)} \bigr)
\\
&&\qquad\leq R\exp\bigl(-(\log r +j)^{1/7}N_{j,r} \bigr)
\\
&&\qquad \leq R\exp
\biggl(-\frac{1}{C} (\log r +j)^{1/7}{ \biggl(\frac{L}{r2^j}
\biggr)}^{d-1} \biggr)
\\
&&\qquad \leq\exp\bigl(-(\log r +j)^{1/7}(\log L)^{29} \bigr),
\end{eqnarray*}
for $r_0\leq r < r_L$, $0\leq j\leq J_1(L)$ and $L$ large enough. In
particular,
\[
\sum_{0\leq j\leq
J_1(L)}\pP\bigl(Y_{j,r}\geq(\log r +
j)^{-3/2}N_{j,r} \bigr) \leq\exp\bigl(-(\log
L)^{28} \bigr). %
\]
Therefore, introducing the set of environments with plenty of bad
points in
the boundary region,
\[
\bbad= \bigcup_{0\leq j\leq J_1(L)} \bigl\{Y_{j,r} \geq
(\log r + j)^{-3/2}N_{j,r} \bigr\}, %
\]
we have proved the following:
%
\begin{lemma}
\label{smoothbad-lemmabbad}
There exists a constant $r_0 >0$ such that if $r\geq r_0$, then
$\ctwo(\delta,L_0,L_1)$ implies that for $L_1\leq L \leq L_1(\log L_1)^2$,
\[
\pP(\bbad) \leq\exp\bigl(-(\log L)^{28} \bigr). %
\]
\end{lemma}

\section{Some important estimates}
\label{estimates}
In this section, we collect estimates on symmetric random walks with kernel
$p\in\mathcal{P}^{\mathrm{s}}_\kappa$ and on
$d$-dimensional Brownian motion with (diagonal) covariance matrix given
by
%
%
\begin{equation}
\label{covariance-matrix-def} \Lambda_p= \bigl(2d p(e_i)
\delta_{i,j} \bigr)_{i,j=1}^d. 
\end{equation}
We can safely use the same letter as for the layers defined in the
foregoing section, since it will always be clear from the context what is
meant. The following statements hold for small $\kappa$, meaning that
there exists $0<\kappa_0<1/(2d)$ such that for $0<\kappa\leq\kappa_0$,
the statements hold true. All constants are then uniform in $p\in
\mathcal{P}^{\mathrm{s}}_\kappa$.

\subsection{Hitting probabilities}
The first two lemmata concern symmetric random walk. The proofs are
provided in the \hyperref[appe]{Appendix}.
%
\begin{lemma}
\label{lemmalawler}
Let $p\in\mathcal{P}^{\mathrm{s}}_\kappa$, and let $0<\eta< 1$.
\begin{longlist}[(ii)]
\item There exists $C=C(\eta)>0$ such that
for all $x\in V_{\eta L}$,
$z\in\partial V_L$,
\[
C^{-1}L^{-d+1} \leq \pi^{(p)}_L(x,z)
\leq CL^{-d+1}. %
\]
\item
There exists $C=C(\eta)>0$ such that for all $x,x'\in V_{\eta L}$,
$z\in\partial V_L$,
\[
\bigl\llvert \pi^{(p)}_L(x,z)-\pi^{(p)}_L
\bigl(x',z\bigr)\bigr\rrvert\leq C\bigl\llvert x-x'
\bigr\rrvert L^{-d}. %
\]
\end{longlist}
\end{lemma}
A good control over hitting probabilities is given by the following:
%
\begin{lemma}
\label{hittingprob}
Let $a \geq1$ and $x,y \in
\mathbb{Z}^d$ with $x\notin V_a(y)$. There exists a constant $C>0$ such
that for $p\in\mathcal{P}^{\mathrm{s}}_\kappa$,
\begin{longlist}[(iii)]
\item
\[
\Prw_{x,p} (T_{V_a(y)} < \infty) \leq C \biggl(
\frac
{a}{\llvert x-y\rrvert } \biggr)^{d-2}. %
\]
\item There exists $C>0$,
independent of $a$, such that when $\llvert x-y\rrvert >7a$,
\[
\Prw_{x,p} (T_{V_a(y)} < \tau_L ) \leq C
\frac{a^{d-2}\max
\{a,\dL(y)\}\max\{1,\dL(x)\}}{\llvert x-y\rrvert ^d}. %
\]
\item There exists $C>0$ such that for all $x\in V_L$, $z\in\partial V_L$,
\[
C^{-1}\frac{\dL(x)}{\llvert x-z \rrvert ^{d}}\leq\pi^{(p)}_L(x,z)
\leq C\frac
{\max\{
1,\dL(x)\}}{\llvert x-z \rrvert ^{d}}. %
\]
\end{longlist}
\end{lemma}
We need to compare exit laws of random walks with different kernels
$p\in\mathcal{P}^{\mathrm{s}}_\kappa$, and we need difference
estimates on smoothed exit measures.
In this direction,\vadjust{\goodbreak} it is easier to work with Brownian motion and then
transfer the results back to the discrete setting. Let us first introduce
some additional notation. Let $p\in{\mathcal{P}^{\mathrm{s}}_\kappa
}$. For a domain
$U\subset\mathbb{R}^d$ with smooth boundary and $x\in U$, denote by
$\pi^{\mathrm{B}(p)}_U(x,\dz z)$ the exit measure from $U$ of a
$d$-dimensional Brownian
motion $W_t$ started at $x$, with diffusion (or covariance) matrix
$\Lambda_p$ defined in~(\ref{covariance-matrix-def}),
that is, $\pE[(W_1-x)^2] = \Lambda_p$. In the case $U=C_L$, we simply write
$\pi^{\mathrm{B}(p)}_L(x,\dz z)$. By a small abuse of notation, we also write
$\pi^{\mathrm{B}(p)}_U(x,z)$ [or $\pi^{\mathrm{B}(p)}_L(x,z)$ if
$U=C_L$] for the (continuous version
of the) density with respect to surface measure on $U$.

In particular, $\pi^{\mathrm{B}(p_o)}_U$ is the exit measure from $U$
of standard
$d$-dimensional Brownian motion with covariance matrix $I_d$. Its density
$\pi^{\mathrm{B}(p_o)}_{L}(x,z)$ is given by the Poisson kernel
%
%
\begin{equation}
\label{est-poissonkernel} \pi^{\mathrm{B}(p_o)}_L(x,z)= \frac{1}{d
\alpha(d) L}
\frac
{L^2-\llvert x \rrvert ^2}{\llvert x-z \rrvert ^d},
\end{equation}
where $\alpha(d)$ is the volume of the unit ball. For general $p\in
\mathcal{P}^{\mathrm{s}}_\kappa$,
there is no explicit expression for the kernel $\pi^{\mathrm
{B}(p)}_L(x,z)$. However,
we have the following:
%
\begin{lemma}
\label{exitdensitybm}
There exists $C>0$ such that for $p\in{\mathcal{P}^{\mathrm
{s}}_\kappa}$ and all $x\in
C_L$, \mbox{$z\in\partial C_L$},
\begin{longlist}[(ii)]
\item
\[
C^{-1}\frac{\dL(x)}{\llvert x-z \rrvert ^{d}}\leq\pi^{\mathrm{B}(p)}_L(x,z)
\leq C\frac{\dL
(x)}{\llvert x-z \rrvert ^{d}}. %
\]
\item For $k\in\mathbb{N}$,
\[
C^{-1}\frac{\dL(x)}{\llvert x-z \rrvert ^{d+k}}\leq\nabla_x^k
\pi^{\mathrm
{B}(p)}_L(x,z) \leq C\frac{\dL(x)}{\llvert x-z \rrvert ^{d+k}}. %
\]
\end{longlist}
\end{lemma}
This lemma gives us immediately the statements corresponding to
Lem\-ma~\ref{lemmalawler} for Brownian motion with covariance matrix
$\Lambda_p$, $p\in\mathcal{P}^{\mathrm{s}}_\kappa$. Clearly,
also Lemma~\ref{hittingprob}
has a
direct analog. In fact, part~(iii) is reformulated for Brownian motion in
the last lemma. For the results corresponding to (i) and (ii), one can
follow the proof of Lemma~\ref{hittingprob} in the \hyperref[appe]{Appendix}, replacing
the random walk estimates by those for Brownian motion. These analogous
results will be used in the \hyperref[appe]{Appendix}.

The following important lemma controls the difference of two Brownian exit
densities on $\partial C_L$, when the corresponding diffusion matrices are
close together.
%
\begin{lemma}
\label{exitbmpbmq}
There exists $C>0$ such that for
$p,q\in\mathcal{P}^{\mathrm{s}}_\kappa$, for all $x\in
C_{(2/3)L}$, $z\in\partial C_L$,
\[
\bigl\llvert\bigl(\pi^{\mathrm{B}(p)}_{L}-\pi^{\mathrm{B}(q)}_{L}
\bigr) (x,z)\bigr\rrvert\leq C\llVert q-p\rrVert_1L^{-(d-1)}.
\]
\end{lemma}
The proof involves techniques from the theory of elliptic PDEs and is given
in the \hyperref[appe]{Appendix}, as well as the proof of the foregoing lemma. It should be
pretty clear that differences of exit probabilities of symmetric random
walks can be bounded in the same way, that is,
\[
\bigl\llvert\bigl(\pi^{(p)}_{L}-\pi^{(q)}_{L}
\bigr) (x,z)\bigr\rrvert\leq C\llVert q-p\rrVert_1L^{-(d-1)}.
\]
However, it seems more difficult to prove this, and in any
case, we will only need a weaker form, which can be readily deduced
from the
last lemma and a coupling argument given in the \hyperref[appe]{Appendix}.
%
\begin{lemma}
\label{php-phq-estimate}
There exists $C>0$ such that for $p,q\in{\mathcal{P}^{\mathrm
{s}}_\kappa}$, for
large $L$, $\psi\in\mathcal{M}_L$, any $x\in\mathcal{U}_L\cap
\mathbb{Z}^d$ and any $z\in\mathbb{Z}^d$,
\[
\bigl\llvert\bigl(\hat{\pi}^{(p)}_{\psi}-{\hat{
\pi}^{(q)}}_{\psi
} \bigr) (x,z)\bigr\rrvert\leq C\llVert q-p
\rrVert_1L^{-d}. %
\]
Moreover, for $x\in V_L$ with $\dL(x) > (1/10)r$,
\[
\bigl\llVert\bigl(\hat{\pi}^{(p)}_{L,r}-{\hat{\pi
}^{(q)}}_{L,r} \bigr) (x,\cdot)\bigr\rrVert_1
\leq C\max\bigl\{h_{L,r}(x)^{-1/4}, \llVert q-p\rrVert
_1\bigr\}. %
\]
\end{lemma}
\begin{pf}
By\vspace*{1.5pt} comparing $\hat{\pi}^{(p)}_{\psi}$ to the analogous Brownian quantity
${\hat{\pi}^{\mathrm{B}(p)}}_{\psi}$ defined in~(\ref{app-cgpsibmdens}),
the first claim follows from Lemma~\ref{exitbmpbmq},
Lemma~\ref{app-kernelest}(vii) from the \hyperref[appe]{Appendix} and the triangle
inequality. The second statement is proved in the same way, with
the choice $\psi(x)= h_{L,r}(x)$. The restriction to $x$ with
$\dL(x)>(1/10)r$ ensures that all exit distributions are taken from balls
that lie completely inside $V_L$.
\end{pf}\vspace*{-10pt}

Let us finish this part by proving the following useful estimate.

%
\begin{lemma}
\label{hittingprob-technical}
Let $a > 0$, $\ell, m\geq1$ and $x\in\mathbb{Z}^d$. Set $R_\ell=
V_\ell\setminus V_{\ell-1}$,
$\alpha=\max\{\llvert \llvert x \rrvert -\ell\rrvert,a \}$. Then for
some constant
$C = C(m) > 0$
\[
\sum_{y\in R_\ell}\frac{1}{ (a+\llvert x-y\rrvert )^m} \leq C \cases{
\ell^{d-(m+1)}, &\quad for $1 \leq m < d-1$,
\vspace*{2pt}\cr
\max\bigl\{\log(\ell/\alpha),
1\bigr\}, &\quad for $m = d-1$,
\vspace*{2pt}\cr
\alpha^{d-(m+1)}, &\quad for $m \geq d$.}
\]
\end{lemma}
\begin{pf}
If $\alpha> \ell$, then the left-hand side
is bounded by
\[
C\ell^{d-1}\alpha^{-m}\leq C\max\bigl\{\alpha^{d-(m+1)},
\ell^{d-(m+1)} \bigr\}. %
\]
If $\alpha\leq
\ell$, we set $A_k = \{y\in R_\ell\dvtx  \llvert x-y\rrvert \in
[(k-1)\alpha,k\alpha)\}$. Then, for all $k\geq1$,
\[
\max_{y\in A_k}\frac{1}{ (a+\llvert x-y\rrvert )^m} \leq2^mk^{-m}
\alpha^{-m}. %
\]
Since for $k\alpha\leq\ell/10$, we have $\llvert A_k \rrvert \leq
C\alpha(k\alpha
)^{d-2}$,
the claim then follows from
\begin{eqnarray*}
&& \sum_{y\in R_\ell}\frac{1}{ (a+\llvert x-y\rrvert )^m}
\\
&&\qquad \leq C \biggl
(\sum
_{1\leq k\leq
\lfloor\ell/(10\alpha)\rfloor}\frac{\alpha(k\alpha
)^{d-2}}{(k\alpha)^m} \biggr) + C
\ell^{d-1}\ell^{-m}
\\
&&\qquad \leq C\alpha^{d-{(m+1)}}\sum_{1\leq k\leq
\lfloor\ell/(10\alpha)\rfloor}k^{d-(m+2)}
+ C\ell^{d-(m+1)}.
\end{eqnarray*}\upqed
\end{pf}

\subsection{Smoothed exit measures}
In order to obtain difference estimates for smoothed exit distributions of
a symmetric random walk, we will compare them to the corresponding
quantities of Brownian motion.

Let $p,q\in\mathcal{P}^{\mathrm{s}}_\kappa$, and let $\psi
= (m_x )\in\mathcal{M}_L$.
The smoothed exit distribution from~$V_L$ of the random
walk (with respect to $p$, $q$, $\psi$) is defined as
\begin{eqnarray*}
\phi_{L,p,\psi,q}(x,z) &=& \pi^{(p)}_L{\hat{
\pi}^{(q)}}_{\psi
}(x,z)
\\
&=& \sum_{y\in\partial
V_L}\pi^{(p)}_L(x,y)
\frac{1}{m_y}\int_{\mathbb{R}_+}\varphi\biggl(\frac
{t}{m_y}
\biggr) \pi^{(q)}_{V_t(y)}(y,z)\,\dt t.
\end{eqnarray*}
For Brownian motion, the smoothing step is defined analogously
to~(\ref{prel-cgpsi}), namely,
\[
{\hat{\pi}^{\mathrm{B}(q)}}_{\psi}(x,\dz z) = \frac{1}{m_x}\int
_{\mathbb{R}_+}\varphi\biggl(\frac{t}{m_x} \biggr)
\pi^{\mathrm{B}(q)}_{C_t(x)}(x,\dz z)\,\dt t. %
\]
The smoothed exit distribution from $C_L$ is then given by
\begin{eqnarray*}
{\phi}^{\mathrm{B}}_{L,p,\psi,q}(x,\dz z) &=& \pi^{\mathrm{B}(p)}_L{
\hat{\pi}^{\mathrm{B}(q)}}_{\psi}(x,\dz z)
\\
&=& \int_{\partial
C_L}\pi^{\mathrm{B}(p)}_L(x,dy)
\frac{1}{m_y}\int_{\mathbb{R}_+}\varphi\biggl(\frac{t}{m_y}
\biggr) \pi^{\mathrm{B}(q)}_{C_t(y)}(y,\dz z)\,\dt t.
\end{eqnarray*}
By\vspace*{1.5pt} ${\phi}^{\mathrm{B}}_{L,p,\psi,q}(x,z)$ we denote the density of
${\phi}^{\mathrm{B}}
_{L,p,\psi,q}(x,\dz z)$
with respect to \mbox{$d$-}dimensional Lebesgue measure.
For the proof of the next lemma, we refer to the \hyperref[appe]{Appendix}.
%
\begin{lemma}
\label{est-phi}
There exists $C>0$ such that for $p,q\in{\mathcal{P}^{\mathrm
{s}}_\kappa}$ and $\psi\in
\mathcal{M}_L$:
\begin{longlist}[(iii)]
\item
\[
\sup_{x\in V_L}\sup_{z\in\mathbb{Z}^d}\bigl\llvert\bigl(
\phi_{L,p,\psi,q}-{\phi}^{\mathrm{B}}_{L,p,\psi,q} \bigr) (x,z)\bigr
\rrvert\leq CL^{-(d+1/4)}, %
\]
\item
\[
\sup_{z\in\mathbb{R}^d}{\bigl\llVert D^i{\phi}^{\mathrm{B}}_{L,p,\psi,q}(
\cdot,z)\bigr\rrVert}_{C_L} \leq CL^{-(d+i)},\qquad i =0,1,2,3,
\]
\item \mbox{}\vspace*{-6pt}\mbox{}
\begin{eqnarray*}
&& \sup_{x,x'\in V_L \cup\partial V_L}\sup_
{z\in\mathbb{Z}^d}\bigl\llvert
\phi_{L,p,\psi,q}(x,z)-\phi_{L,p,\psi,q}\bigl(x',z\bigr)\bigr
\rrvert
\\
&&\qquad\leq C \bigl(L^{-(d+1/4)} + \bigl\llvert x-x' \bigr\rrvert
L^{-(d+1)} \bigr).
\end{eqnarray*}
\end{longlist}
\end{lemma}
The next proposition will be
applied at the end of the proof of Lemma~\ref{lemma-goodpart}. At this
point, the symmetry condition $\mathbf{A1}$ comes into play. We give a
general formulation in terms of a signed\vspace*{1pt} measure $\nu$. Let us introduce
the following notation. For $x=(x_1,\ldots,x_d)\in\mathbb{Z}^d$,
$i=1,\ldots,d$, put
\[
x^{(i)} = (x_1,\ldots,x_{i-1},-x_i,x_{i+1},\ldots,x_d ). %
\]
%
%
\begin{proposition}
\label{est-symmetry}
Let $p,q\in\mathcal{P}^{\mathrm{s}}_\kappa$ and $\ell>0$.
Consider a measure $\nu$ on
$V_\ell$ with
total mass zero satisfying $\nu(x) = \nu(x^{(i)})$ for all $x$ and all
$i=1,\ldots, d$. Then there is a constant $C>0$ such that for $y'\in
V_L$ with $V_\ell(y')
\subset V_L$ and all $z\in\mathbb{Z}^d$, $\psi\in\mathcal{M}_L$,\vspace*{-1pt}
\[
\biggl\llvert\sum_{y\in V_\ell(y')}\nu\bigl(y-y'
\bigr)\phi_{L,p,\psi,q}(y,z)\biggr\rrvert\leq C\llVert\nu\rrVert
_1 \biggl(L^{-(d+1/4)} + \biggl(\frac{\ell}{L}
\biggr)^2L^{-d} \biggr). %
\]
\end{proposition}
\begin{pf}
We simply write $\phi$ for $\phi_{L,p,\psi,q}$ and ${\phi}^{\mathrm
{B}}$ for
${\phi}^{\mathrm{B}}_{L,p,\psi,q}$. Since the proof is the same for
all $y'\in V_L$
with $V_\ell(y') \subset V_L$, we can assume $y'=0$. By
Lemma~\ref{est-phi}(i),
\[
\biggl\llvert\sum_{y}\nu(y)\phi(y,z) - \sum
_{y}\nu(y){\phi}^{\mathrm{B}}(y,z)\biggr\rrvert
\leq C\llVert\nu\rrVert_1L^{-(d+1/4)}. %
\]
Taylor's expansion gives
%
%
\begin{eqnarray}\label{est-symmetry-1}
\qquad&& \sum_{y}\nu(y){
\phi}^{\mathrm{B}}(y,z)\nonumber
\\[-2pt]
&&\qquad =\sum_{y}\nu(y) \bigl[{\phi}^{\mathrm{B}}(y,z)
- {\phi}^{\mathrm{B}}(0,z) \bigr]
\\[-2pt]
&&\qquad =\sum_{y}\nu(y)\nabla_x{
\phi}^{\mathrm{B}}(0,z)\cdot y + \frac{1}{2}\sum
_{y}\nu(y)y\cdot H_x{\phi}^{\mathrm{B}}(0,z)y
+ R(\nu,0,z),\nonumber
\end{eqnarray}
where $\nabla_x{\phi}^{\mathrm{B}}$ is the gradient, $H_x{\phi
}^{\mathrm{B}}$ the Hessian of
${\phi}^{\mathrm{B}}$
with respect to the first variable and $R(\nu,0,z)$ is the remainder
term. Due to the symmetry condition on $\nu$, the first summand on the
right-hand side of~(\ref{est-symmetry-1}) vanishes, and for the second and
third summand one can use Lemma~\ref{est-phi}(ii).
\end{pf}
%
%
\begin{remark}
In~\cite{BZ} and~\cite{BA}, it is assumed that $\mu$, the measure
governing the environment, is isotropic. This\vspace*{1pt} leads us to consider a measure
$\nu$ that is invariant not only under $x\mapsto x^{(i)}$, but also under
$x\mapsto x^{\leftrightarrow(i,j)}$, where for $i<j$,
\[
x^{\leftrightarrow(i,j)} = (x_1,\ldots,x_{i-1},x_j,x_{i+1},
\ldots,x_{j-1},x_i,x_{j+1},\ldots,x_d
). %
\]
In this case, the choice $p=p_o$ results in the sharper bound
\[
\biggl\llvert\sum_{y\in V_\ell(y')}\nu\bigl(y-y'
\bigr)\phi_{L,p_o,\psi,q}(y,z)\biggr\rrvert\leq C\llVert\nu\rrVert
_1 \biggl(L^{-(d+1/4)} + \biggl(\frac{\ell}{L}
\biggr)^3L^{-d} \biggr); %
\]
see Proposition~3.1 in~\cite{BA}. It is then clear from the proof of
Lemma~\ref{lemma-goodpart} that one can work with $p_L = p_o$ for \textit{all} radii $L$, that is, the (isotropic) RWRE exit measure approaches
that of
simple random walk.
\end{remark}

\section{Green's functions for the ball}
\label{super}
One main task of our approach aims at developing good estimates on Green's
functions for the ball of both coarse grained (goodified) RWRE as well
as coarse grained
symmetric random walk in the perturbative regime. The main result is
Lemma~\ref{superlemma}. For the coarse grained symmetric random walk,
the estimates on
hitting probabilities of the last section together with
Proposition~\ref{super-behaviorgreen} yield the right control.

On a certain class of environments, we need to modify
the transition kernels in order to ensure that bad
points are not visited too often by the coarse grained random walks. This
modification will be described in Section~\ref{super-cgmod}.

We work with the same convention concerning the parameter $\kappa$ as
in Section~\ref{estimates}.
\subsection{A local central limit theorem}
\label{clt}
Let $p\in\mathcal{P}^{\mathrm{s}}_\kappa$ and $m\geq1$. Denote
by ${\hat{\pi
}}_{\psi_m}=\hat{\pi}^{(p)}_{\psi_m}$
the coarse grained transition probabilities on $\mathbb{Z}^d$
associated to
the constant function $\psi_m \equiv m$; cf.~(\ref{prel-cgpsi}). We
constantly drop $p$ from notation. Notice that ${\hat{\pi
}}_{\psi_m}$ is centered, and
the covariances satisfy
\[
\sum_{y\in\mathbb{Z}^d}(y_i-x_i)
(y_j-x_j){\hat{\pi}}_{\psi_m}(x,y) =
\lambda_{m,i}\delta_i(j), %
\]
where for large $m$,
$C^{-1}<\lambda_{m,i}/m^2< C$ for some $C>0$.
Define the matrix
\[
\Lambda_m= \bigl(\lambda_{m,i}\delta_i(j)
\bigr)_{i,j=1}^d, %
\]
and let for $x\in\mathbb{Z}^d$
\[
\mJ(x) = \bigl\llvert\Lambda_m^{-1/2}x\bigr\rrvert.
\]
%
%
\begin{proposition}[(Local central limit theorem)]\label{super-localclt}
Let $p\in\mathcal{P}^{\mathrm{s}}_\kappa$, and let $x,y\in
\mathbb{Z}^d$. For $m \geq1$ and
all integers $n\geq1$,
\[
({\hat{\pi}}_{\psi_m} )^n(x,y) = \frac
{1}{(2\pi
n)^{d/2}\det\Lambda_m^{1/2}}\exp
\biggl(-\frac{\mJ
^2(x-y)}{2n} \biggr) + O \bigl(m^{-d}n^{-(d+2)/2}
\bigr). %
\]
\end{proposition}
For the corresponding Green's function ${\hat
{g}}_{m,{\mathbb{Z}^d}}(x,y) =
\sum_{n=0}^\infty({\hat{\pi}}_{\psi_m})^n(x,y)$ we
obtain the following:
%
\begin{proposition}
\label{super-behaviorgreen}
Let $p\in\mathcal{P}^{\mathrm{s}}_\kappa$. Let $x,y\in\mathbb
{Z}^d$, and assume $m \geq m_0
> 0$ large enough.
\begin{longlist}[(ii)]
\item For $\llvert x-y\rrvert < 3m$,
\[
\hat{g}_{m,{\mathbb{Z}^d}}(x,y) = \delta_x(y) + O
\bigl(m^{-d}\bigr). %
\]
\item For $\llvert x-y\rrvert \geq3m$, there exists a constant $c(d)
> 0$ such that
\[
\hat{g}_{m,{\mathbb{Z}^d}}(x,y)= \frac{c(d)\det
\Lambda_m^{-1/2}}{\mJ(x-y)^{d-2}} + O \biggl(\frac{1}{\llvert
x-y\rrvert ^{d}}
\biggl(\log\frac{\llvert x-y\rrvert }{m} \biggr)^d \biggr). %
\]
\end{longlist}
\end{proposition}
Note that the constants in the $O$-notation are independent of $n$, $m$
and $\llvert x-y\rrvert $.

In our applications, $m$ will be a function of $L$. Although these results
look rather standard, we cannot directly refer to the literature
because we
have to keep track of the $m$-dependency. We give a proof of both
statements in the \hyperref[appe]{Appendix}. The last proposition will be used to estimate
the Green's function for the ball $V_L$, ${\hat
{g}}_{m,V_L}(x,y) =
\sum_{n=0}^\infty(1_{V_L}{\hat{\pi}}_{\psi
_m} )^n(x,y)$. Clearly, $\hat{g}_{m,V_L}$ is
bounded from above by~$\hat{g}_{m,{\mathbb{Z}^d}}$.

\subsection{Estimates on coarse grained Green's functions}
\label{super-Ghestimates}
As we will show, the perturbation expansion enables us to control the
goodified Green's function ${\hat{G}^g}_{L,r}$ essentially in terms of
${\hat{g}}^{(p)}_{L,r}$, where $p$ is the kernel corresponding to the
radius $s_L/20$,
stemming from assignment~(\ref{kernelpL}).

The first step in controlling the Green's function is
provided by the following lemma.

%
\begin{lemma}
\label{phgp-php-estimate}
Assume $\ctwo(\delta,L_0,L_1)$, let $L_1\leq L\leq L_1(\log L_1)^2$,
and put
$p=p_{s_L/20}$. Then for all $x\in V_L\setminus\sh_L(2r)$, with
$H(x) =
\max\{L_0, h_{L,r}(x)\}$,
\begin{eqnarray*}
&& \bigl\llVert\bigl({\hat{\Pi}^g}_{L,r} - {\hat{
\pi}^{(p)}}_{L,r} \bigr)\hat{\pi}^{(p)}_{L,r}(x,
\cdot)\bigr\rrVert_1
\\
&&\qquad \leq C\min\bigl\{\log\bigl(s_L/H(x)\bigr) \bigl(\log H(x)
\bigr)^{-9}, \bigl(\log H(x) \bigr)^{-8} \bigr\}
\end{eqnarray*}
and
\[
\bigl\llVert\bigl({\hat{\Pi}^g}_{L,r} - {\hat{
\pi}^{(p)}}_{L,r} \bigr) (x,\cdot)\bigr\rrVert_1
\leq2\delta.
\]
\end{lemma}
\begin{pf}
For $x\in\mathcal{B}^\star_{L,r}$, both left-hand sides are zero.
Now let $x\in
V_L\setminus(\sh_L(2r)\cup\mathcal{B}^\star_{L,r})$ and set
$q=p_{h_{L,r}(x)}$. By
the triangle inequality,
\begin{eqnarray*}
&& \bigl\llVert\bigl({\hat{\Pi}^g}- {\hat{
\pi}^{(p)}} \bigr)\hat{\pi}^{(p)}(x,\cdot)\bigr\rrVert
_1
\\
&&\qquad \leq\bigl\llVert\bigl({\hat{\Pi}^g}- {\hat{\pi}^{(q)}}
\bigr)\hat{\pi}^{(p)}(x,\cdot)\bigr\rrVert_1 + \bigl
\llVert\bigl(\hat{\pi}^{(p)}-{\hat{\pi}^{(q)}} \bigr) (x,
\cdot)\bigr\rrVert_1
\\
&&\qquad \leq\bigl\llVert\bigl({\hat{\Pi}^g}- {\hat{\pi}^{(q)}}
\bigr){\hat{\pi}^{(q)}}(x,\cdot)\bigr\rrVert_1 + 2\sup
_{y\in V_L\setminus\sh(r)}\bigl\llVert\bigl(\hat{\pi}^{(p)}-{\hat{
\pi}^{(q)}} \bigr) (y,\cdot)\bigr\rrVert_1
\\
&&\qquad \leq C \bigl(\bigl(\log H(x)\bigr)^{-9} + \llVert p-q\rrVert
_1 \bigr),
\end{eqnarray*}
where in the last line we used that $x$ is good and
Lemma~\ref{php-phq-estimate}. Now, with $K=\lfloor
\log_2 (s_L/H(x))\rfloor$, Lemma~\ref{diff-admiss} shows
\begin{eqnarray*}
\llVert p-q\rrVert_1&\leq& C\sum_{i=1}^K
\bigl(\log\bigl(2^{-i}s_L\bigr)\bigr)^{-9}
\\
&\leq& C\min\bigl\{K\bigl(\log H(x)\bigr)^{-9}, \bigl(\log H(x)
\bigr)^{-8} \bigr\}.
\end{eqnarray*}
This proves the claim for the smoothed difference. For the nonsmoothed
difference,
\[
\bigl\llVert\bigl({\hat{\Pi}^g}-\hat{\pi}^{(p)} \bigr)
(x,\cdot)\bigr\rrVert_1\leq\bigl\llVert\bigl({\hat{
\Pi}^g}-{\hat{\pi}^{(q)}} \bigr) (x,\cdot)\bigr\rrVert
_1 + \bigl\llVert\bigl(\hat{\pi}^{(p)}-{\hat{\pi}^{(q)}} \bigr) (x,\cdot)\bigr\rrVert_1. %
\]
Since $x$ is good, the first term is bounded by $\delta$, and, by what
we have
just seen, the second term is
bounded by $\delta$ as well if we choose $L_0$ (and so $L$) large enough.
\end{pf}
%
%
\begin{remark}
\label{super-remr}
Notice that the choice of the parameter $r$ depends on
$\delta$. See also the preliminary remarks of Section~\ref{nonsmv-exits}.
\end{remark}
Recall that in the goodifying-procedure introduced in
Section~\ref{smoothbad}, ``bad'' exit distributions inside $V_L$ are
replaced by such of a symmetric random walk with one-step distribution
$p=p_{s_L/20}$. For this $p$ and good points $x$ within the boundary region
$\sh_L(2r)$, we would like to use at least an estimate of the form
\[
\bigl\llVert\bigl(\Ph_{L,r} - \hat{\pi}^{(p)}_{L,r}
\bigr) (x,\cdot)\bigr\rrVert_1 \leq C\delta. %
\]
However, exit measures at points $x$ inside $\sh_L((1/10)r)$ are taken from
intersected balls $V_t(x)\cap V_L$.
We therefore work in this (and only in this) section with slightly modified
transition kernels $\Pt_{L,r}$, ${\tilde{\pi}}_{L,r}$,
${\tilde{\Pi}^g}_{L,r}$ in\vspace*{1pt} the enlarged
ball $V_{L+r}$, taking the exit measure in $\sh_L(2r)$ from uncut balls
$V_t(x)\subset V_{L+r}$, $t\in[h_{L,r}(x),2h_{L,r}(x)]$.

Now, to make things precise, for $q\in{\mathcal{P}^{\mathrm
{s}}_\kappa}$, we set $h_{L,r}(x) =
(1/20)r$ for $x\notin C_L$, and let ${\tilde{\pi
}^{(q)}}_{L,r}$ be the
coarse grained symmetric random walk kernel in $V_{L+r}$ associated to
$\tilde{\psi} =
(h_{L,r}(x))_{x\in V_{L+r}}$,
\[
{\tilde{\pi}^{(q)}}_{L,r}(x,\cdot) = \frac{1}{h_{L,r}(x)}\int
_{\mathbb{R}_+} \varphi\biggl(\frac{t}{h_{L,r}(x)} \biggr)\pi^{(q)}_{V_t(x)\cap V_{L+r}}(x,\cdot)\,\dt t. %
\]
For the corresponding RWRE kernel, we forget about the environment on
$V_{L+r}\setminus V_L$ and set
\[
{\tilde{\Pi}^{(q)}}_{L,r}(x,\cdot) = \cases{ \displaystyle
\frac{1}{h_{L,r}(x)}\int_{\mathbb{R}_+}\varphi\biggl(\frac
{t}{h_{L,r}(x)}
\biggr)\Pi_{V_t(x)}(x,\cdot)\,\dt t, &\quad for $x\in V_L$,
\vspace*{5pt}\cr
\displaystyle {\tilde{\pi}^{(q)}}_{L,r}(x,\cdot), &\quad for $x\in
V_{L+r}\setminus V_L$.} %
\]
For\vspace*{1.5pt} $p=p_{s_L/20}$ and all good $x\in V_L$, we now have $\llVert
({\tilde
{\Pi}^{(p)}}
_{L,r}-{\tilde{\pi}^{(p)}}_{L,r})(x,\cdot)\rrVert _1
\leq\delta$ provided $\varepsilon$ is small enough,
while for $x\in V_{L+r}\setminus V_L$, the difference vanishes anyway.
The goodified version of ${\tilde{\Pi}^{(p)}}_{L,r}$ is then
obtained in an analogous way to~(\ref{smoothbad-goodifiedkernel}),
\[
{\tilde{\Pi}^g}_{L,r}(x,\cdot) = \cases{ {\tilde{
\Pi}^{(p)}}_{L,r}(x,\cdot), &\quad for $x\notin\mathcal {B}^\star_{L,r}$,
\vspace*{2pt}\cr
{\tilde{\pi}^{(p)}}_{L,r}(x,
\cdot), &\quad for $x\in\mathcal{B}^\star_{L,r}$.}
\]
Clearly, for $x\in V_L\setminus\sh_L(2r)$, the first statement of
Lemma~\ref{phgp-php-estimate} holds with the left-hand side there
replaced by
\[
\bigl\llVert\bigl({\tilde{\Pi}^g}_{L,r} - {\tilde{\pi
}^{(p)}}_{L,r} \bigr){\tilde{\pi}^{(p)}}_{L,r}(x,
\cdot)\bigr\rrVert_1. %
\]
But thanks to the modified transition kernels, we now have
\begin{eqnarray*}
\bigl\llVert\bigl({\tilde{\Pi}^g}_{L,r} - {\tilde{\pi
}^{(p)}}_{L,r} \bigr) (x,\cdot)\bigr\rrVert_1 &\leq& 2\delta
\end{eqnarray*}
for \textit{all} $x\in V_L$. Indeed, one just has to notice that
Lemma~\ref{phgp-php-estimate} can now also be applied to points
$x\in\sh_L(2r)$, with the same proof.

We write $\Gt_{L,r},{\tilde{g}}_{L,r}$ and ${\tilde{G}^g}_{L,r}$
for the
Green's functions on $V_{L+r}$ corresponding to $\Pt_{L,r}$,
${\tilde{\pi}}_{L,r}$
and ${\tilde{\Pi}^g}_{L,r}$. Note
%
%
\begin{eqnarray}\label{super-greensfcts-pointwisebdd}
\Gh_{L,r}\leq\Gt_{L,r},\qquad{
\hat{g}}_{L,r} \leq{\tilde{g}}_{L,r},\qquad{
\hat{G}^g}_{L,r} \leq{\tilde{G}^g}_{L,r}
\nonumber\\[-8pt]\\[-8pt]
\eqntext{\mbox{pointwise on }V_{L+r}\times(V_{L+r}\setminus
\partial V_L).}
\end{eqnarray}
Since we do not have exact expressions for ${\tilde
{g}}_{L,r}$ or $\Gt_{L,r}$, we
construct a (deterministic) kernel $\Gamma_{L,r}$ that bounds the Green's
functions from above. For $x\in V_{L+r}$, set
\[
\dtL(x) = \max\biggl(\frac{\dist_{L+r}(x)}{2},3r \biggr),\qquad a(x)
=\min\bigl(
\dtL(x),s_L \bigr). %
\]
Furthermore, for $x,y\in V_{L+r}$, let
\begin{eqnarray*}
\Gamma^{(1)}_{L,r}(x,y) &=& \frac{\dtL(x)\dtL(y)}{a(y)^2(a(y) +
\llvert x-y\rrvert )^d},
\\
\Gamma^{(2)}_{L,r}(x,y) &=& \frac{1}{a(y)^2(a(y) + \llvert x-y\rrvert
)^{d-2}}. %
\end{eqnarray*}
The kernel $\Gamma_{L,r}$ is defined as the pointwise minimum
%
%
\begin{equation}
\label{defgamma} \Gamma_{L,r} = \min\bigl\{\Gamma^{(1)}_{L,r},
\Gamma^{(2)}_{L,r} \bigr\}.
\end{equation}
We cannot derive pointwise estimates on the Green's functions in terms of
$\Gamma_{L,r}$, but we can use this kernel to obtain upper bounds on
neighborhoods $U(x) = V_{a(x)}(x)\cap V_{L+r}$. Call a function $F\dvtx
V_{L+r}\times V_{L+r}\rightarrow\mathbb{R}_+$ a \textit{positive
kernel}. Given two positive kernels $F$ and $G$, we write $F\preceq G$ if
for all $x,y\in V_{L+r}$,
\[
F\bigl(x,U(y)\bigr) \leq G\bigl(x,U(y)\bigr), %
\]
where $F(x,U)$ stands for $\sum_{y\in U\cap\mathbb{Z}^d}F(x,y)$. We write
$F\asymp1$, if there is a constant $C>0$ such that for all $x,y\in
V_{L+r}$,
\[
\frac{1}{C}F(x,y)\leq F(\cdot,\cdot)\leq CF(x,y)\qquad\mbox{on }U(x)\times
U(y). %
\]
We adapt this notation to positive functions of one argument: for
$f\dvtx V_{L+r}\rightarrow\mathbb{R}_+$, $f\asymp1$ means that for some
$C>0$, $C^{-1}f(x)\leq
f(\cdot)\leq Cf(x)$ on any $U(x)\subset V_{L+r}$.
Finally, given $0 <\eta< 1$, we say that a positive kernel $A$ on
$V_{L+r}$ is
\mbox{$\eta$-}\textit{smoothing}, if for all $x\in V_{L+r}$, $A(x,U(x)) \leq
\eta$, and
$A(x,y) = 0$ whenever $y\notin U(x)$.

Now we are in position to formulate our main statement of this
section. Recall our convention concerning constants: they only depend on
the dimension unless stated otherwise.
%
\begin{lemma} \label{superlemma}
\textup{(i)} There exists a constant $C_1>0$ such that for all $q\in{\mathcal
{P}^{\mathrm{s}}_\kappa}$,
\[
{\hat{g}}^{(q)}_{L,r}\preceq C_1
\Gamma_{L,r}\quad\mbox{and}\quad{\tilde{g}}^{(q)}_{L,r}
\preceq C_1\Gamma_{L,r}. %
\]

\textup{(ii)}  Assume $\ctwo(\delta,L_0,L_1)$, and let $L_1\leq L\leq L_1(\log
L_1)^2$. There exists a constant $C>0$ such that for $\delta>0$ small,
\[
{\hat{G}^g}_{L,r} \preceq C\Gamma_{L,r}\quad
\mbox{and}\quad{\tilde{G}^g}_{L,r} \preceq C
\Gamma_{L,r}. %
\]
\end{lemma}
%
%
\begin{remark}
(i) Thanks to~(\ref{super-greensfcts-pointwisebdd}), it suffices to show
the bounds for ${\tilde{g}}_{L,r}$ and~${\tilde{G}^g}_{L,r}$. For
later use, we keep track of
the constant in part~(i) of the lemma.\vspace*{1pt}

(ii) We will later apply part (i) with $q=p_L$. From Lemma~\ref{diff-admiss}
we know that we can assume $p_L\in{\mathcal{P}^{\mathrm{s}}_\kappa
}$ for every choice of $\kappa>0$,
if $L_0$ is large.
\end{remark}
We first prove part~(i), which will be a straightforward consequence of
the estimates on hitting
probabilities in Section~\ref{estimates} and the next lemma.
%
\begin{lemma}
\label{super-greenbound}
There exists a constant $C>0$ such
that for all $q\in\mathcal{P}^{\mathrm{s}}_\kappa$, for all $x\in
V_{L+r}$ and $y \in V_{L}$
with $\dL(y) \geq4s_L$,
\[
{\tilde{g}}^{(q)}_{L,r}(x,y)\leq C \cases{
s_L^{-2}\max\bigl\{\llvert x-y\rrvert,s_L
\bigr\}^{-(d-2)}, &\quad for $y\neq x$,
\vspace*{2pt}\cr
1, &\quad for $y=x$.} %
\]
\end{lemma}
\begin{pf}
The underlying one-step transition kernel is always given by
$q\in\mathcal{P}^{\mathrm{s}}_\kappa$, which we therefore omit
from notation. For example,
${\tilde{g}}={\tilde{g}}^{(q)}$, $\hat{g}_{m,V_L} =
{\hat{g}}^{(q)}_{m,V_L}$, $\Prw_{x} = \Prw_{x,q}$.

If $x=y$, then the claim follows from transience of simple random
walk. Now assume $x\neq y$, and always $\dL(y) \geq4s_L$. Consider first
the case $\llvert x-y\rrvert \leq s_L$. Recall that $\hat{g}_{m,V_L}$
denotes the Green's
function for the ball $V_L$ associated to ${\hat{\pi
}}_{\psi_m}$, where
$\psi_m\equiv m$. With $m= s_L/20$ we have
\[
{\tilde{g}}(x,y) \leq\hat{g}_{m,V_L}(x,y) + \sup_{v\in\sh_L(2s_L)}
\Prw_v (T_{V_{s_L }(y)} < \tau_{V_{L+r}} ) \mathop{\sup_{w\dvtx  w\neq y,}}_{\llvert w-y \rrvert \leq
s_L}{\tilde{g}}(w,y). %
\]
Since
\[
\sup_{v\in\sh_L(2s_L)}\Prw_v (T_{V_{s_L}(y)} < \tau
_{V_{L+r}} ) < 1 %
\]
uniformly in $L$, it follows from Proposition~\ref{super-behaviorgreen}
that
\[
{\tilde{g}}(x,y) \leq C \mathop{\sup_{w\dvtx  w\neq y,}}_{\llvert w-y \rrvert \leq
s_L} \hat{g}_{m,V_L}(w,y)
\leq C \mathop{\sup_{w\dvtx  w\neq y,}}_{\llvert w-y \rrvert \leq s_L} \hat{g}_{m,\mathbb{Z}^d}(w,y) \leq
\frac{C}{s_L^d}. %
\]
If $\llvert x-y\rrvert > s_L$ we use Lemma~\ref{hittingprob}(i) and
the first
case to get
\[
{\tilde{g}}(x,y) \leq\Prw_x (T_{V_{s_L}(y)} < \infty)
\mathop{\sup_{w\dvtx  w\neq y,}}_{\llvert w-y \rrvert \leq s_L}{\tilde{g}}(w,y) \leq\frac
{C}{s_L^2\llvert x-y\rrvert ^{d-2}}. %
\]\upqed
\end{pf}
\begin{pf*}{Proof of Lemma~\ref{superlemma}(\textup{i})}
It suffices to prove the bound
for ${\tilde{g}}$. First we show that there exists a
constant $C > 0$ such that
for all $y \in V_{L+r}$,
%
%
\begin{equation}
\label{super-keyest-g1} \sup_{x\in V_{L+r}}{\tilde{g}} \bigl(x,U(y)
\bigr) \leq
C.
\end{equation}
At first let $\dist_{L+r}(y) \leq6r$. Then $U(y)\subset
\sh_{L+r}(10r)$. We claim that even
%
%
\begin{equation}
\label{super-keyest-g2} \sup_{x\in V_{L+r}}{\tilde{g}} \bigl(x,\sh
_{L+r}(10r) \bigr)\leq C
\end{equation}
for some $C> 0$. Indeed, if $z\in\sh_{L+r}(10r)$, then $
{\tilde{\pi}}(z,\cdot)$ is
an (averaging) exit distribution from balls $V_{\ell}(z)\cap V_{L+r}$, where
$\ell\geq r/20$. Using Lemma~\ref{lemmalawler}(i), we find a
constant $k_1 = k_1(d)$ such that starting at any $z\in\sh_{L+r}(10r)$,
$V_{L+r}$ is left after $k_1$ steps with probability $ > 0$, uniformly in
$z$. This together with the strong Markov property\vspace*{1pt}
implies~(\ref{super-keyest-g2}). Next assume $6r<\dist_{L+r}(y)\leq
6s_L$. Then $U(y) \subset
S(y)=\sh_{L+r} (\frac{1}{2}\dist_{L+r}(y),2\dist
_{L+r}(y) )$. We
claim that
%
%
\begin{equation}
\label{super-keyest-g3} \sup_{x\in V_{L+r}}{\tilde{g}} \bigl(x,S(y)
\bigr)\leq
C.
\end{equation}
For $z\in S(y)$, ${\tilde{\pi}}(z,\cdot)$ is an
averaging exit distribution from
balls $V_l(z)$, where $l\geq\dist_{L+r}(y)/240$. By
Lemma~\ref{lemmalawler}(i), we find some small $0 < c
< 1$ and a constant $k_2(c,d)$
such that after $k_2$ steps, the walk has probability $> 0$ to be in
$\sh_{L+r} (\frac{1-c}{2}\dist_{L+r}(y) )$, uniformly in
$z$ and $y$. But starting in
$\sh_{L+r} (\frac{1-c}{2}\dist_{L+r}(y) )$, an iterative
application of Lemma~\ref{lemmalawler}(i) shows that with probability
$> 0$, the
ball $V_{L+r}$ is left before $S(y)$ is visited
again. Therefore~(\ref{super-keyest-g3}) and
hence~(\ref{super-keyest-g1}) hold in this case.
At last, let $\dist_{L+r}(y) > 6s_L$. Then $\dL(w) \geq4s_L$ for
$w\in U(y)$. Estimating
\[
{\tilde{g}}(x,w) \leq1 + \sup_{v\dvtx  v\neq w} {\tilde{g}}(v,w),
\]
we get with part (i) that
\[
\sup_{w\in U(y)}{\tilde{g}}(x,w) \leq1 + \frac{C}{s_L^d}.
\]
Summing over $w\in U(y)$, (\ref{super-keyest-g1}) follows. Finally, note
that for any $x\in V_{L+r}$,
\[
{\tilde{g}}\bigl(x,U(y)\bigr) \leq\Prw_x (T_{U(y)} < \tau
_{V_{L+r}} )\sup_{w\in U(y)}{\tilde{g}}\bigl(w,U(y)\bigr).
\]
Now ${\tilde{g}}\preceq C\Gamma$ follows from (\ref
{super-keyest-g1}) and
the hitting estimates of
Lemma~\ref{hittingprob}.
\end{pf*}
Let\vspace*{1pt} us now explain our strategy for proving part~(ii). By
version~(\ref{prel-pbe3}) of the perturbation expansion, we can express
${\tilde{G}^g}_{L,r}$ in a series involving ${\tilde{g}}_{L,r}$
and differences of
exit measures. The Green's
function ${\tilde{g}}_{L,r}$ is already controlled by
means of $\Gamma_{L,r}$.
Looking at~(\ref{prel-pbe3}), we thus have to understand what happens if
$\Gamma_{L,r}$ is concatenated with certain smoothing kernels. This
will be
the content of Proposition~\ref{super-concatenating}.

We start with collecting some important properties of $\Gamma_{L,r}$,
which will be
used throughout this text. Define for $j\in\mathbb{N}$,
\begin{eqnarray*}
\mathcal{L}_j &=& \bigl\{y\in V_L\dvtx  j\leq\dL(y) < j+1
\bigr\},
\\
\mathcal{E}_j&=&\bigl\{y\in V_{L+r}\dvtx  \dtL(y)\leq3jr\bigr\}. %
\end{eqnarray*}

%
%
\begin{lemma}[(Properties of $\Gamma_{L,r}$)]\label{super-gammalemma}
\textup{(i)}~Both $\dtL$ and $a$ are Lipschitz with constant $1/2$. Moreover,
for $x,y \in V_{L+r}$,
\[
a(y) + \llvert x-y\rrvert\leq a(x) + \tfrac{3}{2}\llvert x-y\rrvert.
\]

\textup{(ii)}
\[
\Gamma_{L,r}\asymp1. %
\]

\textup{(iii)} For $0\leq j \leq2s_L$, $x\in V_{L+r}$,
\[
\sum_{y\in\mathcal{L}_j} \biggl(\max\biggl\{1,\frac{\dtL
(x)}{a(y)}
\biggr\}\frac{1}{(a(y)+
\llvert x-y\rrvert )^d} \biggr) \leq C\frac{1}{j\vee r}. %
\]

\textup{(iv)} For $1\leq j \leq\frac{1}{3r}s_L$,
\[
\sup_{x\in V_{L+r}}\Gamma_{L,r}(x,\mathcal{E}_j)
\leq C\log(j+1), %
\]
and for $0\leq\alpha< 3$,
\[
\sup_{x\in V_{L+r}}\Gamma_{L,r} \bigl(x,\sh_L
\bigl(s_L, L/(\log L)^\alpha\bigr) \bigr) \leq C(\log\log
L) (\log L)^{6-2\alpha}. %
\]

\textup{(v)} For $x\in V_{L+r}$, in the case of constant $r$,
\[
\Gamma_{L,r}(x,V_L) \leq C\max\biggl\{
\frac{\dtL(x)}{L}(\log L)^6, \biggl(\frac{\dtL(x)}{r} \wedge\log L
\biggr) \biggr\}. %
\]
In the case $r= r_L$,
\[
\Gamma_{L,r_L}(x,V_L) \leq C\max\biggl\{
\frac{\dtL(x)}{L}(\log L)^6, \biggl(\frac{\dtL(x)}{r_L}\wedge\log\log
L \biggr) \biggr\}. %
\]
\end{lemma}

\begin{pf}
(i) The second statement is a direct consequence of the Lipschitz property,
which in turn follows immediately from the definitions of $\dtL$ and
$a$.

(ii) As for $y'\in U(y)$, $\frac{1}{2}a(y)\leq a(y') \leq
\frac{3}{2}a(y)$ and similarly with $a$ replaced by~$\dtL$, it suffices
to show that for $x'\in U(x)$, $y'\in U(y)$,
%
%
\begin{equation}
\label{super-gammaest-ineq} \frac{1}{C} \bigl(a(y)+\llvert x-y\rrvert
\bigr) \leq a
\bigl(y'\bigr)+\bigl\llvert x'-y' \bigr
\rrvert\leq C \bigl(a(y)+\llvert x-y\rrvert\bigr).
\end{equation}
First consider the case $\llvert x-y\rrvert \geq4\max\{a(x), a(y)\}$. Then
\[
a(y) +\llvert x-y\rrvert\leq2a\bigl(y'\bigr) + 2 \bigl(\llvert
x-y\rrvert- a(x) - a(y) \bigr) \leq2 \bigl(a\bigl(y'\bigr) +
\bigl\llvert x'-y' \bigr\rrvert\bigr). %
\]
If $\llvert x-y\rrvert \leq4 a(y)$, then
\[
a(y) +\llvert x-y\rrvert\leq5a(y) \leq5a(y) + \bigl\llvert x'-y'
\bigr\rrvert\leq10 \bigl(a\bigl(y'\bigr) + \bigl\llvert
x'-y' \bigr\rrvert\bigr), %
\]
while for $\llvert x-y\rrvert \leq4 a(x)$, using part~(i) in the
first inequality,
\[
a(y) +\llvert x-y\rrvert\leq a(x)+\tfrac{3}{2}\llvert x-y\rrvert
\leq7a(x) \leq14 \bigl(a\bigl(y'\bigr) + \bigl\llvert
x'-y' \bigr\rrvert\bigr). %
\]
This proves the first inequality in~(\ref{super-gammaest-ineq}). The
second one follows from
\[
a\bigl(y'\bigr) + \bigl\llvert x'-y'
\bigr\rrvert\leq\tfrac{5}{2}a(y) + a(x) + \llvert x-y\rrvert\leq
\tfrac
{7}{2} \bigl(a(y) + \llvert x-y\rrvert\bigr). %
\]

(iii) If $j\leq2s_L$ and $y\in\mathcal{L}_j$, then $a(y)$ is of order
$j\vee r$. By Lemma~\ref{hittingprob-technical} we have
\[
\sum_{y\in\mathcal{L}_j}\frac{1}{(j\vee r + \llvert x-y\rrvert
)^d}\leq C\min\biggl\{
\frac{1}{j\vee r},\frac{1}{\llvert \dLk(x)-(j+r)\rrvert } \biggr\}. %
\]
It remains to show that
%
%
\begin{equation}
\label{super-gammalemma-layerbound2} \max\biggl\{1,\frac{\dtL
(x)}{j\vee r} \biggr\}\min\biggl\{
\frac
{1}{j\vee r}, \frac{1}{\llvert \dLk(x)-(j+r)\rrvert } \biggr\} \leq
C\frac{1}{j\vee r}.
\end{equation}
If\vspace*{1.5pt} $\dtL(x)\leq(j\vee3r)$, this is clear. If $\dtL(x) > (j\vee3r)$,
(\ref{super-gammalemma-layerbound2}) follows from
$\llvert \dLk(x)-(j+r)\rrvert \geq\dtL(x)/2$.

(iv) We follow our convention and write $\Gamma$ instead of
$\Gamma_{L,r}$. If $\dtL(y) \leq3jr$, then $\dL(y) \leq6jr$. Estimating
$\Gamma$ by $\Gamma^{(1)}$, we get
\[
\Gamma(x,\mathcal{E}_j) \leq C\sum_{i=0}^{6jr}
\sum_{y\in\mathcal{L}_i}\frac{\dtL(x)}{a(y)}\frac{1}{(a(y) + \llvert
x-y\rrvert )^d}.
\]
Now the first assertion of (iv) follows from (iii). The second is proved
similarly, so we omit the details.

(v) Set\vspace*{1pt} $B= \{y\in V_L\dvtx  \dtL(y) \leq s_L\vee2\dtL(x)\}$. For $y\in
V_L\setminus B$, it holds that $a(y) = s_L$ and $\llvert x-y\rrvert
\geq
\dtL(y)$. Therefore,
\[
\Gamma(x, V_L\setminus B ) \leq\Gamma^{(1)} (x,
V_L\setminus B )\leq\frac{\dtL(x)}{s_L^2}\sum
_{y\in V_{2L}}\frac{1}{(s_L +
\llvert y \rrvert )^{d-1}}
\leq C\frac{\dtL(x)}{L}(\log
L)^6. %
\]
Furthermore,
\begin{eqnarray*}
\Gamma(x, B )
&\leq&\sum_{i=0}^{2s_L} \sum
_{y\in\mathcal{L}_i}\frac{\dtL(x)}{a(y)}\frac{1}{(a(y) +
\llvert x-y\rrvert )^d} +
\frac{1}{s_L^2}
\mathop{\sum_{y\in V_L\dvtx }}_{s_L\leq\dtL(y)\leq2\dtL
(x)}\frac{1}{(s_L+\llvert x-y\rrvert )^{d-2}}.
\end{eqnarray*}
Lemma~\ref{hittingprob-technical} bounds the second term by
$C(\dtL(x)/L)(\log L)^6$. For the first term, we use twice part~(iii) and
once Lemma~\ref{hittingprob-technical} to get
\[
\sum_{i=0}^{2s_L}\sum
_{y\in\mathcal{L}_i}\frac{\dtL
(x)}{a(y)}\frac{1}{(a(y) + \llvert x-y\rrvert )^d} \leq C\sum
_{i=0}^{5r}\frac{1}{i\vee r} + C\min\Biggl\{
\dtL(x)\sum_{i=5r}^{2s_L}\frac{1}{i^2},
\sum_{i=5r}^{2s_L}\frac
{1}{i} \Biggr\}.
\]
This proves (v).
\end{pf}

%
%
\begin{proposition}\label{super-concatenating}
Let $F, G$ be positive kernels with \mbox{$F\preceq G$}.
\begin{longlist}[(ii)]
\item If $A$ is $\eta$-smoothing and $G\asymp1$, then for some
constant $C = C(d,G) > 0$,
\[
FA \preceq C\eta G. %
\]
\item If $\Phi$ is a positive function on $V_{L+r}$ with $\Phi\asymp
1$, then for some $C = C(d,\Phi) > 0$,
\[
F\Phi\leq CG\Phi. %
\]
\end{longlist}
\end{proposition}

\begin{pf}
(i) Let $y\in V_{L+r}$. As $a$ is Lipschitz with constant $1/2$, we can
choose $K=K(d)$
points $y_k$ out of the set $M=\{y'\in V_{L+r}\dvtx  U(y')\cap U(y) \neq
\varnothing\}$ such that $M$ is covered by the union of the $U(y_k), k =
1,\ldots,K$. Since $A(y',U(y))\neq0$ implies $y'\in M$, we then have
\begin{eqnarray*}
FA\bigl(x,U(y)\bigr) &=& \sum_{y' \in M}F
\bigl(x,y'\bigr) \sum_{y''\in U(y)}A
\bigl(y',y''\bigr) \leq\eta\sum
_{k=1}^KF\bigl(x,U(y_k)\bigr)
\\
&\leq&\eta\sum_{k=1}^KG
\bigl(x,U(y_k)\bigr).
\end{eqnarray*}
Using $G \asymp1$, we get $G(x,U(y_k))\leq C\llvert U(y_k)\rrvert
G(x,y)$. Clearly
$\llvert U(y_k)\rrvert \leq C \llvert U(y)\rrvert $, so that
\[
FA\bigl(x,U(y)\bigr) \leq CK\eta\bigl\llvert U(y)\bigr\rrvert G(x,y).
\]
A second application of $G\asymp1$ yields the claim.

(ii) We can find a constant $K = K(d)$ and a covering of $V_{L+r}$ by
neighborhoods $U(y_k)$, $y_k \in V_{L+r}$, such that every
$y\in V_{L+r}$ is contained in at most $K$ many of the sets $U(y_k)$. Using
$\Phi\asymp1$, it follows that for $x\in V_{L+r}$,
\begin{eqnarray*}
F\Phi(x) &=& \sum_{y\in V_{L+r}}F(x,y)\Phi(y) \leq C\sum
_{k=1}^{\infty}F\bigl(x,U(y_k)\bigr)
\Phi(y_k)
\\
&\leq& C\sum_{k=1}^{\infty}G
\bigl(x,U(y_k)\bigr)\Phi(y_k) \leq C\sum
_{k=1}^{\infty}\sum_{y\in U(y_k)}G(x,y)
\Phi(y)
\\
&\leq& CK\sum_{y\in
V_{L+r}}G(x,y)\Phi(y).
\end{eqnarray*}\upqed
\end{pf}
In terms of our specific kernel $\Gamma_{L,r}$, we obtain the following:
%
\begin{proposition} \label{super-keyest}
Let $A$ be $\eta$-smoothing, and let $F$ be a positive kernel satisfying
$F\preceq\Gamma_{L,r}$.
\begin{longlist}[(ii)]
\item There exists a constant $C_2 > 0$ not depending on $F$ such that
\[
FA \preceq C_2\eta\Gamma_{L,r}. %
\]
\item If additionally $A(x,y) = 0$ for $x\notin V_L$ and $A(x,U(x))
\leq
(\log a(x) )^{-15/2}$ for $x \in V_L\setminus\mathcal
{E}_1$, then
there exists a constant $C_3 > 0$ not depending on $F$ such that for
all $x,z\in V_{L+r}$,
\[
FA\Gamma_{L,r}(x,z) \leq C_3\eta^{1/2}
\Gamma_{L,r}(x,z). %
\]
\end{longlist}
\end{proposition}
\begin{pf}
(i) This is Proposition~\ref{super-concatenating}(i) with $G =
\Gamma$.

(ii) We set $B = V_L\setminus\mathcal{E}_1$ and split into
%
%
\begin{equation}
\label{super-keyest-splitting1} FA\Gamma= F1_{\mathcal{E}_1}A\Gamma+
F1_BA\Gamma.
\end{equation}
Let $x,z\in V_{L+r}$ be fixed, and consider first
$F1_{\mathcal{E}_1}A\Gamma(x,z)$. Using $\Gamma\asymp1$, $A\Gamma
(y,z) \leq
C\eta\Gamma(y,z)$. As $\Gamma(\cdot,z) \asymp1$ and $F1_{\mathcal
{E}_1} \preceq
\Gamma1_{\mathcal{E}_2}$, we get by Proposition~\ref{super-concatenating}(ii),
\[
F1_{\mathcal{E}_1}A\Gamma(x,z) \leq C\eta\Gamma1_{\mathcal
{E}_2}\Gamma(x,z).
\]
Setting $\mathcal{E}_2^1 = \{y\in\mathcal{E}_2\dvtx  \llvert y-z \rrvert
\geq\llvert x-z \rrvert /2\}
$, $\mathcal{E}_2^2=
\mathcal{E}_2\setminus\mathcal{E}_2^1$, we split further into
\[
\Gamma1_{\mathcal{E}_2}\Gamma= \Gamma1_{\mathcal{E}_2^1}\Gamma+
\Gamma1_{\mathcal{E}_2^2}\Gamma. %
\]
If $y\in\mathcal{E}_2^1$, then $\Gamma(y,z) \leq C\Gamma(x,z)$. By
Lemma~\ref{super-gammalemma}(iv), $\Gamma(x, \mathcal{E}_2) \leq
C$. Together we obtain
\[
\Gamma1_{\mathcal{E}_2^1}\Gamma(x,z) \leq C\Gamma(x,z). %
\]
If $y \in\mathcal{E}_2^2$, then $\Gamma(x,y) \leq C\frac
{a(z)^2}{r^2}\Gamma(x,z)$
and $\Gamma^{(1)}(y,z) \leq C\frac{r^2}{a(z)^2}\Gamma^{(1)}(z,y)$, whence
\[
\Gamma1_{\mathcal{E}_2\setminus
\mathcal{E}_2^2}\Gamma(x,z) \leq C\Gamma(x,z)\Gamma^{(1)}(z,
\mathcal{E}_2) \leq C\Gamma(x,z). %
\]
We therefore have shown that
\[
F1_{\mathcal{E}_1}A\Gamma(x,z) \leq C\eta\Gamma(x,z). %
\]
To handle the second summand of~(\ref{super-keyest-splitting1}), set
$\sigma(y) =
\min\{\eta, (\log a(y) )^{-15/2} \}$, $y\in
V_{L+r}$. Clearly, $1_BA\Gamma(y,z) \leq C\sigma(y)\Gamma(y,z)$ and $F1_B
\preceq\Gamma1_{V_L}$. Furthermore, $\sigma(\cdot)\Gamma(\cdot,z)\asymp
1$, so that by Proposition~\ref{super-concatenating}(ii)
\[
F1_BA\Gamma(x,z) \leq C\Gamma1_{V_L}\sigma\Gamma(x,z).
\]
Consider $D^1 = \{y\in V_L\dvtx  \llvert y-z \rrvert \geq\llvert x-z
\rrvert /2\}$, $D^2 =
V_L\setminus D^1$ and split into
\[
\Gamma1_{V_L}\sigma\Gamma= \Gamma1_{D^1}\sigma\Gamma+ \Gamma
1_{D^2}\sigma\Gamma. %
\]
If $y\in D^1$, then $\Gamma(y,z) \leq
C\max\{1,\frac{\dtL(y)}{\dtL(x)} \}\Gamma(x,z)$, implying
$\Gamma1_{D^1}\sigma\Gamma(x,z) \leq C\eta^{1/2}\Gamma(x,z)$ if we prove
%
%
\begin{equation}
\label{super-keyest-toprove} \sum_{y\in
V_L}\max\biggl\{1,
\frac{\dtL(y)}{\dtL(x)} \biggr\}\Gamma(x,y)\sigma(y) \leq C\eta^{1/2}.
\end{equation}
To\vspace*{1pt} this end, we treat the summation over $S^1=\{y \in V_L\dvtx  \dL(y) \leq
2s_L\}$ and $S^2= V_L\setminus S_1$ separately. If $y\in S^2$, then
$a(y) = s_L$. Estimating $\Gamma$ by $\Gamma^{(1)}$ and $\dtL(y)$,
$\dtL(x)$
simply by $L$, we get
%
%
\begin{eqnarray}\label{super-keyest-s2}
&& \sum_{y\in
S^2}\max\biggl\{1,
\frac{\dtL(y)}{\dtL(x)} \biggr\}\Gamma(x,y)\sigma(y)\nonumber
\\
&&\qquad  \leq \frac{C}{(\log L)^{3/2}}\sum
_{y\in V_{2L}}\frac{1}{ (s_L +
\llvert y \rrvert )^d}
\\
&&\qquad \leq\frac{C\log\log L}{(\log L)^{3/2}}.\nonumber
\end{eqnarray}
If $y\in S^1$, we estimate $\Gamma$ again by $\Gamma^{(1)}$ and split the
summation into the layers $\mathcal{L}_j$, $j = 0,\ldots,2s_L$. On
$\mathcal{L}_j$, $\sigma(y) \leq C\min\{\eta,(\log
(j+ 1))^{-15/2} \}$. Thus, by Lemma~\ref{super-gammalemma}(iii),
\begin{eqnarray*}
&& \sum_{y\in S^1}\max\biggl\{1,
\frac{\dtL(y)}{\dtL
(x)} \biggr\} \Gamma(x,y)\sigma(y)
\\
&&\qquad \leq C\sum_{j=0}^{2s_L}\sum
_{y\in\mathcal{L}_j}\max\biggl\{ 1,\frac{\dtL(x)}{a(y)} \biggr\}
\frac{\min\{\eta,(\log
(j+1))^{-15/2} \}}{(a(y)+\llvert x-y\rrvert )^d}
\\
&&\qquad \leq C\sum_{j=0}^{2s_L}\frac{\min\{\eta,(\log(
j+1))^{-15/2} \}}{j\vee r}
\leq C\eta^{1/2}.
\end{eqnarray*}
Together with~(\ref{super-keyest-s2}), we have
proved~(\ref{super-keyest-toprove}). It remains to bound the
term $\Gamma1_{D^2}\sigma\Gamma(x,z)$. But if $y\in D^2$, then
\begin{eqnarray*}
a(y) + \llvert x-y\rrvert &\geq& a(y) + \tfrac{1}{2}\llvert x-z \rrvert
\geq a(z) - \tfrac{1}{2}\llvert y-z \rrvert+ \tfrac{1}{2}\llvert x-z
\rrvert
\\
&\geq&\tfrac{1}{4} \bigl(a(z) + \llvert x-z \rrvert\bigr),
\end{eqnarray*}
whence $\Gamma(x,y) \leq
C\frac{a(z)^2}{a(y)^2}\max\{1,\frac{\dtL(y)}{\dtL(z)}
\}\Gamma(x,z)$.
Using Lemma~\ref{super-gammalemma}(i), we have
\[
\frac{a(z)^2}{a(y)^2}\Gamma(y,z) \leq C\Gamma(z,y), %
\]
so that
$\Gamma1_{D^2}\sigma\Gamma(x,z) \leq C\eta^{1/2}\Gamma(x,z)$ again
follows from~(\ref{super-keyest-toprove}).
\end{pf}
Now we have collected all ingredients to finally prove part~(ii) of our
main Lemma~\ref{superlemma}.

\begin{pf*}{Proof of Lemma~\ref{superlemma}(\textup{ii})}
As already remarked, we only have to prove the statement involving
${\tilde{G}^g}$. We work with the kernel $p=p_{s_L/20}$, but suppress
it from
notation, that is, ${\tilde{\pi}}= {\tilde
{\pi}^{(p)}}$, ${\tilde{g}}={\tilde{g}}^{(p)}$. The perturbation
expansion~(\ref{prel-pbe3}) yields
\[
{\tilde{G}^g}= {\tilde{g}}\sum_{m=0}^\infty(R{
\tilde{g}})^m\sum_{k=0}^\infty
\Delta^k, %
\]
where $\Delta= 1_{V_{L+r}}({\tilde{\Pi}^g}-{\tilde{\pi}})$, $R =
\sum_{k=1}^\infty\Delta^k{\tilde{\pi}}$. With the
constants $C_1$ of
Lem\-ma~\ref{superlemma}(i) and $C_2,C_3$ of
Proposition~\ref{super-keyest} we choose
\[
\delta\leq\frac{1}{32} \biggl(\frac{1}{C_2\vee C_1^2C_3^2} \biggr). %
\]
Recall the properties of ${\tilde{\Pi}^g}$ and ${\tilde{\pi}}$
mentioned after Remark~\ref{super-remr}.
From Lem\-ma~\ref{superlemma}(i) and Proposition~\ref{super-keyest}(i)
with $A = \llvert \Delta\rrvert $, $\eta= 2\delta$ we then deduce that
${\tilde{g}}\llvert \Delta\rrvert
\preceq(C_1/2)\Gamma$, and, by iterating,
\[
\sum_{k=1}^\infty{\tilde{g}}\llvert\Delta
\rrvert^{k-1}\preceq2C_1 \Gamma. %
\]
Furthermore, by part~(ii) of Proposition~\ref{super-keyest} with
$A=\llvert \Delta{\tilde{\pi}}\rrvert $ and Lemma~\ref{superlemma}(i),
\[
\sum_{k=1}^\infty{\tilde{g}}\llvert\Delta
\rrvert^{k-1}\llvert\Delta{\tilde{\pi}}\rrvert{\tilde{g}}\preceq
(C_1/2)\Gamma. %
\]
Iterating this procedure shows that for $m\in\mathbb{N}$,
\[
{\tilde{g}}\bigl(\llvert R \rrvert{\tilde{g}}\bigr)^m \preceq
C_12^{-m}\Gamma. %
\]
Finally, by a further application of Proposition~\ref{super-keyest}(i),
\[
{\tilde{g}}\sum_{m=0}^\infty\bigl(\llvert R
\rrvert{\tilde{g}}\bigr)^m\sum_{k=0}^\infty
\llvert\Delta\rrvert^k \preceq4C_1\Gamma. %
\]
This proves the lemma.
\end{pf*}

\subsection{Modified transitions on environments bad on level 4}\label{super-cgmod}
We shall now describe an environment-depending second version of the coarse
graining scheme, which leads to modified transition kernels $\Pho_{L,r}$,
${\breve{\Pi}^g}_{L,r}$, ${\breve{\pi}}_{L,r}$ on ``really bad''
environments.

We assume that $\ctwo(\delta,L_0,L_1)$ holds, and take $L_1\leq L\leq
L_1(\log L_1)^2$, so that Lemma~\ref{superlemma} can be applied.

\begin{figure}

\includegraphics{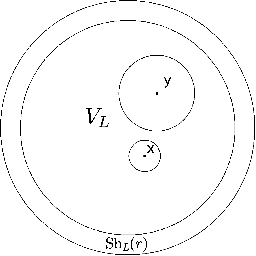}

\caption{$\omega\in\onebad$ bad on level $4$, with $\badP_L\subset
V_{L/2}$. The point $x$ is ``good,'' so the coarse graining radii do
not change at $x$. The point $y$ is ``bad.'' Therefore, at $y$, the
exit distribution is taken from the larger set $V_{t(y)}(y)$, where
$t(y)=K_1h_{L,r}(y)$.}\label{fig4}
\end{figure}

Assume $\omega\in\onebad$ is bad on level $4$, with $\badP_L(\omega
)\subset
V_{L/2}$. Then there exists $D=V_{4h_L(z)}(z)\in\mathcal{D}_L$ with
$\badP_L(\omega) \subset D$, $z\in V_{L/2}$. On $D$, $c r_L\leq
h_{L,r}(\cdot) \leq C r_L$. By Lemma~\ref{superlemma} and the
definition of $\Gamma_{L,r}$, it follows easily that we can find a constant
$K_1\geq2$, depending only on $d$, such that whenever $\llvert
x-y\rrvert \geq
K_1h_{L,r}(y)$ for some $y\in\badP_L$, we have
%
%
\begin{equation}
\label{super-kernelmod} {\hat{G}^g}_{L,r}(x,\badP_L)
\leq C\Gamma_{L,r}(x,D)\leq\tfrac{1}{10}.
\end{equation}
On such $\omega$, we let $t(x) = K_1h_{L,r}(x)$, and define on $V_L$,
\[
\Pho_{L,r}(x,\cdot) = \cases{ \ex _{V_{t(x)}(x)} (x,\cdot;
\Ph_{L,r} ),&\quad for $x\in\badP_L$,
\vspace*{2pt}\cr
\Ph_{L,r}(x,\cdot),&\quad otherwise.} %
\]
By replacing $\Ph$ by ${\hat{\pi}^{(q)}}$ on the right-hand side, we define
${\breve{\pi}^{(q)}}_{L,r}(x,\cdot)$ in an analogous way, for all
$q\in\mathcal{P}^{\mathrm{s}}_\kappa$. More precisely,
\[
{\breve{\pi}^{(q)}}_{L,r}(x,\cdot) = \cases{
\ex _{V_{t(x)}(x)} \bigl(x,\cdot;{\hat{\pi}^{(q)}}_{L,r}
\bigr),&\quad for $x\in\badP_L$,
\vspace*{2pt}\cr
{\hat{\pi}^{(q)}}_{L,r}(x,
\cdot),&\quad otherwise.} %
\]
Note that ${\breve{\pi}^{(q)}}_{L,r_L}$
depends on the environment. See Figure~\ref{fig4} for a visualization of the modified transitions.
We work again with a goodified version
of $\Pho_{L,r}$,
\[
{\breve{\Pi}^g}_{L,r}(x,\cdot) = {\breve{
\Pi}^g}_{L,r}(x,\cdot) \cases{ \ex _{V_{t(x)}(x)}
\bigl(x,\cdot;{\hat{\Pi}^g}_{L,r} \bigr), &\quad for $x\in
\badP_L$,
\vspace*{2pt}\cr
{\hat{\Pi}^g}_{L,r}(x,\cdot),&
\quad otherwise.} %
\]

For all other environments falling not into the above class, we change
nothing and put $\Pho_{L,r}= \Ph_{L,r}$,
${\breve{\Pi}^g}_{L,r}={\hat{\Pi}^g}_{L,r}$, ${\breve{\pi
}}_{L,r}=\hat{\pi}_{L,r}$. This defines
$\Pho_{L,r}$, ${\breve{\Pi}^g}_{L,r}$ and ${\breve{\pi}}_{L,r}$ on
all environments.
We write $\Gho_{L,r}$, ${\breve{G}^g}_{L,r}$, ${\breve
{g}}_{L,r}$ for the Green's
functions corresponding to $\Pho_{L,r}$, ${\breve{\Pi}^g}_{L,r}$ and
${\breve{\pi}}_{L,r}$.

%

\subsubsection{Some properties of the new transition kernels}
The following observations follow from the definition and will be
tacitly used below:
\begin{itemize}
\item On environments which are good or bad on level
at most $3$, the new kernels agree with the old ones, and so do their
Green's functions, that is, $\Gh_{L,r} = \Gho_{L,r}$ and ${\hat
{G}^g}_{L,r} =
{\breve{G}^g}_{L,r}$. On $\good$ with the choice $r=r_L$, we have
equality of all
four Green's functions.
\item If $\omega$ is not bad on level $4$ with $\badP_L\subset V_{L/2}$,
then, with $p=p_{s_L/20}$,
\[
1_{V_L} \bigl(\Pho_{L,r}-{\breve{\Pi}^g}_{L,r}
\bigr) = 1_{V_L} \bigl(\Ph_{L,r}-{\hat{\Pi}^g}_{L,r}
\bigr)= 1_{\mathcal{B}^\star_{L,r}} \bigl(\Ph-\hat{\pi}^{(p)} \bigr).
\]
This will be used in Section~\ref{exit-meas}.
\item In contrast to $\hat{\pi}_{L,r}$, the kernel
${\breve{\pi}}_{L,r}$ depends on the
environment, too. However, $\Pho_{L,r}$, ${\breve{\Pi}^g}_{L,r}$ and
${\breve{\pi}}_{L,r}$
do not change the exit measure from $V_L$, that is, for example,
\[
\ex _{V_L} \bigl(x,\cdot;{\breve{\Pi}^g}_{L,r}
\bigr) = \ex _{V_L} \bigl(x,\cdot;{\hat{\Pi}^g}_{L,r}
\bigr). %
\]
\item The old transition kernels are finer in the sense that the (new) Green's
functions $\Gho$, ${\breve{G}^g}$, ${\breve{g}}$ are pointwise bounded
from above by $\Gh$, ${\hat{G}^g}$ and ${\hat{g}}$,
respectively. In particular, we obtain with the same constants as in
Lemma~\ref{superlemma},
\end{itemize}

%
\begin{lemma} \label{superlemma2}
\textup{(i)}~For all $q\in\mathcal{P}^{\mathrm{s}}_\kappa$,
\[
{\breve{g}}^{(q)}_{L,r}\preceq C_1
\Gamma_{L,r}. %
\]

\textup{(ii)}~Assume $\ctwo(\delta,L_0,L_1)$, and let $L_1\leq L\leq L_1(\log
L_1)^2$. For $\delta>0$ small,
\[
{\breve{G}^g}_{L,r} \preceq C\Gamma_{L,r}.
\]
\end{lemma}

For the new goodified Green's function, we have
%
\begin{corollary}\label{super-cor}
Assume $\ctwo(\delta,L_0,L_1)$, and let $L_1\leq L\leq L_1(\log
L_1)^2$. There exists a constant $C > 0$ such that:
\begin{longlist}[(iii)]
\item On $\onebad$, if $\badP_L\cap\sh_L(r_L)=\varnothing$ or for general
$\badP_L$ in the case $r=r_L$,
\[
\sup_{x\in V_L}{\breve{G}^g}_{L,r}(x,
\badP_L)\leq C. %
\]
On $\onebad$, if $\badP_L\not\subset V_{L/4}$, then, with
$t=\dist(\badP_L,\partial V_L)$,
\[
\sup_{x\in
V_{L/5}}{\breve{G}^g}_{L,r}(x,
\badP_L) \leq C \biggl(\frac{s_L\wedge
(t\vee
r_L)}{L} \biggr)^{d-2}.
\]
\item On $(\bbad)^c$, $\sup_{x\in
V_{2L/3}}{\breve{G}^g}_{L,r} (x,\badP_{L,r}^\partial)
\leq
C(\log
r)^{-1/2}$.
\item For $\omega\in\onebad$ bad on level at most $3$ with
$\badP_L\cap\sh_L(r_L)=\varnothing$, or for $\omega$ bad on level $4$
with $\badP_L\subset V_{L/2}$, putting $\Delta= 1_{V_L}(\Pho
_{L,r}-{\breve{\Pi}^g}_{L,r})$,
\[
\sup_{x\in V_L}\sum_{k=0}^\infty
{\bigl\llVert\bigl({\breve{G}^g}_{L,r}1_{\badP_L}
\Delta\bigr)^k(x,\cdot)\bigr\rrVert}_1\leq C.
\]
\end{longlist}
\end{corollary}

\begin{pf}
(i) The set $\badP_L$ is contained in a neighborhood $D\in\mathcal
{D}_L$. As
${\breve{G}^g}\preceq C\Gamma$, we have
%
%
\begin{equation}
\label{super-cor-1} {\breve{G}^g} (x,\badP_L ) \leq C
\Gamma^{(2)}(x,D).
\end{equation}
From this, the first statement of (i) follows. Now let $x$ be inside
$V_{L/5}$, and $\badP_L\not\subset V_{L/4}$. If the midpoint $z$ of $D$
can be chosen to lie inside $V_L\setminus\sh(r_L)$, $a(\cdot)/h_L(z)$
and $h_L(z)/a(\cdot)$ are bounded on $D$. Then, the second statement
of (i) is
again a consequence of~(\ref{super-cor-1}). If $z\in\sh(r_L)$, we have
\begin{eqnarray*}
{\breve{G}^g} (x,\badP_L ) &\leq& C
\Gamma^{(1)}(x,D) \leq C\sum_{j=0}^{2r_L}
\sum_{y\in\mathcal{L}_j\cap D}\frac
{L}{a(y)L^{d}}
\\
&\leq& CL^{-d+1}\sum_{j=0}^{2r_L}
\frac{r_L^{d-1}}{j\vee r} \leq C(\log L) \biggl(\frac{r_L}{L} \biggr)^{d-1}.
\end{eqnarray*}

(ii) Recall the notation of Section~\ref{s2}. In order to bound
${\breve{G}^g}(x,\badP^\partial_{L,r})$ uniformly in $x\in
V_{2L/3}$, we
look at the
different bad sets $D_{j,r} \in\mathcal{Q}_{j,r}$ of layer $\Lambda_j$,
$0\leq j\leq J_1$. Estimating ${\breve{G}^g}$ by $\Gamma^{(1)}$, we have
\[
{\breve{G}^g} (x,D_{j,r} ) \leq C \bigl(r2^j
\bigr)^{d-1}L^{-d+1}. %
\]
On $(\bbad)^c$, the number of bad sets in layer $\Lambda_j$ is
bounded by
\[
C(\log r +j)^{-3/2}{\bigl(L/\bigl(r2^j\bigr)
\bigr)}^{d-1}. %
\]
Therefore,
\[
{\breve{G}^g} \bigl(x,\badP^\partial_{L,r}\cap
\Lambda_j \bigr) \leq C(\log r + j)^{-3/2}. %
\]
Summing over $0\leq j\leq J_1$, this shows
\[
{\breve{G}^g} \bigl(x,\badP^\partial_{L,r} \bigr)
\leq C(\log r)^{-1/2}. %
\]

(iii) Assume $\omega\in\good$ or $\omega$ is bad on level
$i=1,2,3$. Then
$1_{\badP_L}\Delta= 1_{\badP_L}(\Ph-{\hat{\pi}})$.
Furthermore, if $\badP_L\cap\sh_L(r_L) =
\varnothing$, we have $\llVert {\breve{G}^g}1_{\badP_L}\Delta(x,\cdot
)\rrVert _1
\leq
C\delta$. By choosing $\delta$ small enough, the claim follows. If
$\omega$ is bad on level $4$ and $\badP_L\subset V_{L/2}$, we do not gain
a factor $\delta$ from $\llVert 1_{\badP_L}\Delta(y,\cdot)\rrVert
_1$. However, thanks
to our modified transition kernels, using~(\ref{super-kernelmod}),
$\llVert 1_{\badP_L}\Delta{\breve{G}^g}1_{\badP_L}(y,\cdot)\rrVert _1
\leq1/5$
(recall that
${\breve{G}^g}\leq{\hat{G}^g}$ pointwise), so that (ii) follows in
this case, too.
\end{pf}
%
%
\begin{remark}
\label{super-remark}
All $\delta_0 > 0$ and $L_0$ appearing in the next sections are
understood to be chosen in such a way that if we take $\delta\in
(0,\delta_0]$ and $L\geq L_0$, then the conclusions of
Lemmata~\ref{superlemma},~\ref{superlemma2} and Corollary~\ref
{super-cor} are
valid.
\end{remark}

\section{Exit distributions from the ball}\label{exit-meas}
In this part, we prove the main estimates on exit measures that are
required to propagate condition~$\ctwo(\delta,L_0,\break L_1 )$. First,
in Section~\ref{exit-meas-prel} we collect some preliminary results involving
the kernel $p_L$. Then, in Section~\ref{smv-exits}, we estimate
the total variation norm of the globally smoothed difference
$D_{L,p_L,\psi,q}^\ast$, while\vspace*{1pt} in Section~\ref{nonsmv-exits}, we
prove the
estimates for the nonsmoothed quantity $D_{L,p_L}^\ast$.

We work with both the original kernels $\Ph$, ${\hat{\Pi}^g}$, $
{\hat{\pi}}$ as
well as with the modified kernels $\Pho$, ${\breve{\Pi}^g}$, $
{\breve{\pi}}$ from
Section~\ref{super-cgmod}. For the goodified exit measure from $V_L$,
we write
\[
{\Pi^g}_L= \ex _{V_L} \bigl(x,\cdot;{
\hat{\Pi}^g}_{L,r} \bigr)=\ex _{V_L}
\bigl(x,\cdot;{\breve{\Pi}^g}_{L,r} \bigr). %
\]
Throughout this part, for reasons of readability, we put $p=p_{s_L/20}$.

\subsection{Preliminaries}\label{exit-meas-prel}
We start with a generalization
of Lemma~\ref{phgp-php-estimate} which forms one of the key steps in
transferring
condition~$\ctwo$ from one level to the next.
%
\begin{lemma}
\label{phgp-phpL-estimate}
Assume $\ctwo(\delta,L_0,L_1)$, and let $L_1\leq L\leq L_1(\log L_1)^2$.\break
Then for all $x\in V_L\setminus\sh_L(2r)$, with $H(x) = \max\{L_0,
h_{L,r}(x)\}$,
\begin{eqnarray*}
&& \bigl\llVert\bigl({\hat{\Pi}^g}_{L,r} - {\hat{
\pi}}^{(p_L)}_{L,r} \bigr){\hat{\pi}}^{(p_L)}_{L,r}(x,
\cdot)\bigr\rrVert_1
\\
&&\qquad \leq C\min\bigl\{\log\bigl(s_L/H(x)\bigr) \bigl(\log H(x)
\bigr)^{-9}, \bigl(\log H(x) \bigr)^{-8} \bigr\}
\end{eqnarray*}
and
\[
\bigl\llVert\bigl({\hat{\Pi}^g}_{L,r} - {\hat{\pi
}}^{(p_L)}_{L,r} \bigr) (x,\cdot)\bigr\rrVert_1
\leq3\delta.
\]
\end{lemma}
\begin{pf}
Let $\Delta= 1_{V_L}({\hat{\Pi}^g}-{\hat{\pi}}^{(p_L)})$. With
$B = V_L\setminus\sh_L(2r)$,
\[
1_B\Delta= 1_{B} \bigl({\hat{\Pi}^g}-{
\hat{\pi}^{(p)}} \bigr) + 1_{B} \bigl({\hat{
\pi}^{(p)}} -{\hat{\pi}}^{(p_L)} \bigr). %
\]
Using Lemma~\ref{phgp-php-estimate}, the first term is bounded in total
variation by $2\delta$. For the second term, Lemma~\ref{diff-admiss}
in combination
with Lemma~\ref{php-phq-estimate} yields the bound
$C(\log H(x))^{-9}$. Similarly,
\begin{eqnarray*}
1_{B}\Delta{\hat{\pi}}^{(p_L)}
&=& 1_{B} \bigl[ \bigl({\hat{\Pi}^g}-{\hat{
\pi}^{(p)}} \bigr)\hat{\pi}^{(p)} +\hat{\pi}^{(p)}
\bigl(\hat{\pi}^{(p)}-{\hat{\pi}}^{(p_L)} \bigr)
\\
&&\hspace*{15pt}{}+ \bigl({\hat{
\pi}^{(p)}}-{\hat{\pi}}^{(p_L)} \bigr){\hat{\pi}}^{(p_L)}+{
\hat{\Pi}^g} \bigl({\hat{\pi}}^{(p_L)}-{\hat{
\pi}^{(p)}} \bigr) \bigr].
\end{eqnarray*}
Here, the last three terms on the right are bounded in total variation by
$C(\log H(x))^{-9}$, and for the first one can use
Lemma~\ref{phgp-php-estimate}.
\end{pf}

The next lemma is useful for the globally smoothed exit distributions.
%
\begin{lemma}
\label{est-deltaphi}
Assume $\ctwo(\delta,L_0,L_1)$, and let $L_1\leq L\leq L_1(\log
L_1)^2$. Put $\Delta=
1_{V_L}({\hat{\Pi}^g}_{L,r_L}-{\hat{\pi}}^{(p_L)}_{L,r_L})$.
Then,\vspace*{1pt} for some $C>0$, for all $\psi\in\mathcal{M}_L$, $q\in
\mathcal{P}^{\mathrm{s}}_\iota
$, with
$\phi_{L,p_L,\psi,q}=\pi_L^{(p_L)}{\hat{\pi}}_{\psi
}^{(q)}$ as in Section~\ref{estimates},
\[
\sup_{x\in
V_L}\sup_{z\in\mathbb{Z}^d}\bigl\llvert\Delta
\phi_{L,p_L,\psi,q}(x,z)\bigr\rrvert\leq C(\log L)^{-12}L^{-d}.
\]
\end{lemma}
\begin{pf}
Write ${\hat{\pi}}$ for ${\hat{\pi
}}^{(p_L)}$, $\phi$ for
$\phi_{L,p_L,\psi,q}$. Using $\Delta\phi= \Delta{\hat
{\pi}}\phi$ and the fact
that $\Delta{\hat{\pi}}(x,\cdot)$ sums up to zero,
\begin{eqnarray*}
\bigl\llvert\Delta\phi(x,z)\bigr\rrvert&=& \biggl\llvert\sum
_{y\in
V_L\cup \partial
V_L}\Delta{\hat{\pi}}(x,y) \bigl(\phi(y,z)-\phi(x,z) \bigr)
\biggr\rrvert
\\
&\leq&\bigl\llVert\Delta{\hat{\pi}}(x,\cdot)\bigr\rrVert_1\sup
_{y\dvtx  \llvert \Delta
{\hat{\pi}}(x,y)\rrvert > 0}\bigl\llvert\phi(y,z)-\phi(x,z)\bigr
\rrvert.
\end{eqnarray*}
For $x\in V_L\setminus\sh_L(2r_L)$, we have by
Lemma~\ref{phgp-phpL-estimate}
\[
\bigl\llVert\Delta{\hat{\pi}}(x,\cdot)\bigr\rrVert_1 \leq C\log
\bigl(s_L/h_L(x)\bigr) (\log L)^{-9}.
\]
Moreover, notice that $\llvert \Delta{\hat{\pi}}(x,y)\rrvert > 0$
implies $\llvert y-x \rrvert \leq
Ch_L(x)$. Bounding $\llvert \phi(y,z)-\phi(x,z)\rrvert $ by Lemma~\ref
{est-phi}(iii),
the statement follows for those $x$. If $x\in\sh_L(2r_L)$, we simply bound
$\llVert \Delta{\hat{\pi}}(x,\cdot)\rrVert _1$ by $2$. Now we can
restrict the supremum to
those $y\in V_L$ with $\llvert x-y\rrvert \leq3r_L$, so the claim
again follows from
Lemma~\ref{est-phi}(iii).
\end{pf}

We defined the kernel $p_L$ in terms of averaged variances of
$\pE[\Ph_{L,r}]$. Combined with Lemma~\ref{onestep-exitlaw}, this shows
that the covariances of $\hat{\pi}_{L,r}^{(p_L)}$ agree
with those of
$\pE[\Ph_{L,r}]$ up to an error of order $O(s_L^{-1})$.
The same holds true with $\Ph_{L,r}$ replaced by ${\hat{\Pi}^g}_{L,r}$.
%
\begin{lemma}
\label{lemma-covariancematching}
Assume $\ctwo(\delta,L_0,L_1)$. There exists a constant $C=C(d)$ such that
for $L_1\leq L\leq L_1(\log L_1)^2$, we have
\[
\biggl\llvert\sum_{y\in V_L} \bigl(\pE\bigl[{\hat{
\Pi}^g}_{L,r}\bigr] - {\hat{\pi}}^{(p_L)}_{L,r}
\bigr) (0,y)\frac
{y_i^2}{\llvert y \rrvert ^2}\biggr\rrvert\leq C (\log L)^3L^{-1}
\qquad\mbox{for all }i=1,\ldots,d. %
\]
\end{lemma}
\begin{pf}
Note that under
$\ctwo(\delta,L_0,L_1)$,
\[
\pE\bigl[\bigl\llVert\bigl(\Ph-{\hat{\Pi}^g} \bigr) (0,\cdot)
\bigr\rrVert_1 \bigr]\leq2\pP(0\in\badP_L )\leq C\exp
\bigl(-(1/2) (\log L)^2 \bigr). %
\]
Therefore,
\[
2p_L(e_i)=\sum_{y}\pE
\bigl[{\hat{\Pi}^g}\bigr](0,y)\frac{y_i^2}{\llvert y \rrvert ^2} + O
\bigl(\exp
\bigl(-(1/2) (\log L)^2\bigr) \bigr). %
\]
On the other hand, as in~(\ref{onestep-exitlaw-1}) with $p$ replaced by
$p_L$,
\[
2p_L(e_i)=\sum_{y}{
\hat{\pi}}^{(p_L)}(0,y)\frac
{y_i^2}{\llvert y \rrvert ^2}+ O\bigl(s_L^{-1}
\bigr), %
\]
and the statement follows.
\end{pf}

\subsection{Globally smoothed exits}\label{smv-exits}
Our objective here is to establish the estimates for the smoothed difference
$D_{L,p_L,\psi,q}^\ast$. In this section, we only work with
coarse graining schemes corresponding to the choice $r=r_L$.
Lemma~\ref{lemma-goodpart} compares the ``goodified'' smoothed exit
distribution from the ball of radius $L$ with that of symmetric random walk
with transition kernel $p_L$. In particular, it provides an estimate for
$D_{L,p_L,\psi,q}^\ast$ on $\good$.
%
\begin{lemma}
\label{lemma-goodpart}
Assume $\mathbf{A1}$. For every (small) constant $c_0$, there exist
$\delta_0
> 0$ and $L_0\in\mathbb{N}$ such that if $\delta\in(0,\delta_0]$ and
$L_1\geq L_0$, then $\ctwo(\delta,L_0,L_1)$ implies that for
$L_1\leq L\leq L_1(\log L_1)^2$, for all $\psi\in\mathcal{M}_L$ and all
$q\in\mathcal{P}^{\mathrm{s}}_\iota$,
\[
\pP\Bigl(\sup_{x\in V_L}\bigl\llVert\bigl({\Pi^g}_L
- \pi^{(p_L)}_L \bigr){\hat{\pi}^{(q)}}_{\psi}(x,
\cdot)\bigr\rrVert_1 \geq c_0(\log L)^{-9}
\Bigr) \leq\exp\bigl(-(\log L)^{7/3} \bigr). %
\]
\end{lemma}
\begin{pf}
Clearly,
the claim follows if we show that
%
%
\begin{eqnarray}
\label{smvgoodenv-toshow}
&& \sup_{x\in V_L}\sup_{z\in\mathbb{Z}^d}
\pP\bigl(\bigl\llvert\bigl({\Pi^g}- \pi^{(p_L)} \bigr){\hat{
\pi}^{(q)}}_{\psi
}(x,z)\bigr\rrvert \geq c_0(\log
L)^{-9}L^{-d} \bigr)
\nonumber\\[-8pt]\\[-8pt]\nonumber
&&\qquad \leq\exp\bigl(-(\log L)^{5/2} \bigr).
\end{eqnarray}
We use the abbreviations $\phi= \pi_L^{(p_L)}{\hat{\pi
}^{(q)}}_{\psi}$, $\Delta=
1_{V_L}({\hat{\Pi}^g}-{\hat{\pi}}^{(p_L)})$, ${\hat
{g}}={\hat{g}}^{(p_L)}$ and ${\hat{\pi
}}={\hat{\pi}}^{(p_L)}$.
By the perturbation expansion [cf.~(\ref{readingguide-pert})],
\[
\bigl({\Pi^g}-\pi^{(p_L)} \bigr){\hat{\pi}^{(q)}}_{\psi}
= {\hat{G}^g}\Delta\phi. %
\]
Set $S = \sh_L(2L/(\log L)^2)$, and write
%
%
\begin{equation}
\label{smvgoodenv-splitting} {\hat{G}^g}\Delta\phi= {\hat{G}^g}1_{S}
\Delta\phi+ {\hat{G}^g}1_{S^c}\Delta\phi.
\end{equation}
Using ${\hat{G}^g}\preceq C\Gamma$, Lemma~\ref{super-gammalemma}(iv) (with
$r= r_L$) and Lemma~\ref{est-deltaphi} yield the estimate
\[
\bigl\llvert{\hat{G}^g}1_{S}\Delta\phi(x,z)\bigr
\rrvert\leq\sup_{x\in V_L}{\hat{G}^g}(x,S)\sup
_{y\in
V_L}\bigl\llvert\Delta\phi(y,z)\bigr\rrvert\leq(\log
L)^{-19/2}L^{-d} %
\]
for $L$ large. It remains to bound $\llvert {\hat{G}^g}1_{{S}^c}\Delta
\phi
(x,z)\rrvert $. With
$B = V_L\setminus\sh_L(2r_L)$,
\[
{\hat{G}^g}= {\hat{g}}1_{B}\Delta{\hat{G}^g}+{
\hat{g}}1_{B^c}\Delta{\hat{G}^g}+{\hat{g}}. %
\]
By replacing successively ${\hat{G}^g}$ in the first summand on the right-hand
side,
%
%
\begin{eqnarray}
\label{smvgoodenv-splitting2} {\hat{G}^g}1_{{S}^c}\Delta\phi&=& \Biggl
(\sum
_{k=0}^\infty{ ({\hat{g}}1_B
\Delta)}^k{\hat{g}}+ \sum_{k=0}^\infty{
({\hat{g}}1_B\Delta)}^k{\hat{g}}1_{B^c}
\Delta{\hat{G}^g} \Biggr)1_{{S}^c}\Delta\phi
\nonumber\\[-8pt]\\[-8pt]\nonumber
&=& F1_{{S}^c}\Delta\phi+F1_{B^c}\Delta{\hat{G}^g}1_{{S}^c}
\Delta\phi,
\end{eqnarray}
where $F=\sum_{k=0}^\infty{ ({\hat{g}}1_B\Delta
)}^k{\hat{g}}$.
With $R = \sum_{k=1}^{\infty}(1_B\Delta)^k{\hat{\pi}}$,
expansion~(\ref{prel-pbe3}) shows
\[
F = {\hat{g}}\sum_{m=0}^\infty(R{
\hat{g}})^m\sum_{k=0}^\infty
(1_B\Delta)^k. %
\]
From the proof of Lemma~\ref{superlemma}(ii) we learn that
$\llvert F \rrvert \preceq C\Gamma$.
By Lemma~\ref{super-gammalemma}(iv),\break (v) and again
Lemma~\ref{est-deltaphi}, we see that for large $L$, uniformly in
$x\in
V_L$ and $z\in\mathbb{Z}^d$,
\begin{eqnarray*}
&& \bigl\llvert F1_{B^c}\Delta{\hat{G}^g}1_{S^c}
\Delta\phi(x,z)\bigr\rrvert
\\
&&\qquad \leq C\Gamma\bigl(x,\sh_L(2r_L)\bigr)\sup
_{v\in\sh_L(3r_L)}\Gamma\bigl(v,S^c\cap V_L
\bigr)\sup_{w\in
V_L}\bigl\llvert\Delta\phi(w,z)\bigr\rrvert
\\
&&\qquad \leq(\log L)^{-11}L^{-d}.
\end{eqnarray*}
Thus, the second summand
of~(\ref{smvgoodenv-splitting2}) is harmless. However, with the first
summand one has to be more careful. With $\xi=
{\hat{g}}\sum_{k=0}^\infty(1_B\Delta)^k1_{S^c}\Delta
\phi$, we have
\[
F1_{S^c}\Delta\phi= \xi+ {\hat{g}}\sum_{m=0}^\infty
(R{\hat{g}})^mR\xi= \xi+ F1_{B}\Delta{\hat{\pi}}\xi.
\]
Clearly, $\llvert F1_B\Delta{\hat{\pi}}(x,y)\rrvert \leq C(\log
L)^{-3}$, so it remains
to estimate $\xi(y,z)$, uniformly in $y$ and $z$. Set $N=N(L)=\lceil
\log\log L\rceil$. For small $\delta$, the summands of $\xi$ with
$k\geq N$
are readily bounded by
\begin{eqnarray*}
\sup_{y\in V_L}\sup_{z\in\mathbb{Z}^d}\sum
_{k=N}^\infty\bigl\llvert{\hat{g}}(1_B
\Delta)^k1_{S^c}\Delta\phi(y,z)\bigr\rrvert&\leq& C(\log
L)^6\sum_{k=N}^\infty
\delta^k(\log L)^{-12}L^{-d}
\\
&\leq&(\log L)^{-10}L^{-d}.
\end{eqnarray*}
Now we look at the summands with $k<N$. Since the coarse grained walk
cannot bridge a gap of length $L/(\log L)^2$ in less than $N$ steps, we can
drop the kernel $1_B$. Defining $S' = \sh_L (3L/(\log L)^2
)$, we
thus have
\[
{\hat{g}}(1_B\Delta)^k1_{S^c}\Delta\phi= {
\hat{g}}1_{S'}\Delta^k1_{S^c}\Delta\phi+ {
\hat{g}}1_{S'^c}\Delta^k1_{S^c}\Delta\phi.
\]
The first summand is bounded in the same way as
${\hat{G}^g}1_S\Delta\phi$ from~(\ref{smvgoodenv-splitting}), and
we can
drop the kernel $1_{S^c}$ in the second
summand. Therefore,~(\ref{smvgoodenv-toshow}) follows if we show
\begin{eqnarray*}
&& \sup_{x\in V_L}\sup_{z\in
\mathbb{Z}^d}\pP\Biggl(\Biggl
\llvert\sum_{k=1}^N{\hat
{g}}1_{S'^c}\Delta^k\phi(x,z)\Biggr\rrvert\geq
\frac{1}{2}c_0(\log L)^{-9}L^{-d} \Biggr)
\leq\exp\bigl(-(\log L)^{5/2} \bigr).
\end{eqnarray*}
For $j\in\mathbb{Z}$, consider the interval $I_j
= (jNs_L,(j+1)Ns_L]\cap\mathbb{Z}$. We divide $ S'^c\cap V_L$ into subsets
$W_{\bj} = (S'^c\cap V_L)\cap(I_{j_1}\times\cdots\times
I_{j_d} )$,
where ${\bj} = (j_1,\ldots,j_d) \in\mathbb{Z}^d$.
Let $J$ be the set of those $\bj$ for which $W_{\bj} \neq
\varnothing$. Then we can find a constant $K$ depending only on the dimension
and a disjoint partition of $J$ into sets $J_1,\ldots,J_K$, such that
for any $1\leq
\ell\leq K$,
%
%
\begin{equation}
\label{smvgoodenv-distance} \bj,\bj' \in J_\ell,\qquad \bj\neq
\bj'\quad\Longrightarrow\quad\dist(W_{\bj},W_{\bj'}) >
Ns_L.
\end{equation}
For $x\in V_L$, $z\in\mathbb{Z}^d$, we set
\[
\xi_{\bj} = \xi_{\bj}(x,z) = \sum
_{y\in W_{\bj}}\sum_{k=1}^N{
\hat{g}}(x,y)\Delta^k\phi(y,z), %
\]
and further $t = t(d,c_0,L) = (1/2)c_0(\log L)^{-9}L^{-d}$.
Assume that we can prove
%
%
\begin{equation}
\label{smvgoodenv-expbound} \biggl\llvert\sum_{\bj\in J}\pE[
\xi_{\bj} ]\biggr\rrvert\leq\frac{t}{2}.
\end{equation}
Then
\begin{eqnarray*}
\pP\biggl(\biggl\llvert\sum_{\bj\in J}
\xi_{\bj}\biggr\rrvert\geq t \biggr) &\leq&\pP\biggl(\biggl\llvert
\sum_{\bj\in
J}\xi_{\bj}-\pE[
\xi_{\bj} ]\biggr\rrvert\geq\frac{t}{2} \biggr)
\\
&\leq& K \max_{1\leq\ell\leq K}\pP\biggl(\biggl\llvert\sum
_{\bj\in
J_\ell}\xi_{\bj}-\pE[\xi_{\bj} ]\biggr
\rrvert\geq\frac
{t}{2K} \biggr).
\end{eqnarray*}
Due to~(\ref{smvgoodenv-distance}), the random variables
$\xi_{\bj}-\pE[\xi_{\bj} ]$, $\bj\in J_\ell$,
are independent and centered. Hoeffding's inequality yields, with
${\llVert \xi_{\bj}\rrVert }_{\infty} = \sup_{\omega\in\Omega
}\llvert \xi_{\bj
}(\omega)\rrvert $,
%
%
\begin{equation}
\label{smvgoodenv-hoeffding} \pP\biggl(\biggl\llvert\sum_{\bj\in
J_\ell}
\xi_{\bj}-\pE[\xi_{\bj} ]\biggr\rrvert\geq
\frac{t}{2K} \biggr) \leq2\exp\biggl(-c\frac{L^{-2d}(\log
L)^{-18}}{\sum_{\bj\in
J_\ell}{\llVert \xi_{\bj}\rrVert }^2_{\infty}} \biggr)
\end{equation}
for some constant $c>0$. In order to control the $\sup$-norm of the
$\xi_{\bj}$, we use the estimates
\begin{eqnarray*}
{\hat{g}}(x,W_{\bj}) &\leq& C\Gamma^{(2)}(x,W_{\bj})
\leq\frac{CN^ds_L^d}{s_L^2(s_L + \dist(x,W_{\bj}))^{d-2}}
= CN^d \biggl(1+\frac{\dist(x,W_{\bj})}{s_L} \biggr)^{2-d},
\end{eqnarray*}
and, by Lemmata~\ref{phgp-phpL-estimate} and~\ref{est-deltaphi}, for
$y\in
W_{\bj}$, $\llvert \Delta^k\phi(y,z)\rrvert \leq C\delta
^{k-1}k(\log
L)^{-12}L^{-d}$. Altogether we arrive at
\[
\llVert\xi_{\bj}\rrVert_{\infty} \leq C \biggl(1+
\frac{\dist(x,W_{\bj})}{s_L} \biggr)^{2-d}N^d(\log
L)^{-12}L^{-d}, %
\]
uniformly in $z$. If we put the last display into~(\ref
{smvgoodenv-hoeffding}), we get,
using $d\geq3$ in the last line,
\begin{eqnarray*}
\pP\biggl(\biggl\llvert\sum_{\bj\in
J_\ell}
\xi_{\bj}-\pE[\xi_{\bj} ]\biggr\rrvert\geq
\frac{t}{2K} \biggr)&\leq&2\exp\biggl(-c\frac{(\log L)^{6}}{N^4\sum
_{r=1}^{C(\log
L)^3/N}r^{-d+3}} \biggr)
\\
&\leq&2\exp\biggl(-c\frac{(\log L)^3}{N^3} \biggr).
\end{eqnarray*}
It follows that for $L$ large enough, uniformly in $x$ and $z$,
\[
\pP\biggl(\biggl\llvert\sum_{\bj\in J}
\xi_{\bj}\biggr\rrvert\geq\frac{1}{2}c_0(\log
L)^{-9}L^{-d} \biggr) \leq\exp\bigl(-(\log
L)^{5/2} \bigr). %
\]
It remains to prove~(\ref{smvgoodenv-expbound}). We have
\[
\biggl\llvert\sum_{\bj\in J}\pE[\xi_{\bj} ]
\biggr\rrvert\leq\sum_{y\in S'^c}{\hat{g}}(x,y)\Biggl
\llvert\sum_{y'\in
V_L}\pE\Biggl[\sum
_{k=1}^N\Delta^k{\hat{\pi}}
\bigl(y,y'\bigr) \Biggr]\phi\bigl(y',z\bigr)\Biggr
\rrvert. %
\]
The sum over the Green's function is estimated by ${\hat{g}}(x,S'^c)
\leq C(\log L)^6$. In the innermost sum, we
treat the summands with $k\geq2$ and $k=1$ in different ways. For each
summand with $k\geq2$, we use Proposition~\ref{est-symmetry} applied
to \mbox{$\nu(\cdot) = \pE[\Delta^k{\hat{\pi
}}(y,y+\cdot) ]$}; see Remark~\ref{smv-exits-remark}.
Recalling Lemma~\ref{phgp-phpL-estimate}, we obtain by
choosing $\delta$ small enough, uniformly in $y\in S'^c$,
\begin{eqnarray*}
\biggl\llvert\sum_{y'\in V_L}\pE\bigl[
\Delta^k{\hat{\pi}}\bigl(y,y'\bigr) \bigr]\phi
\bigl(y',z\bigr)\biggr\rrvert&\leq& C\delta^{k-1}(\log
L)^{-9}(ks_L)^2L^{-(d+2)}
\\
&\leq&(c_0/8) (1/2)^{k-1}(\log L)^{-15}L^{-d}.
\end{eqnarray*}
For the summand $\nu(\cdot) = \pE[\Delta{\hat
{\pi}}(y,y+\cdot) ]$
corresponding to $k=1$, the proof of
Proposition~\ref{est-symmetry} shows
\begin{eqnarray*}
&& \biggl\llvert\sum_{y'\in V_L}\pE\bigl[\Delta{\hat{\pi
}}\bigl(y,y'\bigr) \bigr]\phi\bigl(y',z\bigr)\biggr
\rrvert
\\
&&\qquad \leq\Biggl\llvert\frac{1}{2}\sum_{i=1}^d
\frac{\partial^2}{\partial
x_i^2}{\phi}^{\mathrm{B}}(y,z)\sum_{y'}
\nu\bigl(y'\bigr) \bigl(y'_i
\bigr)^2\Biggr\rrvert
 + C(\log L)^{-18}L^{-d}.
\end{eqnarray*}
By translational invariance and Lemma~\ref{lemma-covariancematching} in
the last step,
\begin{eqnarray*}
\biggl\llvert\sum_{y'}\nu\bigl(y'
\bigr) \bigl(y'_i\bigr)^2\biggr\rrvert&=&
\biggl\llvert\sum_{y} \bigl(\pE\bigl[{\hat{
\Pi}^g}\bigr]-{\hat{\pi}} \bigr){\hat{\pi}}(0,y)y_i^2
\biggr\rrvert
 = \biggl\llvert\sum_{y} \bigl(\pE
\bigl[{\hat{\Pi}^g}\bigr]-{\hat{\pi}} \bigr) (0,y)y_i^2
\biggr\rrvert
\\
&\leq& CL/(\log L)^3.
\end{eqnarray*}
Together with Lemma~\ref{est-phi}(ii), this
shows~(\ref{smvgoodenv-expbound}) and completes the proof.
\end{pf}

%
\begin{remark}
\label{smv-exits-remark}
Note that for $y\in S'^c$, the signed
measure $\nu$ fulfills the requirement of
Proposition~\ref{est-symmetry}. Indeed, after $N=\lceil\log\log
L\rceil$ steps away from $y$, the coarse grained walks are still in the
interior part of $V_L$, where the coarse graining radius did not start
to shrink. Therefore, the symmetry condition $\mathbf{A1}$ carries
over to
the signed measure
$\pE[\sum_{k=1}^N(1_{V_L}(\Ph-{\hat{\pi
}}))^k{\hat{\pi}}(y,y+\cdot)]$. Replacing
$\Ph$ by ${\hat{\Pi}^g}$ does not destroy the symmetries of this
measure, so that Proposition~\ref{est-symmetry} can be applied to $\nu$.
\end{remark}

Next, we estimate $D_{L,p_L,\psi,q}^\ast$ on environments with bad
points. Again, we make the choice $r=r_L$.
%
\begin{lemma}
\label{lemma-badpart}
In the setting of Lemma~\ref{lemma-goodpart}, with a possibly smaller
value of $\delta_0$ and a larger $L_0$, we have for $i=1,2,3,4$,
\begin{eqnarray*}
&& \pP\Bigl(\sup_{x\in
V_{L/5}}\bigl\llVert\bigl(\Pi_L-
\pi^{(p_L)}_L \bigr){\hat{\pi}^{(q)}}_{\psi
}(x,
\cdot)\bigr\rrVert_1 > (\log L)^{-9+
9(i-1)/4};
\operatorname{OneBad}_L^{(i)} \Bigr)
\\
&&\qquad \leq\exp\bigl(-(\log L)^{7/3} \bigr).
\end{eqnarray*}
\end{lemma}

\begin{pf}
For the ease of readability, we write
${\hat{\pi}}_\psi$ for ${\hat{\pi}^{(q)}}_\psi$ and let $\phi=
\pi
^{(p_L)}{\hat{\pi}}_{\psi}$. By the
triangle inequality,
%
%
\begin{eqnarray}\label{smvbadpart-1}
&& \bigl\llVert\bigl(\Pi- \pi^{(p_L)} \bigr){\hat{\pi
}}_{\psi}(x,\cdot)\bigr\rrVert_1
\nonumber\\[-8pt]\\[-8pt]\nonumber
&&\qquad \leq\bigl\llVert\bigl(
\Pi-{\Pi^g} \bigr)\hat{\pi}_{\psi}(x,\cdot)\bigr\rrVert
_1 + \bigl\llVert\bigl({\Pi^g}- \pi^{(p_L)}
\bigr)\hat{\pi}_{\psi
}(x,\cdot)\bigr\rrVert_1.
\end{eqnarray}
The second summand on the right can be bounded by
Lemma~\ref{lemma-goodpart}. For the first term we have by the
perturbation expansion [see~(\ref{readingguide-pert})], with
$\Delta= 1_{V_L}(\Pho-{\breve{\Pi}^g})$,
\[
\bigl(\Pi-{\Pi^g} \bigr){\hat{\pi}}_\psi= {
\breve{G}^g}1_{\badP
_L}\Delta\Pi{\hat{\pi}}_\psi.
\]
Note that since we are on $\onebad$, the set $\badP_L$ is contained
in a
small region. First assume that $\badP_L\subset\sh_L(L/(\log
L)^{10})$. Then\vspace*{1pt} $\sup_{x\in V_{L/5}}{\breve{G}^g}(x,\badP_L) \leq
C(\log L)^{-10}$
by Corollary~\ref{super-cor}, which bounds the first summand
of~(\ref{smvbadpart-1}). Next assume $\omega$ bad on level $4$ and
$\badP_L\not\subset V_{L/2}$. Then $\sup_{x\in V_{L/5}}{\breve{G}^g}
(x,\badP_L)
\leq C(\log L)^{-3}$ by the same corollary, which is enough for this case.

It remains to consider the cases $\omega$ bad on level at most $3$ with
$\badP_L\not\subset\sh_L(L/(\log L)^{10})$, or $\omega$ bad on
level $4$ with
$\badP_L\subset V_{L/2}$. We expand
\begin{eqnarray*}
\bigl(\Pi-{\Pi^g} \bigr){\hat{\pi}}_\psi&=& \bigl({
\breve{G}^g} 1_{\badP_L}\Delta\Pi\bigr){\hat{
\pi}}_\psi
\\
&=& \sum_{k=1}^\infty
\bigl({\breve{G}^g}1_{\badP_L}\Delta\bigr)^k{
\Pi^g} {\hat{\pi}}_\psi
\\
&=& \sum_{k=1}^\infty\bigl({
\breve{G}^g}1_{\badP_L}\Delta\bigr)^k\phi+ \sum
_{k=1}^\infty\bigl({\breve{G}^g}1_{\badP_L}
\Delta\bigr)^k \bigl({\Pi^g} -\pi^{(p_L)} \bigr){
\hat{\pi}}_{\psi}
\\
&=& F_1 + F_2.
\end{eqnarray*}
By Corollary~\ref{super-cor},
\begin{eqnarray*}
\bigl\llVert F_1(x,\cdot)\bigr\rrVert_1 &\leq&\sum
_{k=0}^\infty\bigl\llVert\bigl({
\breve{G}^g}1_{\badP_L}\Delta\bigr)^k(x,\cdot)
\bigr\rrVert_1\sup_{v\in
V_L}{\breve{G}^g}(v,
\badP_L)\sup_{w\in
\badP_L}\bigl\llVert\Delta\phi(w,
\cdot)\bigr\rrVert_1
\\
&\leq& C\sup_{w\in
\badP_L}\bigl\llVert\Delta\phi(w,\cdot)\bigr\rrVert
_1.
\end{eqnarray*}
Proceeding as in Lemma~\ref{est-deltaphi}, for $w\in\badP_L$,
%
%
\begin{equation}
\label{smvbadpart-2} \bigl\llVert\Delta\phi(w,\cdot)\bigr\rrVert_1\leq
\bigl\llVert\Delta{\hat{\pi}}^{(p_L)}(w,\cdot)\bigr\rrVert_1
\sup_{w'\dvtx
\llvert \Delta{\hat{\pi}}^{(p_L)}(w,w')\rrvert >
0}\bigl\llVert\phi\bigl(w',\cdot\bigr)-
\phi(w,\cdot)\bigr\rrVert_1.\hspace*{-30pt}
\end{equation}
Lemma~\ref{est-phi}(iii) bounds the second factor on the right by
$C(\log L)^{-3}$. If $\omega$ is bad on level $4$, we simply bound the
first factor by $2$. If $\omega$ is bad on level at most $3$, we have on
$\badP_L$ the equality $\Delta=\Ph-\hat{\pi}^{(p)}$. With
$p'=p_{h_L(w)}$, the
triangle inequality gives, for every $i=1,2,3$,
\begin{eqnarray*}
\bigl\llVert\Delta{\hat{\pi}}^{(p_L)}(w,\cdot)\bigr\rrVert
_1&\leq&\bigl\llVert\bigl(\Ph-{\hat{\pi}}^{(p')} \bigr) {
\hat{\pi}}^{(p')}(w,\cdot)\bigr\rrVert_1 + \bigl\llVert
\bigl({\hat{\pi}}^{(p')}-{\hat{\pi}}^{(p)} \bigr) (w,\cdot)
\bigr\rrVert_1
\\
&&{} + \sup_{y\in
V_L\setminus
\sh_L(r_L)}\bigl\llVert\bigl({\hat{\pi
}}^{(p_L)}-{\hat{\pi}}^{(p')} \bigr) (y,\cdot)\bigr\rrVert
_1.
\end{eqnarray*}
By definition, the first summand is bounded by $C(\log L)^{-9+9i/4}$. For
the second and third summand, we use Lemma~\ref{php-phq-estimate} and the
fact that by Lemma~\ref{diff-admiss}, $\llVert p'-p_L\rrVert _1\leq
C(\log\log
L)(\log L)^{-9}$,
and similarly for $\llVert p-p'\rrVert _1$. For all $i=1,2,3,4$, we
arrive at
$\llVert F_1(x,\cdot)\rrVert _1\leq C(\log L)^{-12+9i/4}$. For $F_2$,
we obtain once
more with Corollary~\ref{super-cor},
\[
\bigl\llVert F_2(x,\cdot)\bigr\rrVert_1 \leq C\sup
_{y\in V_L} \bigl\llVert\bigl({\Pi^g}-
\pi^{(p_L)} \bigr){\hat{\pi}}_{\psi
}(y,\cdot)\bigr\rrVert
_1. %
\]
This term is again estimated by Lemma~\ref{lemma-goodpart} (with
$c_0$ small enough), which completes the proof.
\end{pf}

\subsection{Nonsmoothed and locally smoothed exits}
\label{nonsmv-exits}
Here, we aim at bounding the total variation distance of the exit measures
without additional smoothing (Lemma~\ref{nonsmv-lemma1}), as well as in
the case where the exit measures are convoluted with a kernel of constant
smoothing radius $s$ (Lemma~\ref{nonsmv-lemma2}).

Throughout this part, we work with transition kernels coming from coarse
graining schemes with constant
parameter $r$. We always assume $L$ large enough such that $r< r_L$. The
right choice of $r$ depends on the deviations $\delta$ and $\eta$ we are
shooting for and will become clear from the proofs. In either case, we
choose $r\geq r_0$, where $r_0$ is the constant from Section~\ref{s2}. The
value of $r$ will then also influence the choice of the perturbation
$\varepsilon_0$
in Lemma~\ref{nonsmv-lemma1} and the smoothing radius $\ell$ in
Lemma~\ref{nonsmv-lemma2}, respectively.

We recall the partition of bad points into the sets $\badP_L$,
$\badP_{L,r}$, $\badP_{L,r}^\partial$, $\mathcal{B}^\star_{L,r}$
and the
classification of environments into $\good$, $\onebad$ and $\bbad$ from
Section~\ref{smoothbad}.

The bounds for $\manybad$ (Lemma~\ref{smoothbad-lemmamanybad}) and for
$\bbad$ (Lemma~\ref{smoothbad-lemmabbad}) ensure that we may restrict
ourselves to environments
\[
\omega\in\onebad\cap\,(\bbad)^c.
\]
For such
environments, we introduce two disjoint sets $Q_{L,r}^1(\omega),
Q_{L,r}^2(\omega) \subset V_L$ as follows:
\begin{itemize}
\item If $\badP_L(\omega)\subset V_{L/2}$, set $Q_{L,r}^1(\omega)=
\badP_L(\omega)$ and $Q_{L,r}^2(\omega)= \badP_{L,r}^\partial
(\omega)$.
\item If $\badP_L(\omega)\not\subset V_{L/2}$, set $Q_{L,r}^1(\omega)=
\varnothing$ and $Q_{L,r}^2(\omega)= \mathcal{B}^\star_{L,r}(\omega)$.
\end{itemize}
Of course, on $\good$, we have $Q_{L,r}^1(\omega) = \varnothing$ and $
Q_{L,r}^2(\omega) = \badP_{L,r}^\partial(\omega)$.

Recall that we write
$p$ for $p_{s_L/20}$.
%
\begin{lemma}
\label{nonsmv-lemma1}
There exists $\delta_0 > 0$ such that if $\delta\in(0,\delta_0]$, there
exist $\varepsilon_0 = \varepsilon_0(\delta) > 0$ and $L_0 =
L_0(\delta)>0$ with the following
property: if $\varepsilon\leq\varepsilon_0$ and $L_1\geq L_0$, then
$\mathbf{A0}(\varepsilon)$, $\ctwo(\delta,L_0,L_1)$
imply that for $L_1\leq L\leq L_1(\log L_1)^2$,
\[
\pP\Bigl(\sup_{x\in V_{L/5}}\bigl\llVert\bigl(\Pi_L-
\pi^{(p_L)}_L \bigr) (x,\cdot)\bigr\rrVert_1 >
\delta\Bigr) \leq\exp\biggl(-\frac{9}{5}(\log L)^2 \biggr).
\]
\end{lemma}
\begin{pf}
We choose $\delta_0 > 0$ according to Remark~\ref{super-remark} and
take $\delta\in(0, \delta_0]$. The right choice of $\varepsilon_0$
and $L_0$
will be clear from the course of the proof. From
Lemmata~\ref{smoothbad-lemmamanybad}
and~\ref{smoothbad-lemmabbad} we learn that if we take $L_1\geq L_0$ for
$L_0$ large and $L_1\leq L \leq L_1(\log L_1)^2$, then under $\ctwo
(\delta,L_0,L_1)$,
\[
\pP(\manybad\cup\,\bbad) \leq\exp\bigl(-\tfrac{9}{5}(\log
L)^2 \bigr). %
\]
Therefore, the claim follows if we show that on $\onebad\cap\, (\bbad)^c$, we have for all sufficiently small $\varepsilon$ and
all large $L$,
$x\in V_{L/5}$,
\[
\bigl\llVert\bigl(\Pi- \pi^{(p_L)} \bigr) (x,\cdot)\bigr\rrVert
_1 \leq\delta. %
\]
Let $\omega\in\onebad\cap\,(\bbad)^c$. We use the partition of
$\mathcal{B}^\star_{L,r}$ into the sets $Q^1$, $Q^2$ described
above. With
$\Delta=
1_{V_L}(\Pho-{\breve{\Pi}^g})$, we have inside $V_L$,
\[
\Pi= {\breve{G}^g}1_{Q^1}\Delta\Pi+ {
\breve{G}^g}1_{Q^2}\Delta\Pi+ {\Pi^g}.
\]
By replacing successively $\Pi$ in the first summand on the right-hand
side, we arrive at
\[
\Pi= \sum_{k=0}^\infty{ \bigl({
\breve{G}^g}1_{Q^1}\Delta\bigr)}^k{
\Pi^g}+ \sum_{k=0}^\infty{
\bigl({\breve{G}^g}1_{Q^1}\Delta\bigr)}^k{
\breve{G}^g} 1_{Q^2}\Delta\Pi. %
\]
Since with $\Delta' = 1_{V_L}({\breve{\Pi}^g}-{\breve{\pi
}}^{(p_L)})$, ${\Pi^g}= \pi^{(p_L)} + {\breve{G}^g}\Delta'\pi
^{(p_L)}$, we obtain
%
%
\begin{eqnarray}
\label{nonsmv-lemma1-splitting}
\Pi-\pi^{(p_L)}
&=& \sum_{k=1}^\infty{ \bigl({
\breve{G}^g}1_{Q^1}\Delta\bigr)}^k\pi
^{(p_L)} + \sum_{k=0}^\infty{ \bigl({
\breve{G}^g}1_{Q^1}\Delta\bigr)}^k{
\breve{G}^g} 1_{Q^2}\Delta\Pi\nonumber
\\
&&{} + \sum_{k=0}^\infty{ \bigl({
\breve{G}^g}1_{Q^1}\Delta\bigr)}^k{
\breve{G}^g}\Delta'\pi^{(p_L)}
\\
&=& F_1 + F_2 + F_3.\nonumber
\end{eqnarray}
We will now prove that each of the three parts $F_1$, $F_2$, $F_3$ is bounded
by $\delta/3$. If $Q^1 \neq\varnothing$, then $Q^1 = \badP_L \subset
V_{L/2}$ and $Q^2 =
\badP_{L,r}^\partial$. Using Corollary~\ref{super-cor} in
the second and Lemma~\ref{lemmalawler}(ii) in the third inequality,
%
%
\begin{eqnarray}
\label{nonsmv-lemma1-a1} \qquad\bigl\llVert F_1(x,\cdot)\bigr\rrVert_1
&\leq&\sum_{k=0}^\infty\bigl\llVert\bigl({
\breve{G}^g}1_{\badP_L}\Delta\bigr)^k(x,\cdot)
\bigr\rrVert_1\sup_{y\in
V_L}{\breve{G}^g}(y,
\badP_L)\sup_{z\in\badP_L}\bigl\llVert
\Delta'\pi^{(p_L)}(z,\cdot)\bigr\rrVert_1
\nonumber
\\
&\leq& C \sup_{z\in V_{L/2}}\bigl\llVert\Delta'
\pi^{(p_L)}(z,\cdot)\bigr\rrVert_1
\leq C(\log L)^{-3}
\\
&\leq& C (\log L_0)^{-3} \leq\delta/3\nonumber
\end{eqnarray}
for $L_0 = L_0(\delta)$ large enough, $L\geq L_0$.
Regarding $F_2$, we have in the case $Q^1 \neq\varnothing$
by Corollary~\ref{super-cor}(ii),
\[
\bigl\llVert F_2(x,\cdot)\bigr\rrVert_1 \leq C\sup
_{y\in
V_{2L/3}}{\breve{G}^g}\bigl(y,\badP_{L,r}^\partial
\bigr) \leq C(\log r)^{-1/2}. %
\]
On the other hand, if $Q^1 = \varnothing$, then $\badP_L$ is outside
$V_{L/3}$, so that by Corollary~\ref{super-cor}(i),~(ii)
\[
\bigl\llVert F_2(x,\cdot)\bigr\rrVert_1
\leq2{\breve{G}^g}\bigl(x,\badP^\partial_{L,r}\cup
\badP_L\bigr)\leq C \bigl((\log L)^{-3} + (\log
r)^{-1/2} \bigr). %
\]
Altogether, for all $L \geq L_0$, by choosing $r = r(\delta)$ and $L_0 =
L_0(\delta,r)$ large enough,
%
%
\begin{equation}
\label{smv-lemma1-a2-2} \bigl\llVert F_2(x,\cdot)\bigr\rrVert_1
\leq C \bigl((\log L_0)^{-3} + (\log r)^{-1/2}
\bigr) \leq\delta/3.
\end{equation}
It remains to handle $F_3$. Once again with Corollary~\ref
{super-cor}(iii) for some $C_3 > 0$,
\[
\bigl\llVert F_3(x,\cdot)\bigr\rrVert_1 \leq
C_3\sup_{y\in V_{2L/3}}\bigl\llVert{\hat
{G}^g}\Delta'\pi^{(p_L)}(y,\cdot)\bigr\rrVert
_1. %
\]
We write
%
%
\begin{equation}
\label{nonsmv-lemma1-a3} {\breve{G}^g}\Delta'
\pi^{(p_L)}={\breve{G}^g}1_{V_L\setminus
Q^1}
\Delta'\pi^{(p_L)} + {\breve{G}^g}1_{Q^1}
\Delta'\pi^{(p_L)}.
\end{equation}
As in~(\ref{nonsmv-lemma1-a1}), we deduce that
%
%
\begin{equation}
\label{nonsmv-lemma1-a3-0} \bigl\llVert{\breve{G}^g}1_{Q^1}
\Delta'\pi^{(p_L)}\bigr\rrVert_1\leq
C_3^{-1}\delta/21
\end{equation}
for $L_0$ large enough, $L \geq L_0$. Concerning the first term
of~(\ref{nonsmv-lemma1-a3}), we note that on $V_L\setminus Q^1$, by
definition,
\[
\Delta'\pi^{(p_L)}= \bigl({\hat{\Pi}^g}-{\hat{
\pi}}^{(p_L)} \bigr){\hat{\pi}}^{(p_L)}\pi^{(p_L)}.
\]
For $z\in V_L\setminus(\sh_L(2r_L)\cup Q^1)$, we obtain with
Lemma~\ref{phgp-phpL-estimate},
\[
\bigl\llVert\Delta'\pi^{(p_L)}(z,\cdot)\bigr\rrVert
_1 \leq C(\log\log L) (\log L)^{-9}. %
\]
Since ${\hat{G}^g}(y,V_L) \leq C(\log L)^6$, it follows that
%
%
\begin{equation}
\label{nonsmv-lemma1-a3-1} \sup_{y\in V_{2L/3}}\bigl\llVert{\hat{G}^g}
1_{V_L\setminus(Q^1\cup \sh_L(2r_L))}\Delta'\pi^{(p_L)}(y,\cdot)\bigr
\rrVert
_1 \leq C(\log L)^{-2} \leq C_3^{-1}
\delta/21\hspace*{-30pt}
\end{equation}
for $L$ large. Recall the definition of the layers $\Lambda_j$ from
Section~\ref{s2}. For $z\in
\Lambda_j$, $1\leq j\leq J_1$, we have again by
Lemma~\ref{phgp-phpL-estimate},
\[
\bigl\llVert\Delta'\pi^{(p_L)}(z,\cdot)\bigr\rrVert
_1 \leq C(\log r + j)^{-8}. %
\]
Using Lemma~\ref{super-gammalemma}(iii), it follows that ${\hat{G}^g}
(y,\Lambda_j) \leq C$ for
some constant $C$, independent of $r$ and $j$. Thus
%
%
\begin{equation}
\label{nonsmv-lemma1-a3-2} \sup_{y\in V_{2L/3}}\bigl\llVert{\hat{G}^g}
1_{\bigcup_{j=1}^{J_1}\Lambda_j}\Delta'\pi^{(p_L)}(y,\cdot)\bigr\rrVert
_1 \leq C (\log r)^{-7} \leq C_3^{-1}
\delta/21,
\end{equation}
if $r$ is chosen large enough. For the first layer $\Lambda_0$, there
is a constant $C_0$ with
\[
\sup_{y\in V_{2L/3}}\bigl\llVert{\hat{G}^g}
1_{\Lambda_0}\Delta'\pi^{(p_L)}(y,\cdot)\bigr\rrVert
_1 \leq C_0\sup_{z\in\Lambda_0}\bigl\llVert
\Delta'(z,\cdot)\bigr\rrVert_1. %
\]
Now, with $p_o \equiv1/(2d)$ denoting the transition kernel of simple
random walk,
%
%
\begin{eqnarray}
\label{nonsmv-lemma1-a3-3}
\bigl\llVert\Delta'(z,\cdot)\bigr\rrVert
_1 &\leq& \bigl\llVert\bigl(\Ph-{\hat{\pi}}^{(p_o)} \bigr) (z,\cdot)\bigr
\rrVert_1 + \bigl\llVert\bigl({\hat{\pi}}^{(p_o)}-{\hat{\pi
}}^{(p_L)} \bigr) (z,\cdot)\bigr\rrVert_1
\nonumber\\[-8pt]\\[-8pt]\nonumber
&&{}  + \bigl\llVert
\bigl({\hat{\Pi}^g}-\Ph\bigr) (z,\cdot)\bigr\rrVert_1.
\nonumber
\end{eqnarray}
Concerning the second summand of~(\ref{nonsmv-lemma1-a3-3}), recall that
by Lemma~\ref{diff-admiss}, $\llVert p_o -p_L\rrVert _1\leq(\log L_0)^{-7}$.
Therefore, choosing $L_0$ large enough ($r$ is now fixed), we can
guarantee that
\[
\bigl\llVert\bigl({\hat{\pi}}^{(p_o)}-{\hat{\pi}}^{(p_L)} \bigr)
(z,\cdot)\bigr\rrVert_1 \leq C_0^{-1}C_3^{-1}
\delta/21, %
\]
uniformly in $z\in\Lambda_0$. The third summand
of~(\ref{nonsmv-lemma1-a3-3}) vanishes if
$z\in\Lambda_0\setminus\mathcal{B}^\star_{L,r}$. For $z\in
\mathcal{B}^\star_{L,r}$,
\begin{eqnarray*}
\bigl\llVert\bigl({\hat{\Pi}^g}-\Ph\bigr) (z,\cdot)\bigr\rrVert
_1&=&\bigl\llVert\bigl(\hat{\pi}^{(p)}-\Ph\bigr) (z,
\cdot)\bigr\rrVert_1
\\
&\leq&\bigl\llVert\bigl(\Ph-{\hat{\pi}}^{(p_o)} \bigr) (z,\cdot)\bigr
\rrVert_1+\bigl\llVert\bigl({\hat{\pi}}^{(p_o)}-{\hat{
\pi}}^{(p)} \bigr) (z,\cdot)\bigr\rrVert_1.
\end{eqnarray*}
The last term on the right is bounded in the same way as the second term
of~(\ref{nonsmv-lemma1-a3-3}). For the first summand on the right of the
inequality [and also for the first summand
of~(\ref{nonsmv-lemma1-a3-3})], we may simply choose
$\varepsilon_0=\varepsilon_0(\delta,r)$ small enough such that for
$\varepsilon\leq\varepsilon_0$,
\[
\sup_{z\in\Lambda_0}\bigl\llVert\bigl(\Ph-{\hat{\pi}}^{(p_o)}
\bigr) (z,\cdot)\bigr\rrVert_1\leq C_0^{-1}C_3^{-1}
\delta/21. %
\]
Altogether we have shown that $\llVert F_3(x,\cdot)\rrVert _1\leq
\delta/3$, and
the claim is proved.
\end{pf}
%
%
\begin{remark}
As the proof shows, we do not have to assume condition~$\ctwo(\delta
,L_0,L_1)$ for
the desired deviation $\delta$. We could instead assume
$\ctwo(\delta',L_0,L_1)$ for some arbitrary $0<\delta'\leq\delta
_0$. However,
$L_0$ depends on the chosen $\delta$. This observation will be useful in
the next lemma.
\end{remark}

%
\begin{lemma}
\label{nonsmv-lemma2}
There exists $\delta_0 > 0$ with the following property: for each
\mbox{$\eta>
0$}, there exist a smoothing radius $\ell_0=\ell_0(\eta)$ and $L_0 =
L_0(\eta)$
such that if $L_1\geq L_0$, $\ell\geq\ell_0$ and $\ctwo(\delta
,L_0,L_1)$ holds for
some $\delta\in(0,\delta_0]$, then for $L_1\leq L\leq L_1(\log L_1)^2$
and $\psi\equiv\ell$, for all $q\in{\mathcal{P}^{\mathrm
{s}}_\iota}$ 
%
\[
\pP\Bigl(\sup_{x\in V_{L/5}}\bigl\llVert\bigl(\Pi_L-
\pi^{(p_L)}_L \bigr){\hat{\pi}^{(q)}}_{\psi}(x,
\cdot)\bigr\rrVert_1 > \eta\Bigr) \leq\exp\biggl(-
\frac{9}{5}(\log L)^2 \biggr). %
\]
\end{lemma}
\begin{pf}
The proof is based on a modification of the computations in the foregoing
lemma. Let $\delta_0$ be as in Lemma~\ref{nonsmv-lemma1}. We fix an
arbitrary $0 < \delta\leq\delta_0$ and assume $\ctwo(\delta
,L_0,L_1)$ for
some $L_1 \geq L_0$, where $L_0 = L_0(\eta)$ will be chosen later. In the
following, ``good'' and ``bad'' is always to be understood with respect
to $\delta$. Again, for $L_1\leq L\leq L_1(\log L_1)^2$,
\[
\pP(\manybad\cup\bbad) \leq\exp\bigl(-\tfrac{9}{5}(\log
L)^2 \bigr). %
\]
For $\omega\in\onebad\cap\,(\bbad)^c$, we use the
splitting~(\ref{nonsmv-lemma1-splitting}) of $\Pi-\pi^{(p_L)}$ into
the parts $F_1$, $F_2$, $F_3$. For the summands $F_1$ and $F_2$, we do not
need the additional smoothing by ${\hat{\pi}^{(q)}}_{\psi}$, since
by~(\ref{nonsmv-lemma1-a1}),
\[
\bigl\llVert F_1(x,\cdot)\bigr\rrVert_1 \leq C(\log
L)^{-3}\leq\eta/3, %
\]
and by~(\ref{smv-lemma1-a2-2}),
\[
\bigl\llVert F_2(x,\cdot)\bigr\rrVert_1 \leq C \bigl((
\log L)^{-3} + (\log r)^{-1/2} \bigr) \leq\eta/3, %
\]
if\vspace*{1.5pt} $L\geq L_0$ and $r$, $L_0$ are chosen large enough, depending on $d$
and $\eta$. We turn to~$F_3$. By~(\ref{nonsmv-lemma1-a3-0}),~(\ref{nonsmv-lemma1-a3-1})
and~(\ref{nonsmv-lemma1-a3-2}) we have, with $\Delta'=
1_{V_L}({\breve{\Pi}^g}-{\breve{\pi}}^{(p_L)})$ as in the proof of
Lemma~\ref{nonsmv-lemma1}, writing again ${\hat{\pi
}}_\psi$ for ${\hat{\pi}^{(q)}}_\psi$,
%
%
\begin{eqnarray}
\label{nonsmv-lemma2-a3}
\qquad&& \bigl\llVert F_3{\hat{\pi}}_s(x,
\cdot)\bigr\rrVert _1
\nonumber
\\
&&\qquad \leq C \Bigl(\sup_{y\in V_{2L/3}}\bigl\llVert{\breve
{G}^g}1_{V_L\setminus
\Lambda_0}\Delta'\pi^{(p_L)}(y,
\cdot)\bigr\rrVert_1 + \sup_{z\in\Lambda_0}\bigl\llVert
\Delta'\pi^{(p_L)}{\hat{\pi}}_{\psi}(z,\cdot)\bigr
\rrVert_1 \Bigr)
\\
&&\qquad \leq C \Bigl((\log L)^{-3} +(\log r)^{-8} + \sup
_{z\in
\Lambda_0}\bigl\llVert\Delta'\pi^{(p_L)}{
\hat{\pi}}_{\psi
}(z,\cdot)\bigr\rrVert_1 \Bigr)\nonumber
\\
&&\qquad \leq \eta/6 + C_1\sup_{z\in
\Lambda_0}\bigl\llVert
\Delta'\pi^{(p_L)}{\hat{\pi}}_{\psi
}(z,\cdot)\bigr
\rrVert_1,\nonumber
\end{eqnarray}
if $L\geq L_0$ and $r$, $L_0$ are sufficiently large. Regarding the
second summand of~(\ref{nonsmv-lemma2-a3}), set $m = 3r$ and define for
$K\in\mathbb{N}$
\[
\vartheta_K(z) = \min\bigl\{n\in\mathbb{N}\dvtx  \bigl\llvert
X_n^z -z\bigr\rrvert> Km \bigr\}\in[0,\infty],
\]
where 
$X_n^z$ denotes symmetric random walk with law $\Prw_{z,p_L}$.
By the invariance principle for symmetric random walk, we can clearly
choose $K$ so large such that
\[
\sup_{z\in V_L\dvtx  \dL(z) \leq m}\Prw_{z,q} \bigl(\vartheta_K(z)
\leq\tau_L \bigr) \leq\frac{\eta}{24C_1} %
\]
uniformly in $L\geq L_0$ and $q\in{\mathcal{P}^{\mathrm{s}}_\iota
}$, where $C_1$ is the constant
from~(\ref{nonsmv-lemma2-a3}). If
$z\in\Lambda_0$, $z'\in V_L\cup\partial V_L$ with
$\Delta'(z,z') \neq0$, we have $\dL(z') \leq m$ and $\llvert z-z'
\rrvert \leq
m$. Thus, using Lemma~\ref{app-kernelest}(iii) of the \hyperref[appe]{Appendix} with
$\psi\equiv\ell$,
\begin{eqnarray*}
\hspace*{-3pt}&& C_1\sup_{z\in\Lambda_0}\bigl\llVert
\Delta'\pi^{(p_L)}{\hat{\pi}}_{\psi}(z,\cdot)\bigr
\rrVert_1
\\
\hspace*{-3pt}&&\qquad = C_1\sup_{z\in\Lambda_0}\biggl\llVert\sum
_{z'\in V_L\cup \partial
V_L}\Delta'\bigl(z,z'\bigr)\mathop{\sum
_{w\in\partial V_L\dvtx }}_{\llvert z'-w\rrvert > Km}\pi^{(p_L)}\bigl(z',w
\bigr) \bigl({\hat{\pi}}_{\psi}(w,\cdot)-{\hat{\pi}}_{\psi}(z,
\cdot) \bigr)
\\
\hspace*{-3pt}&&\hspace*{70pt}{}+ \sum_{z'\in V_L\cup\partial
V_L}\Delta'
\bigl(z,z'\bigr)
\\
\hspace*{-3pt}&&\hspace*{141pt}{}\times \mathop{\sum_{w\in\partial
V_L\dvtx }}_{\llvert z'-w\rrvert \leq Km}
\pi^{(p_L)}\bigl(z',w\bigr) \bigl({\hat{\pi}}_{\psi}(w,
\cdot)-{\hat{\pi}}_{\psi}(z,\cdot) \bigr)\biggr\rrVert_1
\\
\hspace*{-3pt}&&\qquad \leq\frac{\eta}{6}+C(K+1)m\frac{\log\ell}{\ell} \leq\eta/3,
\end{eqnarray*}
if we choose $\ell= \ell(\eta,r)$ large enough. This proves the lemma.
\end{pf}

\section{Proofs of the main results}\label{proofmain}

\mbox{}
\begin{pf*}{Proof of Proposition~\ref{main-prop}}
We take $\delta$ small enough and choose $L_0 =
L_0(\delta)$ large enough according to
Remark~\ref{super-remark} and the statements of
Section~\ref{exit-meas}. In the course of this
proof, we might enlarge $L_0$ further.
(ii) is a consequence of Lemma~\ref{nonsmv-lemma2}, so we have to
prove (i). Let $L_1 \geq L_0$, and assume that condition~$\ctwo(\delta,L_0,L_1)$ holds. Then the first point of
$\ctwo(\delta,L_0,L_1(\log L_1)^2)$ is trivially fulfilled. Now let $L_1<
L\leq L_1(\log L_1)^2$, and consider first $L'=L$ and $i=1,2,3$. Take
$\psi\in\mathcal{M}_L$, $q\in\mathcal{P}^{\mathrm{s}}_\iota$.
For simplicity, write
$D_L^\ast$ for
$D_{L,p_L,\psi,q}^\ast$. Then by Lemma~\ref{smoothbad-lemmamanybad},
\begin{eqnarray*}
&& b_i(L,p_L,\psi,q,\delta)
\\
&&\qquad \leq \pP\bigl(D_L^\ast> (\log L)^{-9 +
9(i-1)/4}
\bigr)
\\
&&\qquad \leq \pP(\manybad) + \pP\bigl(D_L^\ast> (\log
L)^{-9 +
9(i-1)/4}; \onebad\bigr)
\\
&&\qquad \leq\exp\bigl(-\tfrac{19}{10}(\log L)^2 \bigr) + \pP
\bigl(D_L^\ast> (\log L)^{-9 + 9(i-1)/4}; \onebad\bigr).
\end{eqnarray*}
For the last summand, we have by
Lemmata~\ref{lemma-goodpart} and \ref{lemma-badpart}, under $\ctwo
(\delta,L_0,\break L_1)$,
\begin{eqnarray*}
&& \pP\bigl(D_L^\ast> (\log L)^{-9 + 9(i-1)/4};
\onebad\bigr)
\\
&&\qquad \leq\pP\bigl(D_L^\ast> (\log L)^{-9};
\good\bigr)
\\
&&\quad\qquad{} + \sum_{j=1}^4 \pP
\bigl(D_L^\ast> (\log L)^{-9 +
9(i-1)/4};
\operatorname{OneBad}_L^{(j)} \bigr)
\\
&&\qquad \leq\exp\bigl(-(\log L)^{7/3} \bigr) + \sum
_{j=1}^i\pP\bigl(D_L^\ast>
(\log L)^{-9 +
9(i-1)/4}; \operatorname{OneBad}_L^{(j)}
\bigr)
\\
&&\quad\qquad{} +\sum_{j=i+1}^4\pP\bigl(
\operatorname{OneBad}_L^{(j)} \bigr)
\\
&&\qquad \leq 4 \exp\bigl(-(\log L)^{7/3} \bigr) + CL^ds_L^d
\exp\bigl(- \bigl((3+i+1)/4 \bigr) \bigl(\log(r_L/20)
\bigr)^2 \bigr).
\end{eqnarray*}
Therefore, for $L$ large,
\[
\pP\bigl(D_L^\ast> (\log L)^{-9 +
9(i-1)/4}; \onebad
\bigr)\leq\tfrac{1}{8}\exp\bigl(- \bigl((3+i)/4 \bigr) (\log L
)^2 \bigr) %
\]
and
\[
b_i(L,p_L,\psi,q,\delta) \leq\tfrac{1}{4}\exp
\bigl(- \bigl((3+i)/4 \bigr) (\log L )^2 \bigr). %
\]
For the case $i=4$, notice that
\[
b_4(L,p_L,\psi,q,\delta) \leq\pP\bigl(D_L^\ast>
(\log L)^{-9/4} \bigr) + \pP\bigl(D_{L,p_L}^\ast>
\delta\bigr). %
\]
The first summand can be estimated as the corresponding terms in the case
$i=1,2,3$, while for the last term we use Lemma~\ref{nonsmv-lemma1}.

It remains to show that for all $L$ with $L_1< L\leq L_1(\log L_1)^2$
and all $L'\in
(L,2L]$, for all $\psi\in\mathcal{M}_{L'}$ and all $q\in{\mathcal
{P}^{\mathrm{s}}_\iota}$, all
$i=1,2,3,4$,
\[
b_i\bigl(L',p_L,\psi,q,\delta\bigr)\leq
\tfrac{1}{4}\exp\bigl(- \bigl((3+i)/4 \bigr) \bigl(\log L'
\bigr)^2 \bigr). %
\]
This, however, is easily deduced from the estimates above, by a slight change
of the coarse graining scheme. Indeed, defining for $L'\in(L,2L]$ the
coarse graining function
$\tilde{h}_{L',r}\dvtx  \overline{C}_{L'} \rightarrow\mathbb{R}_+$ as
\[
\label{smoothbad-hLr} \tilde{h}_{L',r}(x) = \frac{1}{20}\max\biggl
\{s_L h \biggl(\frac{\dt_{L'}(x)}{s_{L'}} \biggr), r \biggr\}, %
\]
it follows that $\tilde{h}_{L',r}(x) = h_{L,r}(0)=s_L/20$ for $x\in V_{L'}$
with $\dt _{L'}(x) \geq2s_{L'}$. With an analogous definition of
good and bad points within $V_{L'}$, using the coarse graining function
$\tilde{h}_{L',r}$ instead of $h_{L',r}$ and the transition kernels
corresponding to $\tilde{h}_{L',r}$, clearly all the statements of
Sections~\ref{estimates} and~\ref{super} remain true. Moreover, we
can use
the kernel $p_L$ to obtain the results of Section~\ref{exit-meas} for the
radius $L'$, noticing that in the proofs at most the constants
change. Arguing now exactly as above for the choice
$L'=L$, we conclude that also the second point of condition~$\ctwo(\delta,L_0,L_1(\log L_1)^2)$ holds true, provided $L_0$ is large,
$L_1\geq L_0$.
\end{pf*}

\begin{pf*}{Proof of Proposition~\ref{main-propkernel}}
According to Proposition~\ref{main-prop}(i), for $\delta$,
$\varepsilon>0$
small and $L_0=L_0(\delta)$ large, conditions $\mathbf{A1}$, $\mathbf{A0}(\varepsilon)$ and
$\ctwo(\delta,L_0,L_0)$ imply $\ctwo(\delta,L_0,L)$ for all $L\geq
L_0$. By shrinking $\varepsilon$ if necessary, we may further
guarantee that
$\ctwo(\delta,L_0,L_0)$ holds true, as it was explained below the
statement of the proposition. Therefore, we can assume that
$\ctwo(\delta,L_0,L)$ is satisfied for all $L\geq L_0$. By
Lemma~\ref{diff-admiss}, the limit $\lim_{L\rightarrow\infty}p_L(e_i)$
exists for each $i=1,\ldots,d$. Now let $\ell=s_L/20$. Using the definition
of $p_L$ in the first, $\ctwo(\delta,L_0,L)$ in the second and fourth and
Lemma~\ref{onestep-exitlaw} in the third equality,
\begin{eqnarray*}
2p_L(e_i) &=& \sum_{y}
\pE\bigl[\Ph_{L,r}(0,y) \bigr]\frac{\llvert y \rrvert ^2}{y_i^2}
\\
&=& \sum
_{y}\hat{\pi}_{L,r}^{(p_\ell)}(0,y)
\frac
{\llvert y \rrvert ^2}{y_i^2} + O\bigl((\log L)^{-9}\bigr)
\\
&=&\sum_{y}\pi_\ell^{(p_\ell)}(0,y)
\frac{\llvert y \rrvert ^2}{y_i^2} + O\bigl((\log L)^{-9}\bigr)
\\
&=&\sum_{y}\pE\bigl[\Pi_\ell(0,y)
\bigr]\frac
{\llvert y \rrvert ^2}{y_i^2}+O\bigl((\log L)^{-9}\bigr).
\end{eqnarray*}
From this, the first statement of the proposition follows. Moreover,
$p_\infty\in\mathcal{P}^{\mathrm{s}}_\iota$, and
recalling that we may first choose $L_0$ as large as we wish and then
choose $\varepsilon$
sufficiently small, we see that $\llVert p_\infty-p_o\rrVert
_1\rightarrow0$ as
$\varepsilon\downarrow0$.
\end{pf*}

\begin{pf*}{Proof of Theorem~\ref{main-theorem1}}
As in the proof of Proposition~\ref{main-propkernel}, if the
parameters are appropriately chosen, $\ctwo(\delta/2,L_0,L)$ holds true
for all $L\geq L_0$. This implies
%
%
\begin{equation}
\label{mainthm1-1} \pP\bigl(D_{L,p_L}^\ast> \delta/2 \bigr)
\leq\exp\bigl(-(\log L)^2 \bigr).
\end{equation}
With $p_\infty$ defined in Proposition~\ref{main-propkernel}, we obtain
by the triangle inequality
%
%
\begin{eqnarray}
\label{mainthm1-2} D_{L,p_\infty}^\ast&=&\sup_{x\in
V_{L/5}}
\bigl\llVert\bigl(\Pi-\pi^{(p_\infty)} \bigr) (x,\cdot)\bigr\rrVert_1
\nonumber\\[-8pt]\\[-8pt]\nonumber
&\leq& D_{L,p_L}^\ast+ \sup_{x\in V_{L/5}}\bigl
\llVert\bigl(\pi^{(p_\infty)}-\pi^{(p_L)} \bigr) (x,\cdot)\bigr\rrVert_1.
\end{eqnarray}
The claim therefore follows if we show that the second summand is bounded
by $\delta/2$ if $L\geq L_0$ and $L_0$ is large enough. To prove this
bound, we move inside $V_L$ according to the coarse grained transition
kernels ${\hat{\pi}}^{(p_L)}_{L,r}$ and $
{\hat{\pi}}^{(p_\infty)}_{L,r}$, respectively,
where similarly to Section~\ref{nonsmv-exits}, $r$ is a large but fixed
number. First, we recall from the proof of
Proposition~\ref{main-propkernel} that
$p_\infty=\lim_{L\rightarrow\infty}p_L$. We therefore deduce from
Lemma~\ref{diff-admiss},
%
%
\begin{equation}
\llVert p_\infty-p_L\rrVert_1 \leq C\sum
_{k=0}^\infty\bigl(\log\bigl(2^kL
\bigr)\bigr)^{-9} \leq C(\log L)^{-8}.
\end{equation}
This implies by Lemma~\ref{php-phq-estimate} that for $x\in V_L$ with
$\dL(x)>(1/10)r$,
%
%
\begin{equation}
\label{mainthm1-4} \bigl\llVert\bigl({\hat{\pi}}^{(p_L)}-{\hat{\pi
}}^{(p_\infty)} \bigr) (x,\cdot)\bigr\rrVert_1 \leq C\max\bigl
\{h_{L,r}(x)^{-1/4}, (\log L)^{-8} \bigr\}.
\end{equation}
Now by the perturbation expansion, with
$\Delta=1_{V_L}({\hat{\pi}}^{(p_L)}-{\hat
{\pi}}^{(p_\infty)})$,
\begin{eqnarray*}
\pi^{(p_\infty)}-\pi^{(p_L)} &=& {\hat{g}}^{(p_L)}\Delta
\pi^{(p_\infty)}
\\
&=& {\hat{g}}^{(p_L)}1_{V_L\setminus\sh_L(2r_L)}\Delta\pi^{(p_\infty)} +{
\hat{g}}^{(p_L)}1_{\sh_L(2r_L)}\Delta\pi^{(p_\infty)}.
\end{eqnarray*}
Since ${\hat{g}}^{(p_L)}(x,V_L)\preceq C\Gamma(x,V_L)\leq
C(\log L)^{6}$, we
obtain with~(\ref{mainthm1-4})
\[
\sup_{x\in
V_L}\bigl\llVert{\hat{g}}^{(p_L)}1_{V_L\setminus\sh
_L(2r_L)}
\Delta\pi^{(p_\infty)}(x,\cdot)\bigr\rrVert_1 \leq C(\log
L)^{-2} \leq\delta/4, %
\]
if $L$ is large enough. Moreover, as
below~(\ref{nonsmv-lemma1-a3-2}), for $x\in V_{L/5}$,
using again~(\ref{mainthm1-4}),
\begin{eqnarray*}
&& \bigl\llVert{\hat{g}}^{(p_L)}1_{\sh_L(2r_L)}\Delta
\pi^{(p_\infty)}(x,\cdot)\bigr\rrVert_1
\\
&&\qquad \leq\bigl\llVert{\hat{g}}^{(p_L)} 1_{\bigcup_{j=1}^{J_1}\Lambda
_j}\Delta(x,\cdot)\bigr
\rrVert_1 + \bigl\llVert{\hat{g}}^{(p_L)}1_{\Lambda_0}
\Delta(x,\cdot)\bigr\rrVert_1
\\
&&\qquad \leq C r^{-1/4} + C_0\sup_{z\in\Lambda_0}\bigl
\llVert\Delta(z,\cdot)\bigr\rrVert_1.
\end{eqnarray*}
We choose $r$ so large such that $Cr^{-1/4}\leq\delta/8$. Now
that $r$ is fixed, the difference over the first layer
$\sh_L(2r)$ is also bounded by $C_0^{-1}\delta/4$ if the difference between
$p_L$ and $p_\infty$ is small enough, that is, if $L$ is large enough. This
proves that the second summand of~(\ref{mainthm1-2}) is bounded by
$\delta/2$, for $L\geq L_0$ and $L_0$ large. Applying~(\ref{mainthm1-1}),
we conclude that
\[
\pP\bigl(D_{L,p_\infty}^\ast>\delta\bigr) \leq\pP
\bigl(D_{L,p_L}^\ast>\delta/2 \bigr) \leq\exp\bigl(-(\log
L)^2 \bigr). %
\]\upqed
\end{pf*}
Since Theorem~\ref{main-theorem2} is proved in a similar way,
using the second part of
Proposition~\ref{main-prop}, we may safely omit the details and turn
now to the
proof of the local estimates.

\begin{pf*}{Proof of Theorem~\ref{local-thm-exitmeas}}
We choose $\delta>0$ small and $L_0(\delta)$ large enough according to
Proposition~\ref{main-prop} such that $\mathbf{A0}(\varepsilon)$ and
$\mathbf{A1}$ imply
$\ctwo(\delta,L_0,L)$ for all $L\geq L_0$. We let $r=r_L$. Recall the
definition of $\good$ from Section~\ref{smoothbad}. With $A_L = \good$,
we note that similar to Lemma~\ref{smoothbad-lemmamanybad}, if $L\geq
L_0$,
\[
\pP\bigl(A_L^c\bigr) \leq\exp\bigl(-(1/2) (\log
L)^2 \bigr). %
\]
For the rest of the proof, we take $\omega\in A_L$. On such
environments, $\Gh$ equals ${\hat{G}^g}$ by our choice of $r=r_L$. We fix
$z\in\partial V_L$ and simply write $W_t$
for $W_t(z)=V_t(z)\cap\partial V_L$.

Let us now prove part (i). We recall that $t\geq r$, and observe that
$W_t$ can be covered by $K\lfloor(t/r)^{d-1}\rfloor$ many neighborhoods
$V_{3r}(y_i)$, $y_i\in\sh_L(r)$, as defined in
Section~\ref{super-Ghestimates}, where $K$ depends on the dimension
only. Applying Lemma~\ref{superlemma}(ii), we deduce that
\begin{eqnarray*}
\Pi_L(x,W_t)&=&{\hat{G}^g}(x,W_t)\leq \sum_{i=1}^{K\lfloor(t/r)^{d-1}\rfloor}{\hat{G}^g}
\bigl(x,V_{3r}(y_i) \bigr)
\\
&\leq& C\sum_{i=1}^{K\lfloor(t/r)^{d-1}\rfloor}\Gamma
\bigl(x,V_{3r}(y_i) \bigr) \leq C \biggl(
\frac{t}{L} \biggr)^{d-1},
\end{eqnarray*}
where for the last inequality we have used that
$\Gamma(x,V_{3r}(y_i) ) \leq C(r/L)^{d-1}$ for some constant
$C=C(d,\eta)$ (recall that we assume $x\in V_{\eta L}$ for some
$0<\eta<1$). From Lemma~\ref{lemmalawler}(i) we know that if $x\in
V_{\eta L}$, then there is a constant $c=c(d,\eta)$ such that
\[
\pi^{(p_o)}(x,z) \geq cL^{-(d-1)}. %
\]
Together with the preceding inequality, this shows (i).

(ii) We recall that for (ii), we assume $t\geq L/(\log L)^6$. If not
explicitly mentioned otherwise, the underlying one-step transition kernel
is in the following given by $p_L$ defined in~(\ref{kernelpL}), which we
omit from notation. Set $\ell=(\log L)^{13/2} r$, and consider the smoothing
kernel ${\hat{\pi}}_\psi$ with $\psi\equiv\ell\in
\mathcal{M}_\ell$. Let
\[
U_\ell(W_t) = \bigl\{y\in\mathbb{Z}^d\dvtx
\dist(y,W_t) \leq2\ell\bigr\}. %
\]
We claim that
%
%
\begin{eqnarray}
\label{local-eq1} \Pi(x,W_t) -\pi(x,W_{t+6\ell}) &\leq&(\Pi-
\pi){\hat{\pi}}_\psi\bigl(x,U_\ell(W_t)
\bigr),
\\
\label{local-eq2} \pi(x,W_{t-6\ell}) -\Pi(x,W_t) &\leq&(\pi-
\Pi){\hat{\pi}}_\psi\bigl(x,U_\ell(W_{t-6\ell})
\bigr).
\end{eqnarray}
Concerning the first inequality,
\[
\Pi{\hat{\pi}}_\psi\bigl(x,U_\ell(W_t)\bigr)
\geq\sum_{y\in
W_t}\Pi(x,y){\hat{\pi}}_\psi
\bigl(y,U_\ell(W_t)\bigr) = \Pi(x,W_t),
\]
since ${\hat{\pi}}_\psi(y,U_\ell(W_t))=1$ for $y \in
W_t$. Also,
\[
\pi{\hat{\pi}}_\psi\bigl(x,U_\ell(W_t)\bigr)
= \sum_{y\in
W_{t+6\ell}}\pi(x,y){\hat{\pi}}_\psi
\bigl(y,U_\ell(W_t)\bigr) \leq\pi(x,W_{t+6\ell}),
\]
since ${\hat{\pi}}_\psi(y,U_\ell(W_t))=0$ for $y\in
\partial V_L\setminus
W_{t+6\ell}$; see Figure~\ref{fig5}. This proves~(\ref{local-eq1}),
while~(\ref{local-eq2}) is entirely similar. In the remainder of this
proof, we often write $\llvert F \rrvert (x,y)$ for $\llvert
F(x,y)\rrvert $. If we show
%
%
\begin{equation}
\label{local-eq3} \bigl\llvert(\pi-\Pi){\hat{\pi}}_\psi\bigr\rrvert
\bigl(x,U_\ell(W_t)\bigr) \leq O \bigl((\log
L)^{-5/2} \bigr)\pi(x,W_t),
\end{equation}
then by~(\ref{local-eq1}),
\begin{eqnarray*}
\Pi(x,W_t) &\leq&\pi(x,W_{t+6\ell}) + O\bigl((\log
L)^{-5/2}\bigr)\pi(x,W_t)
\\
&=& \pi(x,W_t) + \pi(x,W_{t+6\ell}\setminus
W_t) +O\bigl((\log L)^{-5/2}\bigr)\pi(x,W_t)
\\
&=& \pi(x,W_t) \bigl(1+O \bigl(\max\bigl\{\ell/t,(\log
L)^{-5/2}\bigr\} \bigr) \bigr)
\\
&=& \pi(x,W_t) \bigl(1+O\bigl((\log L)^{-5/2}\bigr) \bigr),
\end{eqnarray*}
where in the next-to-last line, we used that $\pi(x,W_{t+6\ell})\leq
C(\ell/t)\pi(x,W_t)$ by Lemma~\ref{lemmalawler}(i).
On the other hand, by~(\ref{local-eq2}) and still assuming~(\ref{local-eq3}),
\begin{eqnarray*}
\Pi(x,W_t) &\geq&\pi(x,W_{t-6\ell})-O\bigl((\log
L)^{-5/2}\bigr)\pi(x,W_t)
\\
&=&\pi(x,W_{t}) \bigl(1-O\bigl((\log L)^{-5/2}\bigr) \bigr),
\end{eqnarray*}
%
%
\begin{figure}

\includegraphics{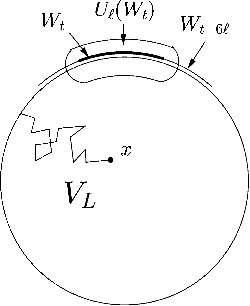}

\caption{On the proof of Theorem~\protect\ref{local-thm-exitmeas}\textup{(ii)}.
If the
walk exits $V_L$ through $\partial V_L\setminus W_{t+6\ell}$, it
cannot enter $U_\ell(W_t)$ in one step with ${\hat{\pi
}}_\psi$, $\psi\equiv\ell$.}\label{fig5}
\end{figure}
so that indeed
\[
\Pi(x,W_t) = \pi(x,W_t) \bigl(1+O \bigl((\log
L)^{-5/2} \bigr) \bigr), %
\]
provided\vspace*{1pt} we prove~(\ref{local-eq3}).
In that direction, set $B =
V_L\setminus\sh_L(5r)$ and write, with
$\Delta=1_{V_L}({\hat{\Pi}^g}-{\hat{\pi}}^{(p_L)})$,
%
%
\begin{equation}
\label{local-eq4} (\pi-\Pi){\hat{\pi}}_\psi= {\hat{G}^g}
\Delta\pi{\hat{\pi}}_\psi= {\hat{G}^g}1_B
\Delta\pi{\hat{\pi}}_\psi+ {\hat{G}^g}1_{\sh_L(5r)}
\Delta\pi{\hat{\pi}}_\psi.
\end{equation}
Looking at the first summand we have
\[
\bigl\llvert{\hat{G}^g}1_B\Delta\pi{\hat{
\pi}}_\psi\bigr\rrvert\bigl(x,U_\ell(W_t)
\bigr) \leq\bigl({\hat{G}^g}1_B\llvert\Delta{\hat{
\pi}}\rrvert\pi\bigr) (x,W_{t+6\ell}). %
\]
By Lemma~\ref{phgp-phpL-estimate}, we have
\[
\bigl\llVert1_B\Delta\pi(x,\cdot)\bigr\rrVert_1\leq
C(\log\log L) (\log L)^{-9}. %
\]
Furthermore, a slight modification of the proof of
Proposition~\ref{super-keyest}(ii) shows
\[
{\hat{G}^g}1_B\llvert\Delta{\hat{\pi}}\rrvert
\Gamma(x,z) \leq C(\log L)^{-5/2}\Gamma(x,z). %
\]
Together with $\pi\preceq C\Gamma$ and $\pi(x,z)\geq c(d,\eta
)L^{-(d-1)}$ this yields the bound
\[
{\hat{G}^g}1_B\llvert\Delta{\hat{\pi}}\rrvert
\pi(x,W_{t+6\ell}) \leq C(\log L)^{-5/2}\Gamma(x,W_t)
\leq C(\log L)^{-5/2}\pi(x,W_t). %
\]
To obtain~(\ref{local-eq3}), it remains to handle the second summand
of~(\ref{local-eq4}); that is, we have to bound
\[
\bigl\llvert{\hat{G}^g}1_{\sh_L(5r)}\Delta\pi{\hat{
\pi}}_\psi\bigr\rrvert\bigl(x,U_\ell(W_t)
\bigr). %
\]
We abbreviate $S= \sh_L(5r)$ and split into
\begin{eqnarray*}
&& {\hat{G}^g}1_S\Delta\pi{\hat{
\pi}}_\psi\bigl(x,U_\ell(W_t)\bigr)
\\
&&\qquad = \sum_{y\in
S}{\hat{G}^g}(x,y)\sum
_{z\in\partial
V_L}\Delta\pi(y,z) \bigl({\hat{
\pi}}_\psi\bigl(z,U_\ell(W_t)\bigr)-{\hat{
\pi}}_\psi\bigl(y,U_\ell(W_t)\bigr) \bigr)
\\
&&\qquad =\sum_{y\in S}{\hat{G}^g}(x,y)\sum
_{z\in
W_{t+6\ell}}\Delta\pi(y,z) \bigl({\hat{
\pi}}_\psi\bigl(z,U_\ell(W_t)\bigr)-{\hat{
\pi}}_\psi\bigl(y,U_\ell(W_t)\bigr) \bigr)
\\
&&\quad\qquad{} -\sum_{y\in S}{\hat{G}^g}(x,y)\sum
_{z\in\partial
V_L\setminus
W_{t+6\ell}}\Delta\pi(y,z){\hat{\pi}}_\psi
\bigl(y,U_\ell(W_t)\bigr).
\end{eqnarray*}
First note that since ${\hat{\pi}}_\psi(y,z')=0$ if
$\llvert y-z'\rrvert > 2\ell$,
\begin{eqnarray*}
&& \biggl\llvert\sum_{y\in S}{
\hat{G}^g}(x,y)\sum_{z\in\partial
V_L\setminus W_{t+6\ell}} \Delta
\pi(y,z){\hat{\pi}}_\psi\bigl(y,U_\ell(W_t)
\bigr)\biggr\rrvert
\\
&&\qquad \leq{\hat{G}^g}1_{U_{2\ell}(W_t)\cap S}\llvert\Delta\pi\rrvert(x,
\partial V_L\setminus W_{t+6\ell}).
\end{eqnarray*}
For $y\in U_{2\ell}(W_t)\cap S$, we apply Lemma~\ref{hittingprob}(iii) together
with Lemma~\ref{hittingprob-technical} and obtain
\begin{eqnarray*}
\llvert\Delta\pi\rrvert(y,\partial V_L\setminus
W_{t+6\ell}) &\leq&\sup_{y'\dvtx
\dist(y',U_{2\ell}(W_t)\cap S)\leq r}\pi\bigl(y',
\partial V_L\setminus W_{t+6\ell}\bigr)
\\
&\leq& C
\frac{r}{\ell}
\leq C(\log L)^{-13/2}.
\end{eqnarray*}
Since ${\hat{G}^g}\preceq\Gamma$ and $\pi(x,z)\geq cL^{-d-1}$,
${\hat{G}^g}(x,U_{2\ell}(W_t)\cap S)\leq C\pi(x,W_t)$, and thus
\[
{\hat{G}^g}1_{U_{2\ell}(W_t)\cap S}\llvert\Delta\pi\rrvert(x,\partial
V_L\setminus W_{t+6\ell}) \leq C(\log L)^{-13/2}
\pi(x,W_t). %
\]
It remains to bound
\[
\biggl\llvert\sum_{y\in S}{\hat{G}^g}(x,y)
\sum_{z\in
W_{t+6\ell}}\Delta\pi(y,z) \bigl({\hat{
\pi}}_\psi\bigl(z,U_\ell(W_t)\bigr)-{\hat{
\pi}}_\psi\bigl(y,U_\ell(W_t)\bigr) \bigr)
\biggr\rrvert. %
\]
Set $D^1(y) = \{z\in W_{t+6\ell}\dvtx  \llvert z-y \rrvert \leq\ell(\log
\ell)^{-4}\}
$. If
$D^1(y)\neq\varnothing$, then $\dist(y,W_t)\leq7\ell$. Using
Lemma~\ref{app-kernelest}
for the difference of the smoothing steps and the standard estimate
on ${\hat{G}^g}$,
\begin{eqnarray*}
&& \sum_{y\in S}{\hat{G}^g}(x,y)\sum
_{z\in
D^1(y)}\bigl\llvert\Delta\pi(y,z)\bigr\rrvert\bigl
\llvert{\hat{\pi}}_\psi\bigl(z,U_\ell(W_t)
\bigr)-{\hat{\pi}}_\psi\bigl(y,U_\ell(W_t)
\bigr)\bigr\rrvert
\\
&&\qquad \leq C\frac{t^{d-1}}{L^{d-1}}(\log\ell)^{-3}\leq C(\log L)^{-5/2}
\pi(x,W_t).
\end{eqnarray*}
The region $W_{t+6\ell}\setminus D^1(y)$ we split into
\begin{eqnarray*}
B_0(y) &=& \bigl\{z\in W_{t+6\ell}\dvtx  \llvert z-y \rrvert\in\bigl(
\ell(\log\ell)^{-4}, t\bigr]\bigr\},
\\
B_i(y) &=& \bigl\{z\in W_{t+6\ell}\dvtx  \llvert z-y \rrvert\in\bigl(it,
(i+1)t\bigr]\bigr\},\qquad i=1,2,\ldots, \lfloor2L/t\rfloor.
\end{eqnarray*}
Moreover, let
\[
S_i = \bigl\{y\in S\dvtx  B_i(y) \neq\varnothing\bigr\},\qquad
i=0,1,\ldots, \lfloor2L/t\rfloor. %
\]
Then
\[
\sum_{y\in S}{\hat{G}^g}(x,y)\sum
_{z\in W_{t+6\ell}\setminus
D^1(y)}\bigl\llvert\Delta\pi(y,z)\bigr\rrvert\leq C
\sum_{i=0}^{\lfloor2L/t\rfloor}{\hat{G}^g}(x,S_i)
\sup_{y\in
S_i}\llvert\Delta\pi\rrvert\bigl(y,B_i(y)
\bigr). %
\]
If $i\geq1$ and $y\in S_i$, then by Lemma~\ref{hittingprob}(iii),
\[
\llvert\Delta\pi\rrvert\bigl(y,B_i(y)\bigr) \leq\sup
_{y'\dvtx \llvert y'-y \rrvert \leq r}\pi\bigl(y',B_i(y)\bigr) \leq C
\frac{rt^{d-1}}{(it)^d}\leq C\frac{r}{i^dt}, %
\]
while in the case $i=0$, using the same lemma and additionally
Lemma~\ref{hittingprob-technical},
\begin{eqnarray*}
\sup_{y'\dvtx \llvert y'-y \rrvert \leq r}\pi\bigl(y',B_i(y)\bigr)&\leq& C r\sum_{z\in\partial
V_L}\frac{1}{((1/2)\ell(\log
\ell)^{-4} + \llvert y-z \rrvert )^d}
\\
&\leq& C\frac{r(\log\ell)^4}{\ell}\leq C (\log L)^{-5/2}.
\end{eqnarray*}
For the Green's function, we use the estimates
\[
{\hat{G}^g}(x,S_0) \leq C\frac{t^{d-1}}{L^{d-1}},\qquad{
\hat{G}^g}\biggl(x,\bigcup_{i
\geq(1/10)L/t}S_i
\biggr) \leq C, %
\]
while for $i=1,2,\ldots,\lfloor(1/10)L/t\rfloor$, it holds that
$\llvert S_i \rrvert \leq Cr(it)^{d-2}t$,
whence
\[
{\hat{G}^g}(x,S_i) \leq C \frac{i^{d-2}t^{d-1}}{L^{d-1}}.
\]
Altogether, we obtain
\begin{eqnarray*}
&& \sum_{i=0}^{\lfloor2L/t\rfloor}{
\hat{G}^g}(x,S_i)\sup_{y\in
S_i}\llvert
\Delta\pi\rrvert\bigl(y,B_i(y)\bigr)
\\
&&\qquad \leq C \Biggl((\log L)^{-5/2}\frac{t^{d-1}}{L^{d-1}} + \Biggl(
\frac{r}{t}\frac{t^{d-1}}{L^{d-1}}\sum_{i=1}^{\lfloor
(1/10)L/t\rfloor}
\frac{1}{i^2} \Biggr) +\frac{t^{d-1}}{L^{d-1}}\frac{r}{L} \Biggr)
\\
&&\qquad \leq C(\log L)^{-5/2}\frac{t^{d-1}}{L^{d-1}}.
\end{eqnarray*}
This finishes the proof of part~(ii).
\end{pf*}

Let us finally show how to obtain transience of the RWRE.
\begin{pf*}{Proof of Corollary~\ref{main-transience}}
As in the proof of Theorem~\ref{local-thm-exitmeas}, we assume that
$\ctwo(\delta,L_0,L)$ holds true for all $L\geq L_0$.
Fix numbers $\rho\geq3$, $\alpha\in(0,(4\rho)^{-1})$ to be specified
below. With these parameters and $n\geq1$, we put
$p=p_{\alpha\rho^n}$ and choose $\psi={ (m_x )}_{x\in
\mathbb{Z}^d}$ constant in
$x$, namely $m_x= \alpha\rho^n$.
Define
\[
A_n= \bigcap_{\llvert x \rrvert \,\leq\,
\rho^{n^{3/2}}}\bigcap
_{t\,\in\,[\alpha\rho^n,2\alpha\rho^n]} \bigl\{ \bigl\llVert\bigl(\Pi_{V_t(x)}-
\pi^{(p)}_{V_t(x)} \bigr)\hat{\pi}^{(p)}_{\psi
}(x,
\cdot)\bigr\rrVert_1 \leq(\log t)^{-9} \bigr\}.
\]
Under $\ctwo(\delta,L_0,L)$, we then have
\[
\pP\bigl(A_n^c \bigr) \leq C\alpha^d
\rho^{(d+1)n^{3/2}}\exp\bigl(- \bigl(\log\bigl(\alpha\rho^n \bigr)
\bigr)^2 \bigr). %
\]
Therefore, for any choice of $\alpha, \rho$ it holds that
\[
\sum_{n=1}^\infty\pP\bigl(A_n^c
\bigr) < \infty, %
\]
whence by Borel--Cantelli,
%
%
\begin{equation}
\label{main-transience-bc} \pP\Bigl(\liminf_{n\rightarrow\infty} A_n
\Bigr)= 1.
\end{equation}
We set $q_{n,\alpha,\rho}= \hat{\pi}^{(p)}_\psi$ and denote the
coarse-grained RWRE
transition kernel by
\[
Q_{n,\alpha,\rho}(x,\cdot) = \frac{1}{\alpha\rho^n}\int_{\mathbb{R}_+}
\varphi\biggl(\frac
{t}{\alpha\rho^n} \biggr)\Pi_{V_t(x)}(x,\cdot)\,\dt t.
\]
If $\llvert x \rrvert \leq\rho^{n^{3/2}}$, we have on $A_n$
\[
\bigl\llVert(Q_{n,\alpha,\rho}-q_{n,\alpha,\rho} )q_{n,\alpha,\rho
}(x,\cdot)
\bigr\rrVert_1 \leq\bigl(\log\bigl(\alpha\rho^n\bigr)
\bigr)^{-9} \leq C(\alpha,\rho)n^{-9}. %
\]
Now assume $\llvert x \rrvert \leq\rho^n+1$. On $A_n$, for $N$ fixed, it follows that for $1\leq M \leq N$,
%
%
\begin{equation}
\label{main-transience-1} \bigl\llVert\bigl( (Q_{n,\alpha,\rho} )^M-
(q_{n,\alpha
,\rho} )^M \bigr)q_{n,\alpha,\rho}(x,\cdot)\bigr\rrVert
_1 \leq C(\alpha,\rho)Mn^{-9}.
\end{equation}
For fixed $\omega$, let $ (\xi_k )_{k\geq0}$ be the
Markov chain
running with transition kernel
$Q_{n,\alpha,\rho}=Q_{n,\alpha,\rho}(\omega)$. Clearly,
$ (\xi_k )_{k\geq0}$ can be obtained by observing the
basic RWRE
$ (X_k )_{k\geq0}$ at randomized stopping times. Then
\begin{eqnarray*}
&& \Prw_{x,\omega} (\xi_{N-1} \in V_{\rho^{n+1}+2\alpha\rho^n} )
\\
&&\qquad \leq(Q_{n,\alpha,\rho} )^{N-1}q_{n,\alpha,\rho} (x,
V_{\rho^{n+1}+4\alpha\rho^n} )
\\
&&\qquad \leq\bigl\llVert\bigl( (Q_{n,\alpha,\rho} )^{N-1}-
(q_{n,\alpha,\rho} )^{N-1} \bigr)q_{n,\alpha,\rho}(x,\cdot)\bigr
\rrVert_1 + (q_{n,\alpha,\rho} )^N(x,V_{2\rho^{n+1}}).
\end{eqnarray*}
Using Proposition~\ref{super-localclt}, we can find
$N=N(\alpha,\rho)\in\mathbb{N}$, depending not on $n$, such that
for any
$x$ with $\llvert x \rrvert \leq\rho^n+1$, it holds that
${(q_{n,\alpha,\rho})}^N(x,V_{2\rho^{n+1}}) \leq
1/10$. With~(\ref{main-transience-1}), we conclude that for such $x$,
$n\geq n_0(\alpha,\rho,N)$ large enough and $\omega\in A_n$,
%
%
\begin{equation}
\label{main-transience-2} \Prw_{x,\omega} (\xi_{N-1} \in V_{\rho
^{n+1}+2\alpha\rho
^n}
) \leq C(\alpha,\rho)Nn^{-9} + 1/10 \leq1/5.
\end{equation}
On the other hand, if $x$ is outside $V_{\rho^{n-1}+2\alpha\rho^n}$,
\begin{eqnarray*}
&& \Prw_{x,\omega} (\xi_{M} \in V_{\rho^{n-1}+2\alpha
\rho^n}\mbox{ for some }0\leq M\leq N-1 )
\\
&&\qquad \leq\sum_{M=1}^{N-1} (Q_{n,\alpha,\rho}
)^Mq_{n,\alpha
,\rho} (x, V_{\rho^{n-1}+4\alpha\rho^n} )
\\
&&\qquad \leq\sum_{M=1}^{N-1}\bigl\llVert\bigl(
(Q_{n,\alpha,\rho} )^M- (q_{n,\alpha,\rho} )^M
\bigr)q_{n,\alpha,\rho
}(x,\cdot)\bigr\rrVert_1 + \sum
_{k=2}^{N} (q_{n,\alpha,\rho} )^k(x,V_{2\rho^{n-1}}).
\end{eqnarray*}
If $\rho^n-1 \leq\llvert x \rrvert $, then
$ (q_{n,\alpha,\rho} )^k(x,V_{2\rho^{n-1}}) = 0$ as long
as $k
\leq(1-3/\rho)/(2\alpha)$. By first choosing $\rho$ large enough, then
$\alpha$ small enough and estimating the higher summands again
with Proposition~\ref{super-localclt}, we deduce that for such $x$ and
all large $n$,
\[
\sum_{k=1}^{\infty} (q_{n,\alpha,\rho}
)^k(x,V_{2\rho
^{n-1}}) \leq1/10. %
\]
Together with~(\ref{main-transience-1}), we have for large $n$,
$\omega
\in A_n$ and $\rho^n-1 \leq\llvert x \rrvert \leq\rho^n+1$,
%
%
\begin{eqnarray}
\label{main-transience-3} && \Prw_{x,\omega} (\xi_{M} \in
V_{\rho^{n-1}+2\alpha
\rho^n}\mbox{ for some } 0\leq M\leq N-1 )
\nonumber\\[-8pt]\\[-8pt]\nonumber
&&\qquad \leq C(\alpha,\rho)N^2n^{-9} + 1/10 \leq1/5.
\end{eqnarray}
Let\vspace*{1pt} $B$ be the event that the walk $ (\xi_k )_{k\geq
0}$ leaves $V_{\rho^{n+1}+2\alpha\rho^n}$ before reaching
$V_{\rho^{n-1}+2\alpha\rho^n}$. From~(\ref{main-transience-2})
and~(\ref{main-transience-3}) we deduce that
$\Prw_{x,\omega} (B ) \geq3/5$, provided $n$ is large enough,
$\omega\in A_n$ and $\rho^n-1\leq\llvert x \rrvert \leq\rho^n +1$.
But on $B$, the
underlying basic\vspace*{2pt} RWRE $ (X_k )_{k\geq0}$ clearly leaves
$V_{\rho^{n+1}}$
before reaching $V_{\rho^{n-1}}$; see Figure~\ref{fig6}. Hence if $\omega\in\liminf A_n$,
there exists $m_0=m_0(\omega)\in\mathbb{N}$ such that
\[
\Prw_{x,\omega} (\tau_{V_{\rho^{n+1}}}< T_{V_{\rho
^{n-1}}} ) \geq3/5
\]
for all $n\geq m_0$, $x$ with $\llvert x \rrvert \geq\rho^n-1$ (of
course, we may
now drop the
constraint $\llvert x \rrvert \leq\rho^n+1$). From this property,
transience easily follows.
Indeed, for $m$, $M$, $k \in\mathbb{N}$ satisfying $M > m\geq m_0$ and
$0\leq
k \leq M+1-m$, set
\[
h_M(k) = \sup_{x\dvtx  \llvert x \rrvert \geq\rho^{m+k}-1}\Prw_{x,\omega}
(T_{V_{\rho^m}}< \tau_{V_{\rho^M}} ). %
\]

%
\begin{figure}

\includegraphics{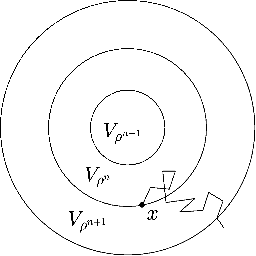}

\caption{On a set of environments with mass $1$, the RWRE started at
any $x$
with $\llvert x \rrvert \geq\rho^n -1$ leaves the ball $V_{\rho
^{n+1}}$ before hitting
$V_{\rho^{n-1}}$ with probability at least $3/5$. This implies transience
of the RWRE.}\label{fig6}
\end{figure}

\noindent
Then $h_M$ solves the difference inequality
\[
h_M(k) \leq\tfrac{2}{5}h_M(k-1) +
\tfrac{3}{5}h_M(k+1) %
\]
with boundary conditions $h_M(0) = 1$, $h_M(M+1-m) = 0$. By
either applying a discrete maximum principle or by a direct
computation, we
see that $h_M \leq\overline{h}_M$, where $\overline{h}_M$
is the solution of the difference equality
%
%
\begin{equation}
\label{main-transience-4} \overline{h}_M(k) = \tfrac{2}{5}
\overline{h}_M(k-1) + \tfrac{3}{5}\overline{h}_M(k+1)
\end{equation}
with boundary conditions $\overline{h}_M(0) = 1$, $\overline
{h}_M(M+1-m)=0$. Solving
(\ref{main-transience-4}), we get
\[
\overline{h}_M(k) = \frac{1}{1- (3/2 )^{M+1-m}} + \frac{1}{1- (2/3 )^{M+1-m}}
\biggl(\frac{2}{3} \biggr)^k. %
\]
Letting $M\rightarrow\infty$, we deduce that for $\llvert x \rrvert
\geq\rho^{m+k}$,
%
%
\begin{equation}
\Prw_{x,\omega} (T_{V_{\rho^m}}< \infty) \leq\lim_{M\rightarrow\infty}
\overline{h}_M(k) = \biggl(\frac
{2}{3} \biggr)^k.
\end{equation}
Together with~(\ref{main-transience-bc}), this proves that for almost all
$\omega\in\Omega$, the random walk is transient under $\Prw_{\cdot
,\omega}$.
\end{pf*}

\setcounter{equation}{0}
\begin{appendix}\label{appe}
\section*{Appendix}
\subsection{Some difference estimates}
In this section we collect some difference estimates of (non)smoothed exit
distributions which are mainly needed to prove Lemma~\ref{est-phi}(i)
and (iii). The first technical lemma connects the exit distributions\vadjust{\goodbreak} of a
symmetric random walk with \mbox{one-}step distribution $p\in{\mathcal
{P}^{\mathrm{s}}_\kappa}$ to
those of
Brownian motion with covariance matrix~$\Lambda_p$. As always, $\kappa>0$
can be chosen arbitrarily small.
%
\begin{lemma}
\label{app-kernelest-technical}
Let $p\in\mathcal{P}^{\mathrm{s}}_\kappa$, and let $\beta,\eta
>0$ with $3\eta<\beta<1$.
For large $L$, there
exists a constant $C > 0$ such that for
$A\subset\mathbb{R}^d$, $A^\beta= \{y\in\mathbb{R}^d\dvtx
\dist(y,A)\leq L^\beta\}$ and $x\in V_L$ with $\dL(x)> L^\beta$,
the following holds:
\begin{longlist}[(ii)]
\item$\pi^{(p)}_L(x,A) \leq
\pi^{\mathrm{B}(p)}_L (x,A^\beta) (1+CL^{-(\beta
-3\eta)}
) + L^{-(d+1)}$.\vspace*{1pt}
\item$\pi^{\mathrm{B}(p)}_L(x,A) \leq\pi^{(p)}_L(x,A^{\beta
}) (1+CL^{-(\beta
-3\eta)} ) + L^{-(d+1)}$.
\end{longlist}
\end{lemma}
\begin{pf}(i) Set $L' = L + L^\eta$, $L''= L +2L^\eta$, and denote by
$A_{'}^\beta$
the image of $\{y\in\partial C_L\dvtx  \dist(y,A)\leq L^\beta/2\}$ on
$\partial C_{L'}$ under the map $y\mapsto(L'/L)y$. We write $\pi
^{\mathrm{B}}_L$
for $\pi^{\mathrm{B}(p)}_L$. Then, denoting by $\mathrm{P}^{\mathrm
{B}}_x$ the law of a $d$-dimensional
Brownian motion $W_t$ with covariance $\Lambda_p$, conditioned on
$W_0=x$, and by $\taubm_L =\inf\{t\geq0\dvtx  W_t \notin
C_{L} \}$
the first exit time from $C_L$,
\[
\pi^{\mathrm{B}}_{L'}\bigl(x,A_{'}^\beta
\bigr)\leq\pi^{\mathrm
{B}}_{L}\bigl(x,A^\beta\bigr) +
\mathrm{P}^{\mathrm{B}}_x \bigl(W_{\taubm_{L'}}\in
A_{'}^\beta, W_{\taubm_L} \notin A^\beta
\bigr). %
\]
Let us first assume that we have already proved
%
%
\begin{equation}
\label{app-kernelest-technical-eq0} \mathrm{P}^{\mathrm{B}}_x \bigl
(W_{\taubm_{L'}}
\in A_{'}^\beta, W_{\taubm_L} \notin A^\beta
\bigr)\leq CL^{\eta-\beta}\pi^{\mathrm
{B}}_{L'}
\bigl(x,A_{'}^\beta\bigr),
\end{equation}
so that
%
%
\begin{equation}
\label{app-kernelest-technical-eq1} \pi^{\mathrm{B}}_L \bigl(x,A^\beta
\bigr) \geq\pi^{\mathrm
{B}}_{L'}\bigl(x,A_{'}^\beta
\bigr) \bigl(1-CL^{-(\beta-\eta)} \bigr).
\end{equation}
As a consequence of~\cite{ZAI}, Theorem~2, for
each $k\in\mathbb{N}$ there exists a positive constant $C_1=C_1(k)$
such that
for each integer $n\geq1$, one can construct on the same probability space
a Brownian motion $W_t$ with covariance matrix $d^{-1}\Lambda_p$ as
well as
a symmetric random walk $X_n$ with one-step probability $p$, both
starting in $x$ and satisfying (with
$\mathbb{Q}$ denoting the probability measure on that space)
%
%
\begin{equation}
\label{app-kernelest-technical-eq2} \mathbb{Q} \Bigl(\max_{0\leq m\leq
n}\llvert
X_m-W_m\rrvert> C_1\log n \Bigr) \leq
C_1 n^{-k}.
\end{equation}
Choose $k >(2/5)(d+1)$, and let $C_1(k)$ be the corresponding constant. The
following arguments hold for sufficiently large $L$. By standard
results on the
oscillation of Brownian paths,
%
%
\begin{equation}
\label{app-kernelest-technical-eq3} \mathbb{Q} \Bigl(\sup_{0\leq t\leq
L^{5/2}}\llvert
W_{\lfloor
t\rfloor}-W_t\rrvert> (5/2)C_1\log L \Bigr)
\leq(1/3)L^{-(d+1)}.
\end{equation}
With
\[
B_1= \Bigl\{\sup_{0\leq t\leq L^{5/2}}\llvert
X_{\lfloor
t\rfloor}-W_t\rrvert\leq5C_1\log L \Bigr\},
\]
we deduce from~(\ref{app-kernelest-technical-eq2})
and~(\ref{app-kernelest-technical-eq3}) that
\[
\mathbb{Q} \bigl({B}^c_1 \bigr) \leq(2/3)L^{-(d+1)}.
\]
Let $B_2= \{\taubm_{L'}\vee\tau_{L''} \leq L^{5/2}\}$.
We claim that
%
%
\begin{equation}
\label{app-kernelest-technical-eq4} \mathbb{Q} \bigl({B}^c_2 \bigr) \leq
(1/3)L^{-(d+1)}.
\end{equation}
By the central limit theorem, one finds a constant $c>0$ with
\[
\mathbb{Q} \bigl(\tau_{L''} \leq{\bigl(L''
\bigr)}^2 \bigr)\geq c\qquad\mbox{for }L\mbox{ large}.
\]
By the Markov property, we obtain $\mathbb{Q} (\tau_{L''} >
L^{5/2} )
\leq(1-c)^{L^{1/3}}$. The probability $\mathbb{Q} (\taubm_{L'} >
L^{5/2} )$ decays in the same way, and~(\ref
{app-kernelest-technical-eq4})
follows. Since $\pi^{\mathrm{B}}_L$ is unchanged if the Brownian
motion is
replaced by
a Brownian motion with covariance $d^{-1}\Lambda_p$, we have ($\Prw_x$
denotes the law of $X_n$)
%
%
\begin{eqnarray}
\label{app-kernelest-technical-eq5}
\quad&&  \pi^{\mathrm{B}}_{L'}\bigl
(x,A_{'}^\beta \bigr)\nonumber
\\
&&\qquad \geq \mathbb{Q} \bigl(X_{\tau_L}\in A, W_{\taubm_{L'}}\in
A_{'}^\beta\bigr)
\\
&&\qquad \geq \Prw_x (X_{\tau_L}\in A ) - \mathbb{Q}
\bigl(X_{\tau_L}\in A, W_{\taubm_{L'}}\notin A_{'}^\beta,
B_1\cap B_2 \bigr)
- L^{-(d+1)}.\nonumber
\end{eqnarray}
Let $U= \{z\in\mathbb{Z}^d\dvtx  \dist(z,(\partial C_{L'}\setminus
A_{'}^\beta)) \leq5C_1\log L \}$. Then
\[
\mathbb{Q} \bigl(X_{\tau_L}\in A, W_{\taubm_{L'}}\notin
A_{'}^\beta, B_1\cap B_2 \bigr)
\leq\Prw_x (X_{\tau_L}\in A, T_U < \tau
_{L''} ). %
\]
By the strong Markov property,
\[
\Prw_x (X_{\tau_L}\in A, T_U <
\tau_{L''} )\leq\Prw_x (X_{\tau_L}\in A )\sup
_{y\in A}\Prw_y ( T_U <
\tau_{L''} ). %
\]
Furthermore, there exists a constant $c>0$ such that for $y\in A$ and
$z\in
U$, we have $\llvert y-z \rrvert \geq cL^{\beta}$ and $\dist_{L''}(z)
\leq\dist_{L''}(y)
\leq2L^\eta$. Therefore, an application of first
Lemma~\ref{hittingprob}(ii) and then
Lemma~\ref{hittingprob-technical} yields
\[
\Prw_y (T_U < \tau_{L''} ) \leq C
L^{2\eta}\sum_{z\in
U}\frac{1}{{\llvert y-z \rrvert }^d} \leq
CL^{2\eta}(\log L) L^{-\beta}\leq CL^{-(\beta-3\eta)}, %
\]
uniformly in $y\in A$. Going back
to~(\ref{app-kernelest-technical-eq5}), we arrive at
\[
\pi^{\mathrm{B}}_{L'} \bigl(x,A_{'}^\beta
\bigr) \geq\pi_L(x,A) \bigl(1-CL^{-(\beta-3\eta)}
\bigr)-L^{-(d+1)}. %
\]
Together with~(\ref{app-kernelest-technical-eq1}), this shows (i), but
we still have to prove~(\ref{app-kernelest-technical-eq0}). First, by
Lemma~\ref{exitdensitybm}(i) (all integrals are surface integrals)
\begin{eqnarray*}
&& \mathrm{P}^{\mathrm{B}}_x \bigl(W_{\taubm_{L'}}\in
A_{'}^\beta, W_{\taubm_L} \notin A^\beta
\bigr)
\\
&&\qquad \leq\int_{\partial C_L\setminus
A^\beta}\pi^{\mathrm{B}}_{L}(x,dy)
\pi^{\mathrm{B}}_{L'}\bigl(y,A^\beta\bigr)
\\
&&\qquad \leq C\dL(x)L^\eta\int_{\partial C_L\setminus
A^\beta}\frac{1}{\llvert x-y\rrvert ^d}
\int_{A_{'}^\beta}\frac{1}{\llvert y-z \rrvert ^d}\,\dt z\,\dt y.
\end{eqnarray*}
Fix $z\in A_{'}^\beta\subset\partial C_{L'}$ and put
\begin{eqnarray*}
D_1 &=& \bigl\{y\in\partial C_L\setminus
A^\beta\dvtx  \llvert y-z \rrvert>\llvert x-z \rrvert/2 \bigr\},
\\
D_2&=&\partial C_L\setminus\bigl(A^\beta\cup
D_1\bigr). %
\end{eqnarray*}
Then, using $\dL(x) \geq L^{\beta}$ and
Lemma~\ref{hittingprob-technical} in the last step,
\[
\int_{D_1}\frac{1}{\llvert x-y\rrvert ^d}\frac{1}{\llvert y-z
\rrvert ^d}\,\dt y \leq
\frac{C}{\llvert x-z \rrvert ^d}\int_{D_1}\frac{1}{\llvert
x-y\rrvert ^d}\,\dt y\leq
\frac
{CL^{-\beta}}{\llvert x-z \rrvert ^d}. %
\]
For the integral over $D_2$ we obtain the same bound, by using that
$\llvert x-y\rrvert \geq(1/2)\llvert x-z \rrvert $ for $y\in D_2$,
the fact that $\llvert y-z \rrvert \geq
cL^\beta$ if $y\in\partial C_L\setminus A^\beta$, $z\in
A_{'}^\beta$ and Lemma~\ref{hittingprob-technical}.
Altogether,
\[
\mathrm{P}^{\mathrm{B}}_x \bigl(W_{\taubm_{L'}}\in
A_{'}^\beta, W_{\taubm_L} \notin A^\beta
\bigr)\leq\frac{C\dL(x)L^{\eta-\beta}}{\llvert x-z \rrvert ^d}. %
\]
Integrating over $z\in A_{'}^\beta$, we obtain with the lower bound of
Lemma~\ref{exitdensitybm}(i),
\[
\mathrm{P}^{\mathrm{B}}_x \bigl(W_{\taubm_{L'}}\in
A_{'}^\beta, W_{\taubm_L} \notin A^\beta
\bigr)\leq CL^{\eta-\beta}\pi^{\mathrm
{B}}_{L'}
\bigl(x,A_{'}^\beta\bigr), %
\]
as claimed.

(ii) One can follow the same steps, interchanging the role of Brownian
motion and the random walk.
\end{pf}

Again, let $p\in\mathcal{P}^{\mathrm{s}}_\kappa$. We write ${\hat
{\pi}^{\mathrm{B}(p)}}_{\psi}(x,z)$ for the
density of
${\hat{\pi}^{\mathrm{B}(p)}}_{\psi}(x,\dz z)$ with respect to
$d$-dimensional Lebesgue measure,
that is, for $\psi= (m_x) \in\mathcal{M}_L$,
%
%
\begin{equation}
\label{app-cgpsibmdens} {\hat{\pi}^{\mathrm{B}(p)}}_{\psi}(x,z)=
\frac{1}{m_x}\varphi\biggl(\frac
{\llvert z-x \rrvert }{m_x} \biggr)\pi^{\mathrm{B}(p)}_{C_{\llvert
z-x \rrvert }}(0,z-x).
\end{equation}
For ease of notation, we write in the following ${\hat{\pi}^{\mathrm
{B}}}_\psi$ for
${\hat{\pi}^{\mathrm{B}(p)}}_{\psi}$ and ${\hat{\pi}}_\psi$ for
$\hat{\pi}^{(p)}
_{\psi}$.
%
\begin{lemma}
\label{app-kernelest} Let $p\in\mathcal{P}^{\mathrm{s}}_\kappa$.
There exists a constant $C > 0$ such that for large $L$,
$\psi=(m_y)\in\mathcal{M}_L$, $x,x'\in\mathcal{U}_L\cap\mathbb
{Z}^d$ and any $z,z'\in\mathbb{Z}^d$:
\begin{longlist}[(vii)]
\item${\hat{\pi}}_{\psi}(x,z) \leq CL^{-d}$;\vspace*{3pt}
\item${\hat{\pi}^{\mathrm{B}}}_{\psi}(x,z)\leq CL^{-d}$;\vspace*{2pt}
\item$\llvert {\hat{\pi}}_{\psi}(x,z) - {\hat{\pi }}_{\psi}(x',z)\rrvert \leq C\llvert x-x' \rrvert L^{-(d+1)}\log L$;\vspace*{2pt}
\item$\llvert {\hat{\pi}}_{\psi}(x,z) - {\hat{\pi}}_{\psi}(x,z')\rrvert \leq C\llvert z-z' \rrvert L^{-(d+1)}\log L$;\vspace*{2pt}
\item$\llvert {\hat{\pi}^{\mathrm{B}}}_{\psi}(x,z) - {\hat{\pi }^{\mathrm{B}}}_{\psi}(x',z)\rrvert \leq C\llvert x-x' \rrvert L^{-(d+1)}\log L$;\vspace*{2pt}
\item$\llvert {\hat{\pi}^{\mathrm{B}}}_{\psi}(x,z) - {\hat{\pi }^{\mathrm{B}}}_{\psi}(x,z')\rrvert \leq C\llvert z-z' \rrvert L^{-(d+1)}\log L$;\vspace*{2pt}
\item$\llvert {\hat{\pi}}_{\psi}(x,z)-{\hat{\pi}^{\mathrm{B}}}_{\psi}(x,z)\rrvert \leq L^{-(d+1/4)}$.
\end{longlist}
\end{lemma}
%
%
\begin{corollary}
\label{app-kernelest-cor}
In the situation of the preceding lemma:
\begin{longlist}[(ii)]
\item
\begin{eqnarray*}
&& \bigl\llvert{\hat{\pi}}_{\psi}(x,z) - {\hat{
\pi}}_{\psi}\bigl(x',z\bigr)\bigr\rrvert
\\
&&\qquad \leq C\min\bigl\{\bigl\llvert x-x' \bigr\rrvert L^{-(d+1)}
\log L, \bigl\llvert x-x' \bigr\rrvert L^{-(d+1)}+
L^{-(d+1/4)} \bigr\};
\end{eqnarray*}
\item
\begin{eqnarray*}
&& \bigl\llvert{\hat{\pi}}_{\psi}(x,z) - {\hat{
\pi}}_{\psi}\bigl(x,z'\bigr)\bigr\rrvert
\\
&&\qquad \leq C\min\bigl\{\bigl\llvert z-z' \bigr\rrvert L^{-(d+1)}
\log L, \bigl\llvert z-z' \bigr\rrvert L^{-(d+1)}+
L^{-(d+1/4)} \bigr\}.
\end{eqnarray*}
\end{longlist}
\end{corollary}

\begin{pf}
Combine (iii)--(vii).
\end{pf}
\begin{pf*}{Proof of Lemma~\ref{app-kernelest}}
(i) and (ii) follow from
the definitions of ${\hat{\pi
}}_{\psi}$ and
${\hat{\pi}^{\mathrm{B}}}_{\psi}$ together with part (i) of
Lemma~\ref{lemmalawler} and Lemma~\ref{exitdensitybm}, respectively.

(iii)~and~(iv) We can restrict ourselves to the case
$\llvert x-x' \rrvert = 1$, as otherwise we take a shortest path
connecting $x$ with $x'$
inside $\mathcal{U}_L\cap\mathbb{Z}^d$, and apply the result for
distance one
$O (\llvert x-x' \rrvert )$ times. We have
\begin{eqnarray*}
&& {\hat{\pi}}_{\psi}(x,z)-{\hat{\pi}}_{\psi}
\bigl(x',z\bigr)
\\
&&\qquad = \biggl(1-\frac{m_x}{m_{x'}} \biggr){\hat{\pi}}_{\psi
}(x,z) +
\frac{1}{m_{x'}}\int_{\mathbb{R}_+} \biggl(\varphi\biggl(
\frac
{t}{m_x} \biggr) -\varphi\biggl(\frac{t}{m_{x'}} \biggr) \biggr)
\pi_{V_t(x)}(x,z)\,\dt t
\\
&&\quad\qquad{}+ \frac{1}{m_{x'}}\int_{\mathbb{R}_+}\varphi\biggl(
\frac
{t}{m_{x'}} \biggr) \bigl(\pi_{V_t(x)}(x,z)-\pi_{V_t(x')}
\bigl(x',z\bigr) \bigr)\,\dt t
\\
&&\qquad = I_1 + I_2 + I_3.
\end{eqnarray*}
Using the fact that $\psi\in\mathcal{M}_L$ and part~(i) for
${\hat{\pi}}_{\psi}(x,z)$, it follows that $\llvert I_1\rrvert \leq
CL^{-(d+1)}$. Using
additionally the smoothness of $\varphi$ and, by
Lemma~\ref{lemmalawler}(i), $\llvert \pi_{V_t(x)}(x,z)\rrvert \leq
CL^{-(d-1)}$,
we also have $\llvert I_2\rrvert \leq CL^{-(d+1)}$. It remains to
handle~$I_3$. By
translation invariance of the random walk, $\pi_{V_t(x)}(x,z) =
\pi_{V_t}(0,z-x)$. In particular, both (iii) and (iv) will follow if we
prove that
%
%
\begin{equation}
\label{a4-1} \biggl\llvert\int_{\mathbb{R}_+}\varphi\biggl(
\frac{t}{m_{x}} \biggr) \bigl(\pi_{V_t}(0,z-x)-\pi_{V_t}
\bigl(0,z-x'\bigr) \bigr)\,\dt t \biggr\rrvert\leq C L^{-d}
\log L
\end{equation}
for $x,x'$ with $\llvert x-x' \rrvert = 1$. By definition of $\mathcal{M}_L$,
$m_x\in
(L/10,5L)$. We may therefore assume that $L/10 < \llvert y-z \rrvert
< 10L$ for
$y=x,x'$. Due to the smoothness of $\varphi$ and the fact that the
integral is over an interval of length at most $2$,~(\ref{a4-1}) will
follow if we show
\[
\biggl\llvert\int_{L/10}^{10L} \bigl(
\pi_{V_t}(0,z-x)-\pi_{V_t}\bigl(0,z-x'\bigr)
\bigr)\,\dt t\biggr\rrvert\leq C L^{-d}\log L. %
\]
We set $J = \{t > 0\dvtx  z-x\in\partial V_t\}$ and $J' = \{t>0\dvtx
z-x'\in\partial V_{t'}\}$, where
\[
t' = t'(t) = \biggl\llvert t\frac{(z-x)}{\llvert z-x \rrvert }-
\bigl(x'-x\bigr)\biggr\rrvert. %
\]
$J$ is an interval\vspace*{1pt} of length at most $1$, and $J'$ has the same length up
to order $O(L^{-1})$. Furthermore, $\llvert (J\setminus J')\cup
(J'\setminus
J)\rrvert $ is of order
$O(L^{-1})$, and $\llvert \frac{\dt }{\dt t}t'\rrvert =
1 +
O(L^{-1})$. Using that both $\pi_{V_t}(0,z-x)$ and $\pi_{V_t}(0,z-x')$
are of order $O(L^{-(d-1)})$, it therefore suffices to prove
%
%
\begin{equation}
\label{a4-2} \biggl\llvert\int_{J\cap
J'} \bigl(
\pi_{V_t(x)}(x,z)-\pi_{V_{t'}(x')}\bigl(x',z\bigr) \bigr)\,\dt t\biggr\rrvert\leq CL^{-d}\log L.
\end{equation}
Write $V$ for $V_t(x)$ and $V'$ for $V_{t'}(x')$. By a first exit
decomposition,
\[
\pi_V(x,z) \leq\pi_{V'}(x,z) + \sum
_{y\in V\setminus V'}\Prw_{x,p} (T_y <
\tau_V )\pi_V(y,z). %
\]
By Lemma~\ref{lemmalawler}(ii), we can replace $\pi_{V'}(x,z)$ by
$\pi_{V'}(x',z) + O(L^{-d})$. For $y\in V\setminus V'$ we have by
Lemma~\ref{hittingprob}(ii) $\pi_V(y,z) = O(\llvert y-z \rrvert
^{-d})$ and
$\Prw_{x,p} (T_y < \tau_V ) = O(L^{-(d-1)})$, uniformly in $t
\in J\cap J'$. Using $\llvert x-x' \rrvert = 1$, we have with $r=
\llvert z-x \rrvert $
\[
\bigcup_{t\in J\cap J'}\bigl(V\setminus V'
\bigr) \subset V_r(x)\setminus V_{r-2}
\bigl(x'\bigr) \subset x + \sh_r(3), %
\]
and for any $y\in x+\sh_r(3)$, it follows by a geometric consideration that
\[
\int_{J\cap J'}1_{\{y\in V\setminus V'\}}\,\dt t \leq C\frac{\llvert
y-z \rrvert }{L}.
\]
Altogether, applying Lemma~\ref{hittingprob-technical} in the last
step,
\begin{eqnarray*}
&& \int_{J\cap J'}\pi_{V}(x,z)\,\dt t
\\
&&\qquad \leq\int_{J\cap J'}\pi_{V'}\bigl(x',z
\bigr)\,\dt t + O\bigl(L^{-d}\bigr) + CL^{-(d-1)}\sum
_{y\in x +
\sh_r(3)}\frac{1}{\llvert y-z \rrvert ^{d}}\frac{\llvert y-z
\rrvert }{L}
\\
&&\qquad \leq\int_{J\cap J'}\pi_{V'}\bigl(x',z
\bigr)\,\dt t + CL^{-d}\log L.
\end{eqnarray*}
The reverse inequality, proved in the same way, then implies~(\ref
{a4-2}).

(v) and (vi) are the analogous statements of (iii) and (iv) for Brownian
motion with covariance matrix $\Lambda_p$ and can be proved in the
same way.

(vii) Fix $\alpha= 2/3$, $\beta= 1/3$, and let $0<\eta<1/40$. Set $A =
C_{L^\alpha}(z)$ and $A^{\mathbb{Z}}= A\cap\mathbb{Z}^d$. By part~(iv), we
have
%
%
\begin{equation}
\label{a7-1} {\hat{\pi}}_{\psi}(x,z) \leq\frac{1}{\llvert A^{\mathbb
{Z}}\rrvert }{\hat{
\pi}}_{\psi} \bigl(x,A^{\mathbb{Z}} \bigr)+ CL^{-(d+1-\alpha)}\log L.
\end{equation}
Moreover,
%
%
\begin{equation}
\label{a7-2} {\hat{\pi}}_{\psi} \bigl(x,A^{\mathbb{Z}} \bigr) =
\frac{1}{m_x}\int_{L/10}^{10L} \varphi\biggl(
\frac{t}{m_x} \biggr)\pi_{V_t(x)} \bigl(x,A^{\mathbb
{Z}} \bigr)\,\dt
t.
\end{equation}
By Lemma~\ref{app-kernelest-technical}(i), it follows that for $t
\in(L/10,10L)$,
\[
\pi_{V_t(x)} \bigl(x,A^{\mathbb{Z}} \bigr) \leq\pi^{\mathrm{B}}_{V_t(x)}
\bigl(x,A^\beta\bigr) \bigl(1+CL^{-(\beta-3\eta
)} \bigr) +
CL^{-(d+1)}, %
\]
where $A^\beta= C_{L^\alpha+ L^\beta}(z)$ and the constant $C$ is
uniform in $t$. If we plug the last line into~(\ref{a7-2}), and use
parts (ii)~and~(vi), we arrive at
\begin{eqnarray*}
{\hat{\pi}}_\psi\bigl(x,A^{\mathbb{Z}} \bigr) &\leq&{\hat{
\pi}^{\mathrm{B}}}_\psi\bigl(x,A^\beta\bigr)
\bigl(1+CL^{-(\beta
-3\eta)} \bigr) + CL^{-(d+1)}
\\
&\leq&{\hat{\pi}^{\mathrm{B}}}_\psi(x,A) \bigl(1+CL^{-(\beta-3\eta
)}
\bigr) + CL^{-d}L^{(d-1)\alpha+\beta}
\\
&\leq&\llvert A \rrvert\cdot{\hat{\pi}^{\mathrm{B}}}_\psi(x,z) +
CL^{d\alpha
}L^{-(d+\beta-3\eta)}.
\end{eqnarray*}
Notice that in our notation, $\llvert A \rrvert $ is the volume of
$A$, while
$\llvert A^{\mathbb{Z}}\rrvert $ is
the cardinality of~$A^{\mathbb{Z}}$. From Gauss's circle problem we
have learned that $\llvert A \rrvert =
\llvert A^{\mathbb{Z}}\rrvert + O (L^{(d-1)\alpha} )$. Going back
to~(\ref{a7-1}),
this implies
\[
{\hat{\pi}}_{\psi}(x,z) \leq{\hat{\pi}^{\mathrm{B}}}_{\psi
}(x,z)
+ L^{-(d+1/4)}, %
\]
as claimed. To prove the reverse inequality, we can follow the same steps,
replacing the random walk estimates by those of Brownian motion and vice
versa.
\end{pf*}

\subsection{Proof of Lemma~\texorpdfstring{\protect\ref{est-phi}}{4.7}}
For simplicity, let us write $\phi$ for
$\phi_{L,p,\psi,q}$ and ${\phi}^{\mathrm{B}}$ for ${\phi}^{\mathrm
{B}}_{L,p,\psi,q}$.
\begin{pf*}{Proof of Lemma~\ref{est-phi}}
(i) Set $\alpha= 2/3$, $\beta=
1/3$ and $\eta= \dist(x,\partial
V_L)$. Choose $y_1\in\partial V_L$ such that $\llvert x-y_1\rrvert =
\eta$. First
assume $\eta\leq L^{\beta}$. The following estimates are valid for $L$
large. First,
\begin{eqnarray*}
\phi(x,z) &=& \mathop{\sum_{y\in\partial V_L\dvtx }}_{\llvert y-y_1\rrvert \leq
L^\alpha}\pi^{(p)}_L(x,y){
\hat{\pi}^{(q)}}_{\psi}(y,z) + \mathop{\sum
_{y\in\partial
V_L\dvtx }}_{\llvert y-y_1\rrvert >
L^\alpha}\pi^{(p)}_L(x,y){\hat{
\pi}^{(q)}}_{\psi}(y,z)
\\
&=&I_1 + I_2.
\end{eqnarray*}
For $I_2$, notice that $\llvert y-y_1\rrvert > L^\alpha$ implies
$\llvert y-x \rrvert >
L^{\alpha}/2$. Using Lemmata~\ref{app-kernelest}(i) and
\ref{hittingprob}(iii) in the first and
Lemma~\ref{hittingprob-technical} in the second inequality, we have
%
%
\begin{equation}
\label{est-phi-eq-1-1} I_2 \leq C\eta L^{-d}
\mathop{\sum_{y\in\partial V_L\dvtx }}_{\llvert y-y_1\rrvert >
L^\alpha}\frac{1}{\llvert x-y\rrvert ^d} \leq C\eta L^{-(d+\alpha)}
\leq
L^{-(d+1/4)}.
\end{equation}
For $I_1$, we first use Lemma~\ref{app-kernelest}(iii) to
deduce
\[
{\hat{\pi}^{(q)}}_{\psi}(y,z) \leq{\hat{\pi}^{(q)}}_{\psi
}(y_1,z)
+ CL^{-(d+1-\alpha)}\log L. %
\]
Therefore by part (vii),
\[
I_1 \leq{\hat{\pi}^{(q)}}_{\psi}(y_1,z)
+ L^{-(d+1/4)} \leq{\hat{\pi}^{\mathrm{B}(q)}} _{\psi
}(y_1,z)
+ 2L^{-(d+1/4)}. %
\]
A similar argument as in~(\ref{est-phi-eq-1-1}),
using Lemma~\ref{exitdensitybm}(i), yields
\[
\int_{y\in\partial C_L\dvtx  \llvert y-y_1\rrvert > L^\alpha}\pi^{\mathrm
{B}(p)}_L(x,dy) \leq
L^{-1/4}. %
\]
Using Lemma~\ref{app-kernelest}(ii) in the first and (v) in the
second inequality, we obtain
\begin{eqnarray*}
{\hat{\pi}^{\mathrm{B}(q)}}_{\psi}(y_1,z) &\leq&{\hat{\pi
}^{\mathrm{B}(q)}}_{\psi}(y_1,z) \int_{y\in
\partial C_L\dvtx
\llvert y-y_1\rrvert \leq L^\alpha}
\pi^{\mathrm{B}(p)}_L(x,dy) + CL^{-(d+1/4)}
\\
&\leq&\int_{y\in\partial C_L\dvtx
\llvert y-y_1\rrvert \leq L^\alpha}\pi^{\mathrm{B}(p)}_L(x,dy){
\hat{\pi}^{\mathrm{B}(q)}}_{\psi
}(y,z) + CL^{-(d+1/4)}
\\
&\leq&{\phi}^{\mathrm{B}}(x,z) + CL^{-(d+1/4)}.
\end{eqnarray*}
Together with~(\ref{est-phi-eq-1-1}), we have shown that $\phi(x,z)
\leq{\phi}^{\mathrm{B}}(x,z) + CL^{-(d+1/4)}$ when $\eta\leq
L^\beta$.

Now we look at the case $\eta> L^\beta$. We take a cube $U_1$ of radius
$L^\alpha$, centered at~$y_1$, and set $W_1 = U_1\cap\partial V_L$.
Then we can find a partition of $\partial V_L\setminus W_1$ into
disjoint sets $W_i = U_i\cap\partial V_L$, $i=2,\ldots,k_L$, where $U_i$
is a cube such that for some $c_1,c_2> 0$ depending only on $d$,
\[
c_1L^{\alpha(d-1)}\leq\llvert W_i\rrvert\leq
c_2L^{\alpha(d-1)}. %
\]
For $i \geq2$, we fix an arbitrary $y_i\in W_i$. Let $W_i^\beta=
\{y\in\mathbb{R}^d\dvtx  \dist(y,W_i) \leq L^\beta\}$. Applying first
Lemma~\ref{app-kernelest}(iii) and then
Lemma~\ref{app-kernelest-technical}(i) gives
%
%
\begin{eqnarray}
\label{est-phi-eq-1-2} \phi(x,z) &\leq&\sum_{i=1}^{k_L}
\pi^{(p)}_L(x,W_i){\hat{\pi
}^{(q)}}_{\psi}(y_i,z)+ L^{-(d+1/4)}
\nonumber\\[-8pt]\\[-8pt]\nonumber
&\leq&\sum_{i=1}^{k_L}
\pi^{\mathrm{B}(p)}_L\bigl(x,W_i^\beta\bigr){
\hat{\pi}^{(q)}}_{\psi}(y_i,z)
\bigl(1+L^{-1/4} \bigr) +L^{-(d+1/4)}.
\end{eqnarray}
As the $W_i^\beta$ overlap, we refine them as follows: set $\tilde{W}_1
= W_1^\beta\cap\partial C_L$, and split $\partial C_L\setminus
\tilde{W}_1$ into a collection of disjoint measurable sets $\tilde{W}_i
\subset W_i^\beta\cap\partial C_L$, $i=2,\ldots,k_L$, such that
$\bigcup_{i=1}^{k_L} \tilde{W}_i = \partial C_L$ and
$\llvert (W_i^\beta\cap\partial C_L)\setminus\tilde{W_i}\rrvert \leq
C_1L^{\alpha(d-2)+\beta}$ for some $C_1=C_1(d)$. By construction we can
find constants $c_3, c_4 > 0$ such that $\llvert \tilde{W}_i\rrvert
\geq
c_3L^{\alpha(d-1)}$ and, for $i=2,\ldots,k_L$,
\[
\inf_{y\in W_i^\beta}\llvert x-y\rrvert\geq c_4\sup
_{y\in\tilde
{W}_i}\llvert x-y\rrvert, %
\]
which implies by Lemma~\ref{exitdensitybm}(i) that
\[
\sup_{y\in W_i^\beta}\pi^{\mathrm{B}(p)}_L(x,y) \leq C\inf
_{y\in
\tilde{W}_i}\pi^{\mathrm{B}(p)}_L(x,y). %
\]
For $i=1,\ldots, k_L$ we then have
\[
\pi^{\mathrm{B}(p)}_L\bigl(x,W_i^\beta\bigr)
\leq\pi^{\mathrm{B}(p)}_L(x,\tilde{W}_i)
\bigl(1+CL^{\beta-\alpha
} \bigr)\leq\pi^{\mathrm{B}(p)}_L(x,
\tilde{W}_i) \bigl(1+L^{-1/4} \bigr). %
\]
Plugging the last line into~(\ref{est-phi-eq-1-2}),
\[
\phi(x,z) \leq\sum_{i=1}^{k_L}
\pi^{\mathrm{B}(p)}_L(x,\tilde{W}_i){\hat{
\pi}^{(q)}}_{\psi
}(y_i,z) \bigl(1+L^{-1/4}
\bigr)+ L^{-(d+1/4)}. %
\]
A reapplication of Lemma~\ref{app-kernelest}(iii), (vii) and then
(ii) yields
\begin{eqnarray*}
\phi(x,z) &\leq&\sum_{i=1}^{k_L}\int
_{\tilde{W}_i}\pi^{\mathrm{B}(p)} _L(x,dy){\hat{
\pi}^{(q)}}_{\psi}(y,z)+ L^{-(d+1/4)}
\nonumber
\\
&\leq&\sum_{i=1}^{k_L}\int
_{\tilde{W}_i}\pi^{\mathrm
{B}(p)}_L(x,dy){\hat{
\pi}^{\mathrm{B}(q)}}_{\psi
}(y,z) \bigl(1+L^{-1/4} \bigr)
+L^{-(d+1/4)}
\\
&=& {\phi}^{\mathrm{B}}(x,z) + CL^{-(d+1/4)}.
\end{eqnarray*}
The reverse inequalities in both the cases $\eta\leq L^\beta$ and
$\eta>
L^{\beta}$ are obtained similarly.

(ii) Let $\psi= (m_y)\in\mathcal{M}_L$ and $z\in\mathbb{Z}^d$. For
$y\in\mathcal{U}_L$ we set
%
%
\begin{equation}
\label{app-phidens} g(y,z) = \frac{1}{m_y}\varphi\biggl(\frac{\llvert
z-y \rrvert }{m_y}
\biggr)\pi^{\mathrm{B}(q)} _{C_{\llvert z-y \rrvert }}(0,z-y).
\end{equation}
Then
\[
{\phi}^{\mathrm{B}}(x,z) = \int_{\partial C_L}\pi^{\mathrm
{B}(p)}_L(x,dy)g(y,z).
\]
Choose a cutoff function $\chi\in C^{\infty} (\mathbb
{R}^d )$
with compact support in $\{x\in\mathbb{R}^d\dvtx  1/2 < \llvert x \rrvert
< 2\}$ such that
$\chi\equiv1$ on $\{x\in\mathbb{R}^d\dvtx  2/3\leq\llvert x \rrvert \leq
3/2\}$. Setting
$m_v = 1$ for $v\notin\mathcal{U}_L$, we define
\[
\tilde{g}(y,z) = g(Ly,z)\chi(y),\qquad y,z\in\mathbb{R}^d. %
\]
By Brownian scaling,
\[
\tilde{g}(y,z) = \frac{1}{m_{Ly}}\varphi\biggl(\frac{\llvert
z-Ly\rrvert }{m_{Ly}} \biggr)
\frac
{1}{\llvert z-Ly\rrvert ^{d-1}} \pi^{\mathrm{B}(q)}_{C_1} \biggl(0,
\frac
{z-y}{\llvert z-y \rrvert } \biggr)\chi(y). %
\]
Notice that $\tilde{g}(\cdot,z)\in C^4 (\mathbb{R}^d )$, with
$\tilde{g}(y,z) = 0$ if $\llvert y \rrvert \notin(1/2,2 )$ or
$\llvert z-Ly\rrvert
\notin
(L/10,10L)$. Let $\mathcal{L} = p(e_1)\partial^2/\partial x_1^2
+\cdots+p(e_d)\partial^2/\partial x_d^2$. Then $u(\overline{x},z) =
{\phi}^{\mathrm{B}}(x,z)$, $x = L\overline{x}$, solves
\[
\label{app-DP} \cases{ \mathcal{L}u(\cdot,z)=0, &\quad in $C_1$,
\vspace*{2pt}\cr
u(\cdot,z)= \tilde{g}(\cdot, z), &\quad on $\partial C_1$.}
\]
By Corollary~6.5.4 of Krylov~\cite{KRY}, $u(\cdot,z)$ is smooth on
$\overline{C}_1$.
Write
\[
{\bigl\llvert u(\cdot,z)\bigr\rrvert}_k = \sum
_{i=0}^{k}\bigl\llVert D^iu(\cdot,z)
\bigr\rrVert_{C_1}. %
\]
By Theorem~6.3.2 in the same book, there exists $C>0$ independent of $z$
such that
\[
{\bigl\llvert u(\cdot,z)\bigr\rrvert}_3 \leq C{\bigl\llvert
\tilde{g}(\cdot,z)\bigr\rrvert}_4. %
\]
A direct calculation shows that $\sup_{z\in\mathbb{R}^d}{\llvert \tilde
{g}(\cdot,z)\rrvert }_4 \leq
CL^{-d}$. Now the claim follows from
\[
\bigl\llVert D^i{\phi}^{\mathrm{B}}(\cdot,z)\bigr\rrVert
_{C_L} = L^{-i}\bigl\llVert D^iu(\cdot,z)\bigr
\rrVert_{C_1}. %
\]

(iii) Let $x, x'\in V_L\cup\partial V_L$. Choose $\tilde{x}\in V_L$ next
to $x$ so that $\llvert \tilde{x}-x\rrvert = 1$ if $x\in\partial
{V}_L$ and
$\tilde{x}=x$ otherwise. Similarly, choose $\tilde{x}'\in V_L$
next to $x'$. By the triangle inequality,
%
%
\begin{eqnarray}
\label{app-est-phi-eq1} && \bigl\llvert\phi(x,z) - \phi\bigl(x',z
\bigr)\bigr\rrvert
\nonumber\\[-8pt]\\[-8pt]\nonumber
&&\qquad \leq\bigl\llvert\phi(x,z) - \phi(\tilde{x},z)\bigr\rrvert+\bigl
\llvert\phi(
\tilde{x},z) - \phi\bigl(\tilde{x}',z\bigr)\bigr\rrvert
+\bigl
\llvert\phi\bigl(\tilde{x}',z\bigr) - \phi\bigl(x',z
\bigr)\bigr\rrvert.\hspace*{-30pt}
\end{eqnarray}
By parts (i) and (ii) combined with the mean value theorem, we get for the
middle term
\begin{eqnarray*}
&& \bigl\llvert\phi(\tilde{x},z) - \phi\bigl(\tilde{x}',z
\bigr)\bigr\rrvert
\\
&&\qquad \leq\bigl\llvert\phi(\tilde{x},z) - {\phi}^{\mathrm{B}}(\tilde
{x},z)\bigr
\rrvert+ \bigl\llvert{\phi}^{\mathrm{B}}(\tilde{x},z) - {\phi
}^{\mathrm
{B}}
\bigl(\tilde{x}',z\bigr)\bigr\rrvert+ \bigl\llvert{
\phi}^{\mathrm{B}}\bigl(\tilde{x}',z\bigr) - \phi\bigl(\tilde
{x}',z\bigr)\bigr\rrvert
\\
&&\qquad \leq C \bigl(L^{-(d+1/4)} + \bigl\llvert x-x' \bigr\rrvert
L^{-(d+1)} \bigr).
\end{eqnarray*}
If $x\in\partial{V_L}$, then $\phi(x,z) = {\hat{\pi}^{(q)}}_{\psi
}(x,z)$, so
that the first term of~(\ref{app-est-phi-eq1}) can be written as
\[
\bigl\llvert\phi(x,z) - \phi(\tilde{x},z)\bigr\rrvert= \biggl\llvert
\sum
_{y\in\partial V_L}\pi^{(p)}_L(
\tilde{x},y) \bigl({\hat{\pi}^{(q)}}_{\psi
}(y,z) - {\hat{
\pi}^{(q)}}_{\psi}(x,z) \bigr)\biggr\rrvert. %
\]
Set $A = \{y\in\partial V_L\dvtx  \llvert x-y\rrvert > L^{1/4}\}$. Then by
Lemmata~\ref{hittingprob}(iii) and~\ref{hittingprob-technical},
\[
\pi^{(p)}_L(\tilde{x},A) \leq C\sum
_{y\in A}\frac{1}{\llvert x-y\rrvert ^d} \leq CL^{-1/4}. %
\]
For all $y\in\partial V_L$, we have by Lemma~\ref{app-kernelest}(i)
$\llvert {\hat{\pi}^{(q)}}_{\psi}(y,z)-{\hat{\pi}^{(q)}}_{\psi
}(x,z)\rrvert
\leq CL^{-d}$. If $y\in
\partial
V_L\setminus A$, then part (iii) gives
$\llvert {\hat{\pi}^{(q)}}_{\psi}(y,z)-{\hat{\pi}^{(q)}}_{\psi
}(x,z)\rrvert
\leq CL^{-(d+3/4)}\log L$. Altogether,
\[
\bigl\llvert\phi(x,z) - \phi(\tilde{x},z)\bigr\rrvert\leq CL^{-(d+1/4)}.
\]
The third term of~(\ref{app-est-phi-eq1}) is treated in exactly the
same way.
\end{pf*}

\subsection{Proofs of Lemmata~\texorpdfstring{\protect\ref{lemmalawler}}{4.1}
and~\texorpdfstring{\protect\ref{hittingprob}}{4.2}}
We let
\[
V_L^{p}(y) = \bigl\{x\in\mathbb{Z}^d\dvtx
\bigl\llvert\Lambda_p^{-1/2}(x-y)\bigr\rrvert\leq L\bigr\}
\]
and
$V_L^p= V_L^p(0)$. Note that
%
%
\begin{equation}
\label{ellipses-ball} V_{(1+2d\kappa)^{-1/2}L}^p\subset V_L
\subset V_{(1-2d\kappa)^{-1/2}L}^p.
\end{equation}
We will make use of the following:
%
\begin{lemma} \label{lemmalawler2}
\textup{(i)}
Let $0 < \ell< L$ and $x \in V_L^p$ with $\ell< \llvert \Lambda
_p^{-1/2}x\rrvert
< L$. Then
\[
\Prw_{x,p} (\tau_{V_L^p} < T_{V_\ell^p} ) =
\frac{\ell
^{-d+2}-\llvert \Lambda_p^{-1/2}x\rrvert ^{-d+2} + O(\ell
^{-d+1})}{\ell^{-d+2}-L^{-d+2}}. %
\]

\textup{(ii)} There exists $C>0$ such that for all $\theta\in\mathbb{R}^d$
with $\llvert \theta \rrvert =1$ and $\ell\geq0$,
\[
\Prw_{0,p}(X_n\cdot\theta\geq-\ell\mbox{ for all }0\leq n
\leq\tau_L)\leq C (\ell+1)L^{-1}. %
\]
\end{lemma}
\begin{pf}(i) For the case of simple random walk, that is, $p=p_o$,
this is
Proposition~1.5.10 of~\cite{LAW}. For the case of general $p\in
\mathcal{P}^{\mathrm{s}}_\kappa$,
one can use Proposition 6.3.1 of~\cite{LawLim} for the Green's function,
and then the proof is exactly the same.

(ii) This is a version of the gambler's ruin estimate; see, for
example,~\cite{LawLim}, Exercise 7.5.
\end{pf}
We turn to the proof of Lemma~\ref{lemmalawler}. We constantly
write $\Prw_{x}$ for $\Prw_{x,p}$ and $\pi_L$ for $\pi^{(p)}_L$,
but stress
that $V_L^p$ is not the same set as $V_L$.

\begin{pf*}{Proof of Lemma~\ref{lemmalawler}}
(i) $\pi_L(\cdot,z)$ is
$p$-harmonic inside $V_L$; that is, for \mbox{$x\in V_L$},
\[
\pi_L(x,z) = \sum_{e\in\mathbb{Z}^d\dvtx  \llvert e \rrvert =1}p(e)
\pi_L(x+e,z). %
\]
Applying a discrete Harnack inequality, as, for example, provided by
Theorem~6.3.9 in the book of
Lawler and Limic~\cite{LawLim}, we obtain
$C^{-1}\pi_L(0,z)\leq\pi_L(\cdot,z)\leq C\pi_L(0,z)$ on $V_{\eta
L}$ for
some $C=C(d,\eta)$, and it remains to show that $\pi_L(0,z)$ has the
right order of magnitude. Note that we cannot directly apply Lemma~6.3.7
in this book, since we look at the exit distribution from $V_L$, not
from $V_L^p$. However, by a last-exit decomposition as in Lemma~6.3.7,
with $g_{V_L}(x,y)=\sum_{k=0}^\infty(1_{V_L}p)^k(x,y)$ and
$\tilde{\tau}_A = \inf\{n\geq1\dvtx  X_n\notin A\}$,
\[
\pi_L(0,z) = \sum_{y\in
V_{L/2}}g_{V_L}(0,y)
\Prw_{z} (X_{\tilde{\tau}_{V_L\setminus
V_{L/2}}}=y ). %
\]
Using~(\ref{ellipses-ball}), we have for $y\in V_{L/2}$, with
$L_1=(1+2d\kappa)^{-1/2}L$ and
$L_2= (1-2d\kappa)^{-1/2}L$,
\[
g_{V^p_{L_1}}(0,y)\leq g_{V_L}(0,y)\leq g_{V^p_{L_2}}(0,y).
\]
For $y\in V_{L/2}\cap\partial(V_L\setminus V_{L/2})$ both outer Green's
functions are of order $L^{-d+2}$, by Proposition~6.3.5
of~\cite{LawLim}. Now by translation invariance and
Lemma~\ref{lemmalawler2}(ii), with $\theta= -z/\llvert z\rrvert $,
\[
\Prw_{z} (X_{\tilde{\tau}_{V_L\setminus
V_{L/2}}}\in V_{L/2} ) \leq
\Prw_{0} (X_n\cdot\theta\geq0\mbox{ for all }0\leq n\leq
\tau_{L/2} )\leq CL^{-1}, %
\]
which proves that $\pi_L(0,z)\leq CL^{-d+1}$. For the lower bound, if
$\kappa$ is small enough, we find an ellipsoid $V^p_{L_3}(y)$ with
$L_3\geq
(9/10)L$, centered at some $y\in V_L$ and lying completely inside $V_L$
such that $z\in\partial V^p_{L_3}(y)\cap\partial V_L$. Also,
$V^p_{L_3/5}(y)\subset V_{L/2}$ if $\kappa$ is small. Therefore,
\[
\Prw_{z} (X_{\tilde{\tau}_{V_L\setminus
V_{L/2}}}\in V_{L/2} ) \geq
\Prw_{z} \bigl(X_{\tilde{\tau
}_{V_{L_3}^p(y)\setminus
V^p_{L_3/5}(y)}}\in V^p_{L_3/5}(y)
\bigr). %
\]
With positive probability, the random walk starting at $z$ enters
$V^p_{L_3}(y)$ in the next step and then visits a point $w$ with
$\dL(w)\geq1$, staying inside $V^p_{L_3}(y)$. Thus
\[
\Prw_{z} \bigl(X_{\tilde{\tau}_{V_{L_3}^p(y)\setminus
V^p_{L_3/5}(y)}}\in V^p_{L_3/5}(y)
\bigr) \geq c\Prw_{w} (T_{V^p_{L_3/5}(y)}<\tau_{V_{L_3}^p(y)} ) \geq
cL^{-1}, %
\]
where the last inequality follows from bounding the expression obtained in
Lemma~\ref{lemmalawler2}(i).
With the estimate on $g_{V_L}$, this proves that $\pi_L(0,z)$ is
bounded from below by $C^{-1}L^{-d+1}$, and (i) follows.

(ii) By the triangle inequality,
\[
\bigl\llvert\pi_L(x,z)-\pi_L\bigl(x',z
\bigr)\bigr\rrvert\leq C\bigl\llvert x-x' \bigr\rrvert\max
_{u,v\in
V_{\eta
L}\dvtx  \llvert u-v\rrvert \leq1}\bigl\llvert\pi_L(u,z)-
\pi_L(v,z)\bigr\rrvert. %
\]
For $u\in V_{\eta L}$, the function $\pi_L(u+\cdot,z)$ is $p$-harmonic
inside $V_{(1-\eta)L}$. The claim now follows from (6.19) of Theorem~6.3.8 in~\cite{LawLim},
together with (i).
\end{pf*}
Before we start with the proof of Lemma~\ref{hittingprob}, we prove a
further auxiliary lemma, which already includes the upper bound of part
(iii).
%
\begin{lemma} Let $x \in V_L$, $y \in\partial V_L$, and set $t =
\llvert x-y\rrvert $.
\label{lemma1}
\begin{longlist}[(ii)]
\item
\[
\Prw_x (X_{\tau_L} = y ) \leq C {\dL(x)}^{-d+1}.
\]
\item
\[
\Prw_x(X_{\tau_L} = y) \leq C\frac{\max\{1,\dL(x)\}}{t}\max
_{x'\in\partial
V_{t/3}(y)\cap V_L}\Prw_{x'}(X_{\tau_L} = y). %
\]
\item
\[
\Prw_x(X_{\tau_L} = y) \leq C\frac{\max\{1,\dL(x)\}}{\llvert
x-y\rrvert ^d}.
\]
\end{longlist}
\end{lemma}
\begin{pf}(i) We can assume that $s = \dL(x) \geq6$. If $s' = \lfloor
s/3\rfloor$,
then\break $\partial V_{s'}(x) \subset V_{L-s'}$. Using
Lemma~\ref{lemmalawler2}(ii), we compute for any $y' \in V_L$ with
$\llvert y-y' \rrvert = 1$, $\theta= -y'/\llvert y'\rrvert $,
\[
\Prw_{y'} (T_{\partial V_{s'}(x)} < \tau_L ) \leq
\Prw_0 (X_n\cdot\theta\geq-1\mbox{ for all }0\leq n\leq
\tau_{V_{s'}} )\leq Cs^{-1}. %
\]
By Lemma~\ref{hittingprob}(i) it follows that uniformly in
$z\in\partial{V}_{s'}(x)$,
\[
\Prw_z (T_x < \tau_L ) \leq
\Prw_z (T_x < \infty) \leq C \bigl(s'
\bigr)^{-d+2}\leq C s^{-d+2}. %
\]
Thus, by the strong Markov property at $T_{\partial V_{s'}(x)}$,
\[
\Prw_{y'} (T_x < \tau_L ) \leq
Cs^{-d+1}. %
\]
Since by time reversibility of symmetric random walk
\begin{eqnarray*}
\Prw_x (X_{\tau_L} = y ) &=& \mathop{\sum
_{y'\in V_L,}}_{\llvert y'-y \rrvert =
1}\Prw_x \bigl(X_{\tau_L} = y,
X_{\tau_L-1} = y' \bigr)
\\
&=& \frac{1}{2d}\mathop{\sum_{y'\in V_L,}}_{\llvert y'-y \rrvert = 1}
\Prw_{y'} (T_x < \tau_L ),
\end{eqnarray*}
the claim is proved.

(ii) We may assume that $t=\llvert x-y\rrvert >100d$ and $\dL(x) <
t/100$. Choose a
point $x'\in\mathbb{Z}^d$ such that $V_{t/10}(x')\cap V_L =\varnothing
$ and
$\llvert x-x' \rrvert \leq\dL(x)+t/10 + \sqrt{d}$. Then $\llvert
x-x' \rrvert \leq
t/5$. Furthermore, since $\llvert x'-y\rrvert \geq4t/5$,
\[
\bigl(V_{t/4}\bigl(x'\bigr)\cup\partial
V_{t/4}\bigl(x'\bigr) \bigr) \cap V_{t/3}(y) =
\varnothing. %
\]
We apply twice the strong Markov property and obtain
\[
\Prw_x(X_{\tau_L} = y) \leq\Prw_x (
\tau_{V_{t/4}(x')} < \tau_L ) \max_{z\in\partial V_{t/3}(y)\cap V_L}
\Prw_z(X_{\tau_L}=y). %
\]
Arguing much as in (i), Lemma~\ref{lemmalawler2}(ii) shows
\[
\Prw_x (\tau_{V_{t/4}(x')} < \tau_L )\leq C
\frac{\max\{
1,\dL(x)\}}{t}, %
\]
which completes the proof of part (ii).

(iii) By (ii) it suffices to prove that for some constant $K$ and for all
$\ell\geq1$,
%
%
\begin{equation}
\label{eq1} \max_{z\in\partial V_{\ell/3}(y)\cap V_L}\Prw_z(X_{\tau_L}=y)
\leq K\ell^{-d+1}.
\end{equation}
Let $c_1$ and $c_2$ be the constants from (i) and (ii), respectively.
Define $\eta= 3^{-d}c_2^{-1}$ and $K = \max\{3^{d(d-1)}c_2^{d-1},
c_1\eta^{-d+1}\}$. For $\ell\leq3^{d}c_2$, there is nothing to prove since
$K\ell^{-d+1} \geq1$. Thus let $\ell> 3^{d}c_2$, and choose $\ell
_0$ with $\ell_0
< \ell\leq2\ell_0$. Assume that~(\ref{eq1}) is proved for all $\ell
' \leq
\ell_0$. We show that~(\ref{eq1}) also holds for $\ell$. For $z$ with
$\dL(z) \geq\eta\ell$, it follows from (i) that
\[
\Prw_z(X_{\tau_L}=y) \leq c_1\eta^{-d+1}
\ell^{-d+1} \leq K\ell^{-d+1}. %
\]
If $1 \leq\dL(z) < \eta\ell$, then by (ii) and the fact that $\ell
/3 \leq
\ell_0$
\begin{eqnarray*}
\Prw_z(X_{\tau_L}=y) &\leq& c_2
\frac{\max\{1,\dL(z)\}}{\llvert z-y \rrvert }\max_{z'\in\partial
V_{t/9}(y)\cap
V_L}\Prw_z(X_{\tau_L}=y)
\nonumber
\\
&\leq& c_23\eta K (\ell/3 )^{-d+1} \leq K
\ell^{-d+1}.
\end{eqnarray*}
If $\dL(z) < 1$, then again by (i)
\[
\Prw_z(X_{\tau_L}=y) \leq c_23
\ell^{-1} K (\ell/3 )^{-d+1} \leq K\ell^{-d+1}.
\]
This proves the claim.
\end{pf}

\begin{pf*}{Proof of Lemma~\ref{hittingprob}}(i) follows from
Proposition~6.4.2 of~\cite{LawLim}.

(ii) We consider different cases. If $\llvert x-y\rrvert \leq\dL
(y)/2$, then $\dL
(x) \geq
\dL(y)/2$, and thus by Lemma~\ref{hittingprob}(i),
\[
\Prw_x (T_{V_a(y)} < \tau_L ) \leq
\Prw_x (T_{V_a(y)}<\infty) \leq C{ \biggl(
\frac{a}{\llvert x-y\rrvert } \biggr)}^{d-2}\leq C\frac{a^{d-2}\dL
(y)\dL
(x)}{\llvert x-y\rrvert ^d}. %
\]
For the rest of the proof we assume that $\llvert x-y\rrvert > \dL
(y)/2$. Set $a' =
\dL(y)/5$. First we argue that for the case $1\leq a\leq a'$, we only have
to prove the bound for~$a'$. Indeed, if $\dL(y)/6 \leq a < a'$, we get an
upper bound by replacing $a$ by $a'$. For $1 \leq a < \dL(y)/6$, the
strong Markov property together with Lemma~\ref{hittingprob}(i)
yields
\begin{eqnarray*}
\Prw_x (T_{V_a(y)} < \tau_L )&\leq&\max
_{z\in\partial(\mathbb{Z}^d\setminus
V_{a'}(y))}\Prw_z (T_{V_a(y)} <
\tau_L )\Prw_x (T_{V_{a'}(y)} <
\tau_L )
\\
&\leq& C{ \biggl(\frac{a}{a'-1} \biggr)}^{d-2} \frac{{(a')}^{d-2}\dL
(y)\max\{1,\dL(x)\}}{\llvert x-y\rrvert ^d}
\\
&\leq& C\frac{a^{d-2}\dL(y)\max\{1,\dL(x)\}}{\llvert x-y\rrvert ^d}.
\end{eqnarray*}
Now we prove the claim for $a = \dL(y)/5$. We take a point $y'\in
\partial
V_L$ closest to $y$. If $\llvert x-z \rrvert \geq\llvert x-y\rrvert
/2$ for all $z\in V_a(y')$, then
by Lemma~\ref{lemma1}(iii),
\[
\max_{z\in V_a(y')}\Prw_x (X_{\tau_L}=z ) \leq C
2^d\frac{\max\{1,\dL(x)\}}{\llvert x-y\rrvert ^d}. %
\]
As a subset of $\mathbb{Z}^d$, $V_{a}(y')\cap\partial V_L$ contains on
the order of $\dL(y)^{d-1}$ points. Therefore, by
Lemma~\ref{lemmalawler}(i), we deduce that there exists some $\delta>
0$ such that
\[
\min_{x'\in V_a(y)}\Prw_{x'} \bigl(X_{\tau_L} \in
V_{a}\bigl(y'\bigr) \bigr) \geq\delta. %
\]
We conclude that
%
%
\begin{eqnarray}
\label{eq2} && \frac{a^{d-1}\max\{1,\dL(x)\}}{\llvert x-y\rrvert ^d}
\nonumber
\\
&&\qquad \geq c\Prw_{x} \bigl(X_{\tau_L} \in V_{a}
\bigl(y'\bigr) \bigr) \geq c\Prw_{x}
\bigl(X_{\tau_L} \in V_{a}\bigl(y'\bigr),
T_{V_a(y)} < \tau_L \bigr)
\nonumber\\[-8pt]\\[-8pt]\nonumber
&&\qquad =c\sum_{x'\in V_a(y)}\Prw_{x}
\bigl(X_{T_{V_a(y)}} = x', T_{V_a(y)} <
\tau_L \bigr)\Prw_{x'} \bigl(X_{\tau_L} \in
V_a\bigl(y'\bigr) \bigr)
\nonumber
\\
&&\qquad \geq c \delta\cdot\Prw_{x} (T_{V_a(y)} <
\tau_L ).\nonumber
\end{eqnarray}
On the other hand, if $\llvert x-z \rrvert < \llvert x-y\rrvert /2$
for some $z\in V_a(y')$, then
\[
\llvert x-y\rrvert\leq\llvert x-z \rrvert+ \bigl\llvert z-y'
\bigr\rrvert+ \bigl\llvert y'-y \bigr\rrvert\leq2\dL(y) + \llvert
x-y\rrvert/2 %
\]
and thus
%
%
\begin{equation}
\label{eq3} \dL(y)/2 < \llvert x-y\rrvert\leq4\dL(y).
\end{equation}
If $\dL(x) \geq4\dL(y)/5$, we use again Lemma~\ref{hittingprob}(i). For
$\dL(x) < 4\dL(y)/5$, we get by Lemma~\ref{lemmalawler2}(ii),
\[
\Prw_x (T_{V_a(y)} < \tau_L ) \leq
\Prw_x (T_{V_{L-4\dL(y)/5}}<\tau_L ) \leq C
\frac{\max\{1,\dL(x)\}}{\dL(y)}. %
\]
Together with~(\ref{eq3}), this proves the claim for $a=\dL(y)/5$. It
remains to handle the case $\max\{1,\dL(y)/5\} \leq a$. If $z\in
V_{6a}(y)$, we have
\[
\llvert x-y\rrvert\leq\llvert x-z \rrvert+ 6a %
\]
and thus, using $\llvert x-y\rrvert > 7a$,
\[
\llvert x-y\rrvert\leq7\llvert x-z \rrvert. %
\]
Therefore Lemma~\ref{lemma1}(iii) yields
\[
\max_{z\in V_{6a}(y)}\Prw_x (X_{\tau_L}=z ) \leq C
\frac{\max\{1,\dL(x)\}}{\llvert x-z \rrvert ^d} \leq7^dC\frac{\max\{
1,\dL
(x)\}}{\llvert x-y\rrvert ^d}. %
\]
Again by Lemma~\ref{lemmalawler}(i), we find some $\delta> 0$ such
that
\[
\min_{x'\in V_a(y)}\Prw_{x'} \bigl(X_{\tau_L} \in
V_{6a}(y) \bigr) \geq\delta. %
\]
A similar argument as in~(\ref{eq2}), with $V_a(y')$ replaced by
$V_{6a}(y)$, finishes the proof of (ii).

(iii) It only remains to prove the lower bound. Let $t=\llvert x-z
\rrvert $. First assume
$t\geq L/2$. By replacing $V_L$ and $V_{2L/3}$ by appropriate
ellipsoids as
in the proof of the lower bound of Lemma~\ref{lemmalawler}(i), we deduce
with part~(i) of Lemma~\ref{lemmalawler2} that
\[
\Prw_x (T_{V_{2L/3}}<\tau_L ) \geq c
\frac{\dL(x)}{t}. %
\]
The claim then follows from the strong Markov property at time
$T_{V_{2L/3}}$ and Lemma~\ref{lemmalawler}(i). Now assume $t<L/2$. We
can restrict ourselves to the case $10\sqrt{d}\leq t<L/2$. We find
$t'\in[t,t+\sqrt{d}]$ and $x'\in V_L$ such that $V_{t'}(x')\subset
V_L$ and
$z\in\partial V_{t'}(x')$. If $\dL(x)> t/2$, a simple geometric
consideration and Lemma~\ref{lemmalawler}(i) show that there is a
strictly positive probability to exit the ball $V_{t/2}(x)$ within
$V_{2t/3}(x')$. Since by the same lemma,
%
%
\begin{equation}
\label{eq4} \inf_{y\in V_{2t/3}(x')}\Prw_y (
\tau_L=z ) \geq ct^{-(d-1)},
\end{equation}
we obtain the claim in this case again by applying the strong Markov
property. Finally, assume $\dL(x) \leq t/2$. Once more by
Lemma~\ref{lemmalawler2}(i),
\[
\Prw_x (T_{V_{L-t/3}}<\tau_L ) \geq c
\frac{\dL(x)}{t} %
\]
and
%
%
\begin{eqnarray}
\label{eq5} \Prw_x (\tau_L=z )&\geq&
\Prw_x (\tau_L=z, T_{L-t/3}<
\tau_L, T_{V_{2t/3}(x')}<\tau_L )\nonumber
\\
&\geq& c\frac{\dL(x)}{t}\Prw_x (\tau_L=z \mid
T_{L-t/3}<\tau_L, T_{V_{2t/3}(x')}<\tau_L )
\\
&&{} \times\Prw_x (T_{V_{2t/3}(x')}<\tau_L \mid
T_{L-t/3}<\tau_L ).\nonumber
\end{eqnarray}
By~(\ref{eq4}), we deduce
\begin{eqnarray*}
\Prw_x (\tau_L=z \mid T_{L-t/3}<
\tau_L, T_{V_{2t/3}(x')}<\tau_L ) &=&
\Prw_x (\tau_L=z \mid T_{V_{2t/3}(x')}<
\tau_L )
\\
&\geq& ct^{-(d-1)}.
\end{eqnarray*}
Finally, a geometric argument and Lemma~\ref{lemmalawler}(i) show that
the probability in~(\ref{eq5}) is bounded from below by some
$\delta>0$.
\end{pf*}

\subsection{Proofs of Lemmata~\texorpdfstring{\protect\ref{exitdensitybm}}{4.3}
and~\texorpdfstring{\protect\ref{exitbmpbmq}}{4.4}}

\mbox{}

\begin{pf*}{Proof of Lemma~\ref{exitdensitybm}}
By Brownian scaling, we may restrict ourselves to the case $L=1$.
Let
%
%
\begin{equation}
\label{Ep} E_p =\bigl\{y\in\mathbb{R}^d\dvtx  \bigl
\llvert\Lambda_p^{1/2}y\bigr\rrvert< 1\bigr\},
\end{equation}
and fix $x\in C_1$, $z\in\partial C_1$. Set $z'=\Lambda_p^{-1/2}z$,
$x'=\Lambda_p^{-1/2}x$. Again by Brownian scaling and the fact that
$\det
\Lambda_p^{-1/2} = 1 + O(\kappa)$, the (continuous versions of the)
densities satisfy
\[
\pi^{\mathrm{B}(p)}_{C_1}(x,z)= \pi^{\mathrm
{B}(p_o)}_{E_p}
\bigl(x',z'\bigr) \bigl( 1+O(\kappa) \bigr).
\]
Now (i) follows from Theorem~1 of Krantz~\cite{KRA} (with $\Omega=
E_p$), and for (ii), one can use the derivative
estimates in Section~2 of the same paper.
\end{pf*}

%
%
\begin{figure}

\includegraphics{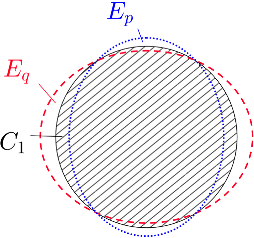}

\caption{The ball $C_1$ (shaded), and the ellipsoids $E_p$ and $E_q$.}\label{fig7}
\end{figure}

\begin{pf*}{Proof of Lemma~\ref{exitbmpbmq}}
We can assume $L=1$. Let $\eta= \llVert q-p\rrVert _1$. Define $E_p$ as
in~(\ref{Ep}), and similarly $E_q$, see Figure~\ref{fig7}. Let $x\in C_{2/3}$, $z\in\partial
C_1$, and put $x'=\Lambda_p^{-1/2}x$, $x''=\Lambda_q^{-1/2}x$ and
$z'=\Lambda_p^{-1/2}z$, $z''=\Lambda_q^{-1/2}z$. If $\kappa$ is small,
both $x'$ and $x''$ lie in $C_{3/4}\subset C_{4/5}\subset E_p\cap E_q$.
For the rest of the proof, we write $\pi^{\mathrm{B}}_{E_p}$ instead of
$\pi^{\mathrm{B}(p_o)}_{E_p}$, and similarly $\pi^{\mathrm
{B}}_{E_q}$ for $\pi^{\mathrm{B}(p_o)}
_{E_q}$. If
$\gamma$ is a parametrization of the unit sphere $\partial C_1$, then
$\Lambda_p^{-1/2}\circ\gamma$ and $\Lambda_q^{-1/2}\circ\gamma$ are
parametrizations of $E_p$ and $E_q$, respectively. Since the coefficients
of the covariance matrices satisfy $(\Lambda_p)_{i,i}=
(\Lambda_q)_{i,i}+O(\eta)$, we obtain by Brownian scaling,
%
%
\begin{equation}
\label{bmpbmq-1}  \quad\qquad\bigl\llvert\pi^{\mathrm{B}(p)}_{C_1}(x,z)-
\pi^{\mathrm{B}(q)}_{C_1}(x,z)\bigr\rrvert\leq C\bigl\llvert
\pi^{\mathrm{B}}_{E_p}\bigl(x',z'\bigr)-
\pi^{\mathrm{B}} _{E_q}\bigl(x'',z''
\bigr)\bigr\rrvert+O(\eta).
\end{equation}
Clearly
%
%
\begin{equation}
\label{z-z} \bigl\llvert z'-z''
\bigr\rrvert= \bigl\llvert\bigl(\Lambda_p^{-1/2}-
\Lambda_q^{-1/2}\bigr)z\bigr\rrvert\leq C\eta,
\end{equation}
and also $\llvert x'-x''\rrvert \leq C\eta$. By the derivative
estimate of
Lemma~\ref{exitdensitybm}(ii),
\[
\bigl\llvert\pi^{\mathrm{B}}_{E_q}\bigl(x',z''
\bigr) - \pi^{\mathrm{B}}_{E_q}\bigl(x'',z''
\bigr)\bigr\rrvert\leq C\eta. %
\]
With~(\ref{bmpbmq-1}), the claim will therefore follow if we show that
\[
\bigl\llvert\pi^{\mathrm{B}}_{E_p}\bigl(x,z'\bigr)-
\pi^{\mathrm{B}}_{E_q}\bigl(x,z''\bigr)
\bigr\rrvert\leq C\eta%
\]
uniformly in $x\in C_{3/4}$, $z'\in\partial E_p$, $z''\in\partial
E_q$ with $\llvert z'-z'' \rrvert \leq
C\eta$. In this direction, recall that the Green's function on
$\mathbb{R}^d$ of standard $d$-dimensional Brownian motion is given by
\[
\Phi(x,y) = \frac{c_d}{\llvert x-y\rrvert ^{d-2}}, %
\]

\noindent where $c_d={\bolds\Gamma}(d/2-1)/(2\pi^{d/2})$; cf.~\cite{LawLim},
page 241. The
Green's function of standard Brownian motion killed outside
$E_p$ is given by (see, e.g., Evans~\cite{EV})
\[
\Phi(x,y) - \Phi_x^{(p)}(y),\qquad x,y\in
E_p, x\neq y, %
\]
where the corrector function $\Phi_x^{(p)}$ solves the Dirichlet problem
($x$ is fixed)
%
%
\begin{equation}
\label{corrector-dirichlet} \cases{ \Delta\Phi^{(p)}_x= 0, &\quad in
$E_p$,
\vspace*{2pt}\cr
\Phi^{(p)}_x= \Phi(x,\cdot), &\quad
on $\partial E_p$.}
\end{equation}
Furthermore, the density
$\pi^{\mathrm{B}}_{E_p}(x,z')$ with respect to surface measure on
$\partial E_p$ is
given by the normal derivative of the Green's function in the direction
of the inward unit normal vector $\nu_p=\nu_p(z')$ on $\partial E_p$,
that is,
\begin{eqnarray*}
\pi^{\mathrm{B}}_{E_p}\bigl(x,z'\bigr) &=&
\partial_{\nu_p} \bigl(\Phi\bigl(x,z'\bigr)-
\Phi^{(p)}_x\bigl(z'\bigr) \bigr)
\\
&=& \nabla_{z'} \bigl(\Phi\bigl(x,z'\bigr)-
\Phi^{(p)}_x\bigl(z'\bigr) \bigr)\cdot
\nu_p\bigl(z'\bigr),\qquad z'\in\partial
E_p.
\end{eqnarray*}
With $\nu_q$ denoting the inward unit normal on $\partial E_q$, we
therefore have to show that
%
%
\begin{equation}
\label{corrector-estimate} \bigl\llvert\partial_{\nu_p} \bigl(\Phi
\bigl(x,z'\bigr)-\Phi^{(p)}_x
\bigl(z'\bigr) \bigr)-\partial_{\nu_q} \bigl(\Phi
\bigl(x,z''\bigr)-\Phi^{(q)}_x
\bigl(z''\bigr) \bigr)\bigr\rrvert\leq C\eta
\end{equation}
uniformly in $x\in C_{3/4}$, $z'\in\partial E_p$, $z''\in\partial
E_q$ with $\llvert z'-z'' \rrvert \leq
C\eta$. First, note that
\begin{eqnarray*}
&& \bigl\llvert\partial_{\nu_p}\Phi\bigl(x,z'
\bigr)- \partial_{\nu_q} \Phi\bigl(x,z''\bigr)
\bigr\rrvert
\\
&&\qquad \leq\bigl\llvert\nabla_{z'}\Phi\bigl(x,z'\bigr)-
\nabla_{z''}\Phi\bigl(x,z''\bigr)\bigr
\rrvert+ \bigl\llvert\nabla_{z'} \Phi\bigl(x,z'\bigr)
\cdot\bigl(\nu_p\bigl(z'\bigr)-\nu_q
\bigl(z''\bigr)\bigr)\bigr\rrvert.
\end{eqnarray*}
Using~(\ref{z-z}) and
\[
\partial_{z'_i}\Phi\bigl(x,z'\bigr)= c_d(2-d)
\frac{(x-z')_i}{\llvert x-z'\rrvert ^d}, %
\]
we easily obtain
\[
\bigl\llvert\nabla_{z'}\Phi\bigl(x,z'\bigr)-
\nabla_{z''}\Phi\bigl(x,z''\bigr)\bigr
\rrvert\leq C\eta. %
\]
Moreover, with $\Lambda=\Lambda_q^{-1/2}\Lambda_p^{1/2}$, we have
\[
\nu_q\bigl(z''\bigr) =
\frac{\Lambda^{-1}\nu_p(z')}{\llvert \Lambda^{-1}\nu_p(z')\rrvert }. %
\]
The coefficients of the diagonal matrix $\Lambda$ are of order
$1+O(\eta)$,
which shows
%
%
\begin{equation}
\label{nup-nuq} \bigl\llvert\nabla_{z'} \Phi\bigl(x,z'
\bigr)\cdot\bigl(\nu_p\bigl(z'\bigr)-
\nu_q\bigl(z''\bigr)\bigr)\bigr\rrvert
\leq C \bigl\llvert\nu_p\bigl(z'\bigr)-
\nu_q\bigl(z''\bigr)\bigr\rrvert\leq C
\eta.
\end{equation}
In view of~(\ref{corrector-estimate}), it remains to prove that
\[
\bigl\llvert\partial_{\nu_p}\Phi_x^{(p)}
\bigl(z'\bigr)- \partial_{\nu_q} \Phi_x^{(q)}
\bigl(z''\bigr)\bigr\rrvert\leq C\eta. %
\]
First recall that $\Phi_x^{(p)}$ solves the Dirichlet
problem~(\ref{corrector-dirichlet}). The boundary function
$\Phi(x,\cdot)$ is smooth in a tubular neighborhood of $\partial
E_p$, with
bounded derivatives up to arbitrary order. By multiplying with an
appropriate smooth cutoff function equal to $1$ near the boundary, we
obtain a smooth function in $\mathbb{R}^d$. By Corollary~6.5.4 of
Krylov~\cite{KRY}, we see that $\Phi_x^{(p)}$ is a smooth function in
$\overline{E}_p$, and similarly $\Phi_x^{(q)}$ is smooth in
$\overline{E}_q$.

Now, with $\Lambda$ as before, we have equality of the sets
$\Lambda\overline{E}_p = \overline{E}_q$. Fix $x\in C_{3/4}\subset
E_p\cap E_q$,
and let
\[
u(y) = \Phi_x^{(q)}(\Lambda y)-\Phi_x^{(p)}(y),
\qquad y\in\overline{E}_p. %
\]
With $f(\cdot) = \Delta_y\Phi_x^{(q)}(\Lambda\cdot)$ and
$g(\cdot) = \Phi(x,\Lambda\cdot)-\Phi(x,\cdot)$, $u$ solves
\[
\label{dirichlet} \cases{ \Delta u=f, &\quad in $E_p$,
\vspace*{2pt}\cr
u=g, &\quad
on $\partial E_p$.} %
\]
Recalling that the coefficients of $\Lambda$ are of order $1+O(\eta
)$, we
use harmonicity of~$\Phi_x^{(q)}$ and boundedness of the derivatives
to obtain
\[
{\bigl\llVert D^0f\bigr\rrVert}_{E_p}+{\bigl\llVert
D^1f\bigr\rrVert}_{E_p}\leq C\eta. %
\]
The function $g$ is smooth in a tubular neighborhood $U$ of $\partial E_p$,
and a similar (explicit) calculation as above gives
\[
\sum_{i=0}^3{\bigl\llVert
D^ig\bigr\rrVert}_{U}\leq C\eta. %
\]
We extend $g$ to the interior of $E_p$ such that
$\sum_{i=0}^3{\llVert D^ig\rrVert }_{\overline{E}_p}\leq C\eta$.
Then, applying a
Schauder estimate as given by Theorem~6.3.2 of Krylov~\cite{KRY},
we deduce that the derivatives of $u$ up to order $2$ are uniformly bounded
by $C\eta$. But, similarly to above~(\ref{nup-nuq}),
\begin{eqnarray*}
\bigl\llvert\partial_{\nu_p}\Phi_x^{(p)}
\bigl(z'\bigr)- \partial_{\nu_q} \Phi_x^{(q)}
\bigl(z''\bigr)\bigr\rrvert&\leq&\bigl\llvert
\nabla_{z'}\Phi_x^{(p)}\bigl(z'
\bigr)- \nabla_{z''}\Phi_x^{(q)}
\bigl(z''\bigr)\bigr\rrvert+ C\eta
\\
&\leq&\bigl\llvert\nabla_{z'}u\bigl(z'\bigr)\bigr
\rrvert+C\eta
\\
&\leq& C\eta,
\end{eqnarray*}
where in the next-to-last step we used $z''=\Lambda z'$ and
\[
\bigl\llvert\nabla_{z''}\Phi_x^{(q)}
\bigl(z''\bigr) - \nabla_{z'}
\Phi_x^{(q)}\bigl(z''\bigr)
\bigr\rrvert\leq C\eta.
\]\upqed
\end{pf*}

\subsection{Proofs of
Propositions~\texorpdfstring{\protect\ref{super-localclt}}{5.1}
and~\texorpdfstring{\protect\ref{super-behaviorgreen}}{5.2}}
By a small abuse of notation, we will in this part write $
{\hat{\pi}}_m$ for ${\hat{\pi}}_{\psi_m}$.
Since ${\hat{\pi}}_m(x,y) = {\hat{\pi
}}_m(0,y-x)$, it suffices to look at
${\hat{\pi}}_m(x)={\hat{\pi}}_m(0,x)$ and
$\hat{g}_{m,\mathbb{Z}^d}(x) =
\hat{g}_{m,\mathbb{Z}^d}(0,x)$. Recall the definitions
of $\lambda_{m,i}$ and
$\Lambda_m$ from Section~\ref{clt}.

\begin{pf*}{Proof of Proposition~\ref{super-localclt}}
For bounded $m$, that is, $m\leq m_0$ for some $m_0$, the result is a
special case of Theorem~2.1.1 in~\cite{LawLim}. Also, for $n\leq
n_0$ and all $m$, the statement follows from Lemma~\ref
{app-kernelest}(i). We therefore have to prove the proposition only for
large $n$ and $m$.
To this end, let
\[
B_m = [-\sqrt{\lambda_{m,1}} \pi, \sqrt{
\lambda_{m,1}} \pi]\times\cdots\times[-\sqrt{\lambda_{m,d}}
\pi, \sqrt{\lambda_{m,d}} \pi], %
\]
and for $\theta\in B_m$ set
\[
\phi_m(\theta) = \sum_{y\in\mathbb{Z}^d}
\erm^{i\theta\cdot
\Lambda_m^{-1/2}y}{\hat{\pi}}_m(y). %
\]
The Fourier inversion formula gives
\[
{\hat{\pi}}_m^n(x) = \frac{1}{(2\pi)^d\det\Lambda_m^{1/2}} \int
_{B_m}\erm^{-ix\cdot\Lambda_m^{-1/2}\theta}\bigl[\phi_m(\theta)
\bigr]^n\,\dt\theta. %
\]
We decompose the integral into
\[
(2\pi)^d\det\Lambda_m^{1/2}n^{d/2}{
\hat{\pi}}_m^n(x) = I_0(n,m,x)+\cdots+
I_3(n,m,x), %
\]
where, with $\beta= \sqrt{n} \theta$,
\begin{eqnarray*}
I_0(n,m,x) &=& \int_{\mathbb{R}^d}\erm^{-ix\cdot\Lambda_m^{-1/2}\beta
/\sqrt{n}}
\erm^{-\llvert \beta\rrvert ^2/2}\,\dt\beta,
\\
I_1(n,m,x) &=& \int_{\llvert \beta\rrvert \leq
n^{1/4}}\erm^{-ix\cdot\Lambda_m^{-1/2}\beta/\sqrt{n}}
\bigl(\bigl[\phi_m(\beta/\sqrt{n})\bigr]^n-
\erm^{-\llvert \beta\rrvert ^2/2} \bigr)\,\dt\beta,
\\
I_2(n,m,x) &=& -\int_{\llvert \beta\rrvert >n^{1/4}}\erm^{-ix\cdot
\Lambda
_m^{-1/2}\beta/\sqrt{n}}
\erm^{-\llvert \beta\rrvert ^2/2}\,\dt\beta,
\\
I_3(n,m,x) &=& n^{d/2}\int_{n^{-1/4}<\llvert \theta \rrvert, \theta\in
B_m}
\erm^{-ix\cdot\Lambda_m^{-1/2}\theta}\bigl[\phi_m(\theta)\bigr]^n\,\dt
\theta.
\end{eqnarray*}
By completing the square in the exponential, we get
\[
I_0(n,m,x)= (2\pi)^{d/2}\exp\biggl(-\frac{\mJ^
2(x)}{2n}
\biggr). %
\]
For $I_1$ and $\llvert \beta\rrvert \leq n^{1/4}$, we expand $\phi
_m$ in a series around
the origin
%
%
\begin{eqnarray}
\label{app-lclt1} \phi_m(\beta/\sqrt{n}) &=& 1 - \llvert\beta\rrvert
^2/2n + \llvert\beta\rrvert^4O \bigl(n^{-2}
\bigr),
\nonumber\\[-8pt]\\[-8pt]\nonumber
\log\phi_m(\beta/\sqrt{n}) &=& - \llvert\beta\rrvert
^2/2n + \llvert\beta\rrvert^4O \bigl(n^{-2}
\bigr).
\end{eqnarray}
Therefore,
\[
{\bigl[\phi_m(\beta/\sqrt{n})\bigr]}^n =
\erm^{-\llvert \beta\rrvert ^2/2} \bigl(1+\llvert\beta\rrvert^4O
\bigl(n^{-1} \bigr) \bigr), %
\]
so that
\[
\bigl\llvert I_1(n,m,x)\bigr\rrvert\leq O \bigl(n^{-1}
\bigr)\int_{\llvert \beta\rrvert \leq
n^{1/4}}\erm^{-\llvert \beta\rrvert ^2/2}\llvert\beta\rrvert
^4\,\dt\beta= O \bigl(n^{-1} \bigr). %
\]
Similarly, $I_2$ is bounded by
\[
\bigl\llvert I_2(n,m,x)\bigr\rrvert\leq C \int_{n^{1/4}}^\infty
r^{d-1}\erm^{-r^2/2}\,\dt r = O \bigl(n^{-1} \bigr).
\]
Concerning $I_3$, we follow closely the proof of~\cite{BZ}, Proposition~\textup{B1}, and split the integral further into
\begin{eqnarray*}
n^{-d/2}I_3(n,m,x) &=& \int_{n^{-1/4}<\llvert \theta \rrvert \leq a}
+ \int
_{a<\llvert \theta \rrvert \leq A} + \int_{A<\llvert \theta
\rrvert \leq m^{\alpha}} + \int
_{m^{\alpha}<\llvert \theta \rrvert, \theta\in B_m}
\\
&=& (I_{3,0} + I_{3,1} + I_{3,2} +
I_{3,3} ) (n,m,x),
\end{eqnarray*}
where $0<a<A$ and $\alpha\in(0,1)$ are constants that will be chosen
in a
moment, independently of $n$ and $m$. By~(\ref{app-lclt1}), we can find
$a> 0$ such that for $\llvert \beta\rrvert \leq a\sqrt{n}$, $\log\phi
_m(\theta)
\leq
-\llvert \theta \rrvert ^2/3$ (recall that $\beta= \sqrt{n} \theta
$). Then
\[
\bigl\llvert I_{3,0}(n,m,x)\bigr\rrvert\leq C\int
_{n^{-1/4}}^\infty r^{d-1}\erm^{-nr^2/3}\,\dt
r = O \bigl(n^{-(d+2)/2} \bigr). %
\]
As a consequence of Lemma~\ref{lemmalawler}(i) and of our coarse
graining, it follows that for any $0<a < A$, one has for some $0 < \rho=
\rho(a,A) < 1$, uniformly in $m$ (and $p\in{\mathcal{P}^{\mathrm
{s}}_\kappa}$),
\[
\sup_{a\leq\llvert \theta \rrvert \leq A}\bigl\llvert\phi_m(\theta
)\bigr
\rrvert\leq\rho. %
\]
Using this fact,
\[
\bigl\llvert I_{3,1}(n,m,x)\bigr\rrvert\leq CA^{d}
\rho^n = O \bigl(n^{-(d+2)/2} \bigr). %
\]
To deal with the last two integrals is more delicate
since we have to take into account the $m$-dependency. First,
\[
\bigl\llvert I_{3,2}(n,m,x)\bigr\rrvert\leq\int_{A<\llvert \theta
\rrvert \leq
m^\alpha}
\bigl\llvert\phi_m(\theta)\bigr\rrvert^n\,\dt\theta.
\]
We bound the integrand pointwise. Since ${\hat{\pi
}}_m(\cdot)$ is invariant under
the maps $e_i\mapsto-e_i$, $e_j\mapsto e_j$ for $j\neq i$, it suffices to
look at $\theta$ with all components positive. Assume $\theta_1 =
\max\{\theta_1,\ldots,\theta_d \}$. Set $M=\lfloor
2\pi\sqrt{\lambda_{m,1}}/\theta_1\rfloor$ and $K = \lfloor
5m/M\rfloor$. Notice that ${\hat{\pi}}_m(x) > 0$ implies
$\llvert x \rrvert < 2m$. By taking $A$
large enough, we can assume that on the domain of integration, $M \leq
m$. First,
\begin{eqnarray*}
\phi_m(\theta)&=&\sum_{ (x_2,\ldots,x_d )} \Biggl(\exp
\Biggl(\frac{i}{\sqrt{\lambda_{m,1}}}\sum_{s=2}^dx_s
\theta_s \Biggr)
\\
&& \hspace*{39pt}{} \times \sum_{j=1}^K\sum_{x_1=-2m
+ (j-1)M}^{-2m +
jM-1}\exp\biggl(
\frac{ix_1\theta_1}{\sqrt{\lambda_{m,1}}} \biggr){\hat{\pi}}_m(x) \Biggr).
\end{eqnarray*}
Inside the $x_1$-summation, we write for each $j$ separately
\[
{\hat{\pi}}_m(x) = {\hat{\pi}}_m(x)-{\hat{
\pi}}_m\bigl(x^{(j)}\bigr)+{\hat{\pi}}_m
\bigl(x^{(j)}\bigr), %
\]
where $x^{(j)} = (-2m + (j-1)M,x_2,\ldots,x_d)$. By
Corollary~\ref{app-kernelest-cor},
\[
\bigl\llvert{\hat{\pi}}_m(x)-{\hat{\pi}}_m
\bigl(x^{(j)}\bigr)\bigr\rrvert\leq C\biggl\llvert\frac{x_1+2m
-(j-1)M}{m}
\biggr\rrvert^{1/2}m^{-d}. %
\]
Thus
\[
\Biggl\llvert\sum_{x_1=-2m+ (j-1)M}^{-2m
+jM-1}\exp\biggl(
\frac{ix_1\theta_1}{\sqrt{\lambda_{m,1}}} \biggr) \bigl({\hat{\pi
}}_m(x)-{\hat{\pi
}}_m\bigl(x^{(j)}\bigr) \bigr)\Biggr\rrvert\leq C
\theta_1^{-3/2}m^{-d+1} %
\]
and
\[
\Biggl\llvert\sum_{j=1}^K\sum
_{x_1=-2m
+ (j-1)M}^{-2m +
jM-1}\exp\biggl(\frac{ix_1\theta_1}{\sqrt{\lambda_m}} \biggr)
\bigl({\hat{\pi}}_m(x)-{\hat{\pi}}_m
\bigl(x^{(j)}\bigr) \bigr)\Biggr\rrvert\leq C\theta_1^{-1/2}m^{-d+1}.
\]
On our domain of integration, $0< (\theta_1/\sqrt{\lambda
_{m,1}} ) \leq
Cm^{\alpha-1} < 2\pi$ for large $m$. Therefore,
\begin{eqnarray*}
&&\Biggl\llvert\sum_{j=1}^K{
\hat{\pi}}_m\bigl(x^{(j)}\bigr)\sum
_{x_1 = -2m +(j-1)M}^{-2m +
jM-1}\exp\biggl(\frac{ix_1\theta_1}{\sqrt{\lambda_{m,1}}} \biggr)
\Biggr\rrvert\\
&&\qquad
\leq CKm^{-d}\biggl\llvert\frac{1-\exp(i\theta_1M/\sqrt{\lambda
_{m,1}} )}{1-\exp(i\theta_1/\sqrt{\lambda_{m,1}}
)}\biggr\rrvert
\\
&&\qquad\leq C
\llvert\theta\rrvert m^{-d},
\end{eqnarray*}
and altogether for sufficiently large $A$, $m$ and $n$,
\[
\int_{A<\llvert \theta \rrvert \leq
m^\alpha}\bigl\llvert\phi_m(\theta)\bigr
\rrvert^n\,\dt\theta\leq C_1^n\int
_{A<\llvert \theta \rrvert \leq m^\alpha} \biggl(\frac{1}{\sqrt
{\llvert \theta \rrvert }} +\frac{\llvert \theta \rrvert }{m}
\biggr)^n\,\dt\theta= O \bigl(n^{-(d+2)/2} \bigr). %
\]
For $I_{3,3}$ we again assume all components of $\theta$ positive and as
before
$\theta_1 = \max\{\theta_1,\ldots,\theta_d\}$. Since
\[
{\hat{\pi}}_m(x) = \sum_{y=-2m}^{x_1}
\bigl( {\hat{\pi}}_m(y,x_2,\ldots,x_d)-{
\hat{\pi}}_m(y-1,x_2,\ldots,x_d) \bigr),
\]
we have
\begin{eqnarray*}
\bigl\llvert\phi_m(\theta)\bigr\rrvert&\leq& Cm^{d-1}
\Biggl\llvert\sum_{x_1=-2m}^{2m}\exp\biggl(
\frac{ix_1\theta_1}{\sqrt
{\lambda_{m,1}}} \biggr)
\\
&&\hspace*{32pt}{}\times \sum_{y=-2m}^{x_1}
\bigl({\hat{\pi}}_m(y,x_2,\ldots,x_d)-{\hat{
\pi}}_m(y-1,x_2,\ldots,x_d) \bigr)\Biggr
\rrvert
\\
&\leq& Cm^{d-1}\sum_{y=-2m}^{2m}
\bigl\llvert{\hat{\pi}}_m(y,x_2,\ldots,x_d)-{
\hat{\pi}}_m(y-1,x_2,\ldots,x_d)\bigr\rrvert
\\
&&{} \times\Biggl\llvert\sum_{x_1=y}^{2m}
\exp\biggl(\frac{ix_1\theta
_1}{\sqrt{\lambda_{m,1}}} \biggr)\Biggr\rrvert.
\end{eqnarray*}
The sum over the exponentials is estimated by $C m/\llvert \theta
\rrvert $, so that again
with Corollary~\ref{app-kernelest-cor},
\[
\bigl\llvert\phi_m(\theta)\bigr\rrvert\leq C_2
m^{1/2}\llvert\theta\rrvert^{-1}. %
\]
Hence, for $\alpha$ close to $1$ and large $n$, $m$,
\[
\int_{m^{\alpha}<\llvert \theta \rrvert, \theta\in
B_m}\bigl\llvert\phi_m^n(
\theta)\bigr\rrvert\,\dt\theta\leq C_2^nm^{n/2+\alpha
(d-n)}=
O \bigl(n^{-(d+2)/2} \bigr). %
\]\upqed
\end{pf*}
For Proposition~\ref{super-behaviorgreen}, we need additionally a
large deviation estimate.
%
\begin{lemma}[(Large deviation estimate)]
\label{app-ldestim}
Let $p\in\mathcal{P}^{\mathrm{s}}_\kappa$. There exist constants
$c_1,c_2 > 0$ such that for
$\llvert x \rrvert \geq3m$,
\[
{\hat{\pi}}_m^n(x) \leq c_1m^{-d}
\exp\biggl(-\frac
{\llvert x \rrvert ^2}{c_2nm^2} \biggr). %
\]
\end{lemma}
\begin{pf}
Write $\Prw$ for $\Prw_{0,{\hat{\pi}}_m}$ and $\Erw$ for
the expectation with respect to $\Prw$, and denote by $X_n^j$ the $j$th
component of the random walk $X_n$ under $\Prw$. For $r > 0$, since $p$
is symmetric,
\begin{eqnarray*}
\sum_{y\dvtx  \llvert y \rrvert \geq r}{\hat{\pi}}_m^n(y)
&\leq&\sum_{j=1}^d\Prw\bigl(\bigl\llvert
X_n^j\bigr\rrvert\geq d^{-1/2}r\bigr)
\nonumber
\\
&\leq&2d\max_{j=1,\ldots,d}\Prw\bigl(X_n^j
\geq d^{-1/2}r\bigr).
\end{eqnarray*}
We claim that
\[
\Prw\bigl(X_n^j \geq d^{-1/2}r\bigr) \leq\exp
\biggl(-\frac{r^2}{8dnm^2} \biggr). %
\]
By the martingale maximal inequality for all $t, \lambda> 0$,
\[
\Prw\bigl(X_n^j \geq\lambda\bigr) \leq
\erm^{-t\lambda}\Erw\bigl[\exp\bigl(tX_n^j\bigr)
\bigr] = \erm^{-t\lambda} \bigl(\Erw\bigl[\exp\bigl(tX_1^j
\bigr) \bigr] \bigr)^n. %
\]
Since $X_1^j \in(-2m,2m)$ and $x\mapsto\erm^{tx}$ is convex, it
follows that
\[
\exp\bigl(tX_1^j\bigr)\leq\frac{1}{2}
\frac{(2m-X_1^j)}{2m}\erm^{-2tm} + \frac{1}{2}\frac{(2m+X_1^j)}{2m}
\erm^{2tm}. %
\]
Therefore, using again symmetry of $X_1^j$,
\[
\Erw\bigl[\exp\bigl(tX_n^j \bigr) \bigr] \leq
\bigl(\tfrac{1}{2}\erm^{-2tm}+\tfrac{1}{2}
\erm^{2tm} \bigr)^n = \cosh^n(2tm)\leq
\erm^{2nt^2m^2} %
\]
and
\[
\Prw\bigl(X_n^j \geq d^{-1/2}r\bigr)\leq
\erm^{-td^{-1/2}r}\erm^{2nt^2m^2}. %
\]
Putting $t = r/(4\sqrt{d} nm^2)$ we get
\[
\Prw\bigl(X_n^j \geq d^{-1/2}r\bigr)\leq\exp
\biggl({-\frac
{r^2}{8dnm^2}} \biggr). %
\]
From this it follows that
\begin{eqnarray*}
{\hat{\pi}}_m^n(x) &=& \sum
_{y\dvtx \llvert y \rrvert \geq\llvert x \rrvert -2m} {\hat{\pi
}}_m^{n-1}(y){\hat{
\pi}}_m(x-y) \leq\frac{c_1}{m^d}\exp\biggl({-
\frac{(\llvert x \rrvert -2m)^2}{8d(n-1)m^2}} \biggr)
\\
&\leq&\frac{c_1}{m^d}\exp\biggl({-\frac{\llvert x \rrvert
^2}{c_2nm^2}} \biggr).
\end{eqnarray*}\upqed
\end{pf}
\begin{pf*}{Proof of Proposition~\ref{super-behaviorgreen}}(i) This
follows from Proposition~\ref{super-localclt}.

(ii) We set
\[
N = N(x,m) = \frac{\llvert x \rrvert ^2}{m^2} \biggl(\log\frac
{\llvert x \rrvert ^2}{m^2}
\biggr)^{-2}. %
\]
We split $\hat{g}_{m,\mathbb{Z}^d}(x)$ into
\[
\hat{g}_{m,\mathbb{Z}^d}(x) = \sum_{n=1}^\infty
{\hat{\pi}}_m^n(x) = \sum_{n=1}^{\lfloor N\rfloor}{
\hat{\pi}}_m^n(x)+ \sum_{n=\lfloor N\rfloor+1}^{\infty}{
\hat{\pi}}_m^n(x). %
\]
For the first sum on the right, we use the large deviation estimate from
Lemma~\ref{app-ldestim},
\[
\sum_{n=1}^{\lfloor N\rfloor}{\hat{
\pi}}_m^n(x) \leq c_1m^{-d}\sum
_{n=1}^{\lfloor
N\rfloor}\exp\biggl(-\frac{\llvert x \rrvert ^2}{c_2 nm^2}
\biggr) = O \bigl(\llvert x \rrvert^{-d} \bigr). %
\]
In the second sum, we replace the transition probabilities by the
expressions obtained in Proposition~\ref{super-localclt}.
The error terms are estimated by
\[
\sum_{n=\lfloor N\rfloor+1}^\infty O \bigl(m^{-d}n^{-(d+2)/2}
\bigr) = O \biggl(\llvert x \rrvert^{-d} \biggl(\log
\frac{\llvert x \rrvert ^2}{m^2} \biggr)^d \biggr). %
\]
Putting $t_n = 2n\mJ^{-2}(x)$, we obtain for the main part
\begin{eqnarray*}
&& \sum_{n=\lfloor
N\rfloor+
1}^\infty
\frac{1}{(2\pi n)^{d/2}\det\Lambda_m^{1/2}}\exp\biggl(-\frac{\mJ
^2(x)}{2n} \biggr)
\\
&&\qquad = \frac{\mJ^{-d+2}(x)}{2\pi^{d/2}\det\Lambda_m^{1/2}}\sum_{n=\lfloor
N\rfloor+1}^\infty
t_n^{-d/2}\exp(-1/t_n) (t_n-t_{n-1})
\nonumber
\\
&&\qquad = \frac{\mJ^{-d+2}(x)}{2\pi^{d/2}\det\Lambda_m^{1/2}}\int_0^\infty
t^{-d/2}\exp(-1/t)\,\dt t + O \bigl(\llvert x \rrvert^{-d}
\bigr).
\end{eqnarray*}
This proves the statement for $\llvert x \rrvert \geq3m$ with
\[
c(d) = \frac{1}{2\pi^{d/2}}\int_0^\infty
t^{-d/2}\exp(-1/t)\,\dt t. %
\]\upqed
\end{pf*}
\end{appendix}

\section*{Acknowledgments}
We are indebted to Ofer Zeitouni for sharing his insights and for
many helpful discussions, which had a considerable influence on this paper.

Moreover, we would like to thank an anonymous referee for his thorough
reading of our manuscript and for many valuable comments, which greatly
improved the presentation.


%

\printaddresses
\end{document}